\documentclass[12pt,pdftex]{article}
\pdfoutput=1

\usepackage{amsbsy,amsfonts,amsmath,amssymb,amsthm,enumerate,framed}
\usepackage{graphicx}

\theoremstyle{plain}
\allowdisplaybreaks[1]
\numberwithin{equation}{section}

\newtheorem{thm}{Theorem}[section]

\newtheorem{lmm}[thm]{Lemma}
\newtheorem{prp}[thm]{Proposition}

\usepackage{bbm}
\newcommand{\indic}{\mathbbm{1}}
\newcommand{\ind}[1]{\indic\raisebox{-1pt}{$\sst\{#1\}$}}
\newcommand{\indn}[1]{\indic\raisebox{-2pt}{$\sst #1$}}

\newcommand{\bL}[2]{\mathbf{L}_{#1}^{#2}}
\newcommand{\bM}[4]{\sideset{_{#1}^{#2}}{_{#3}^{#4}}{\mathop{\mathbf{M}}}}
\newcommand{\bR}[3]{\sideset{_{#1}^{#2}}{_{#3}}{\mathop{\mathbf{R}}}}

\newcommand{\cC}{\mathcal{C}}

\newcommand{\cG}{\mathcal{G}}

\newcommand{\cL}{\mathcal{L}}

\newcommand{\conn}{\longleftrightarrow}

\newcommand{\dc}{d_\mathrm{c}}
\newcommand{\db}{\Longleftrightarrow}
\newcommand{\dpst}{\displaystyle}
\newcommand{\even}{^{\sss\mathrm{even}}}

\newcommand{\lbeq}[1]{\label{eq:#1}}
\newcommand{\Ld}{\mL^{\!d}}

\newcommand{\mE}{\mathbb{E}}
\newcommand{\mL}{\mathbb{L}}

\newcommand{\mN}{\mathbb{N}}
\newcommand{\mP}{\mathbb{P}}

\newcommand{\mT}{\mathbb{T}}
\newcommand{\mZ}{\mathbb{Z}}
\newcommand{\nn}{\nonumber}
\newcommand{\odd}{^{\sss\mathrm{odd}}}
\newcommand{\pc}{p_\mathrm{c}}

\newcommand{\piv}{{\tt piv}}
\newcommand{\Proof}[1]{\medskip\noindent\textit{#1}~}
\newcommand{\QED}{\hspace*{\fill}\rule{7pt}{7pt}\medskip}

\newcommand{\ssc}[1]{^{\sss(#1)}}
\newcommand{\sss}{\scriptscriptstyle}
\newcommand{\sst}{\scriptstyle}

\newcommand{\Td}{\mT^d}

\newcommand{\vep}{\varepsilon}
\newcommand{\vno}{\varnothing}

\newcommand{\Zd}{\mZ^d}

\newcommand{\Zp}{\mZ_+}

\newcommand{\pisawone}[1]{~\mathop{\raisebox{-0.8pc}{\includegraphics[scale
 =0.2]{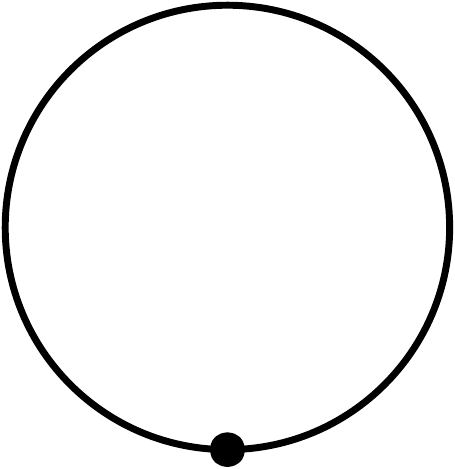}}}_{#1}~}
\newcommand{\pisawtwo}[2]{~\mathop{\raisebox{-0.9pc}{\includegraphics[scale
 =0.2]{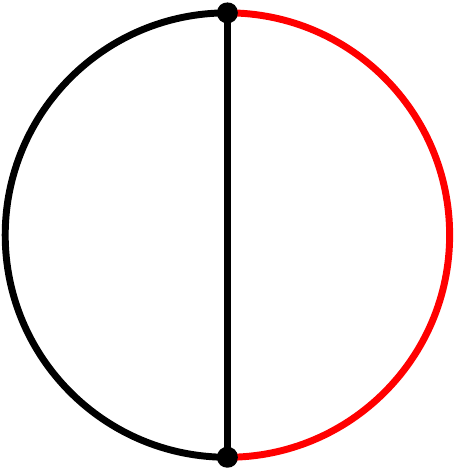}}}_{#1}^{#2}~}
\newcommand{\pisawthree}[2]{\sideset{_{#1}}{_{#2}}{\mathop{\raisebox{-1pc}
 {\includegraphics[scale=0.2]{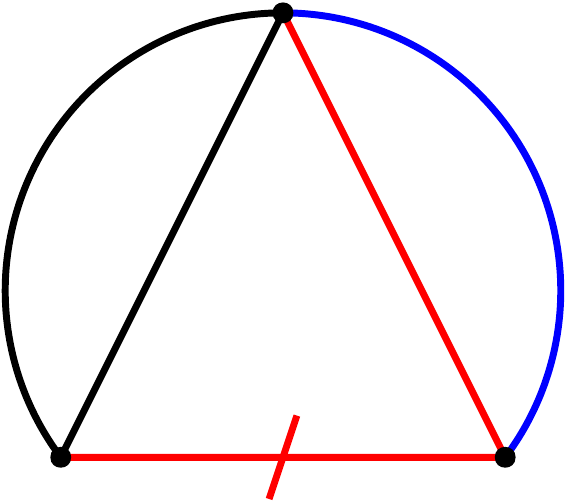}}}}}
\newcommand{\pisawfour}[3]{\sideset{_{#1}}{^{#2}}{\mathop{\raisebox{-1pc}
 {\includegraphics[scale=0.2]{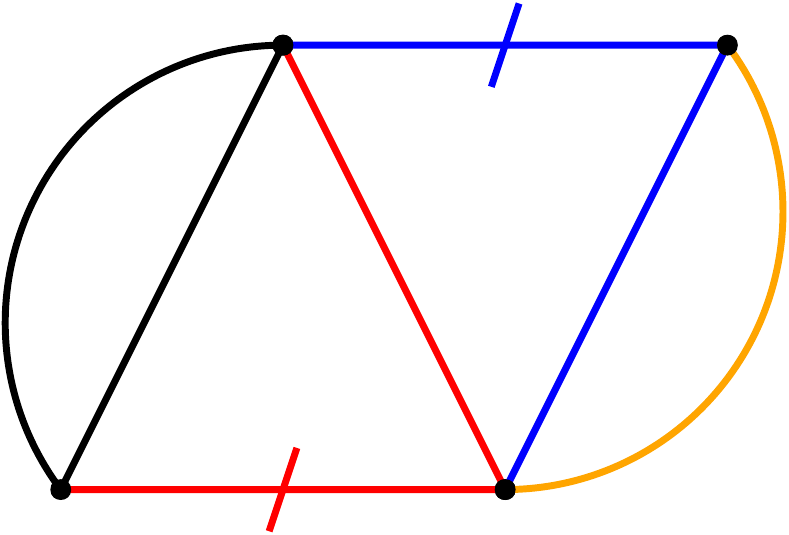}}}}_{\raisebox{5pt}{$\sst #3$}}}
\newcommand{\pisawfourw}[3]{\sideset{_{#1}}{^{#2}}{\mathop{\raisebox{-1pc}
 {\includegraphics[scale=0.2]{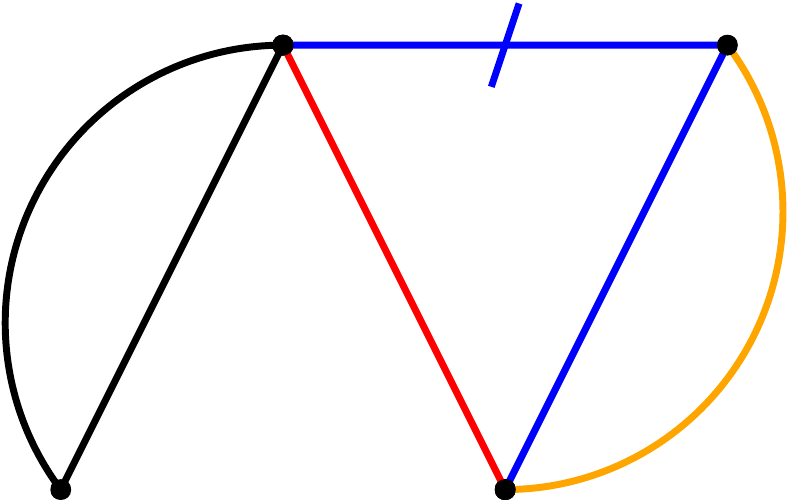}}}}_{#3}}
\newcommand{\pisawfourww}[3]{\sideset{_{#1}}{^{#2}}{\mathop{\raisebox{-1pc}
 {\includegraphics[scale=0.2]{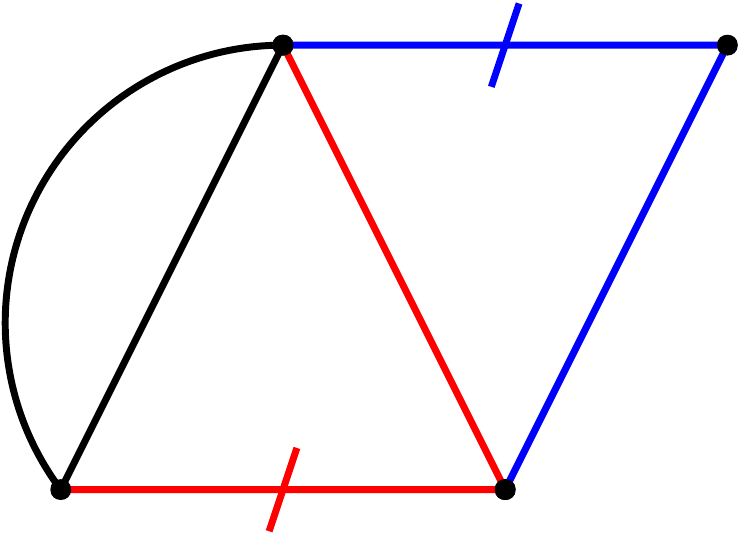}}}}_{\raisebox{5pt}{$\sst #3$}}}
\newcommand{\pisawfive}[3]{\sideset{_{#1}}{_{#2}}{\mathop{\raisebox{-1pc}
 {\includegraphics[scale=0.2]{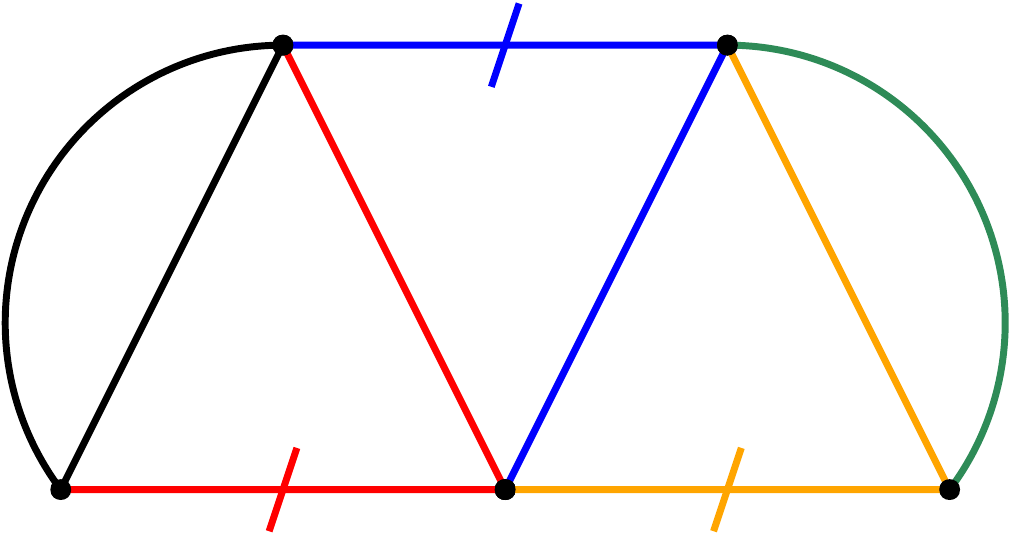}}}}_{\raisebox{5pt}{$\sst #3$}}}
\newcommand{\pisawfivew}[3]{\sideset{_{#1}}{_{#2}}{\mathop{\raisebox{-1pc}
 {\includegraphics[scale=0.2]{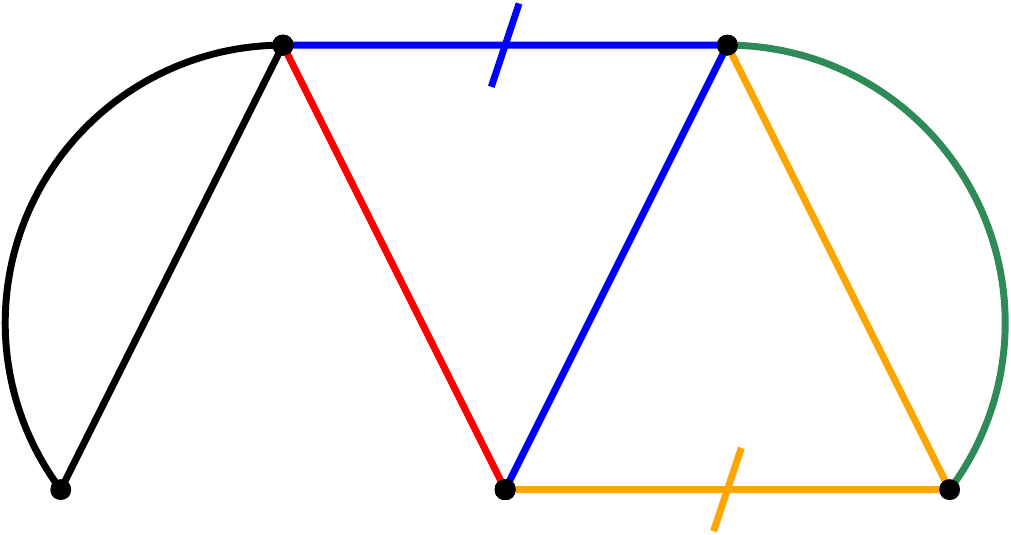}}}}_{\raisebox{5pt}{$\sst #3$}}}
\newcommand{\pisawfiveww}[3]{\sideset{_{#1}}{_{#2}}{\mathop{\raisebox{-1pc}
 {\includegraphics[scale=0.2]{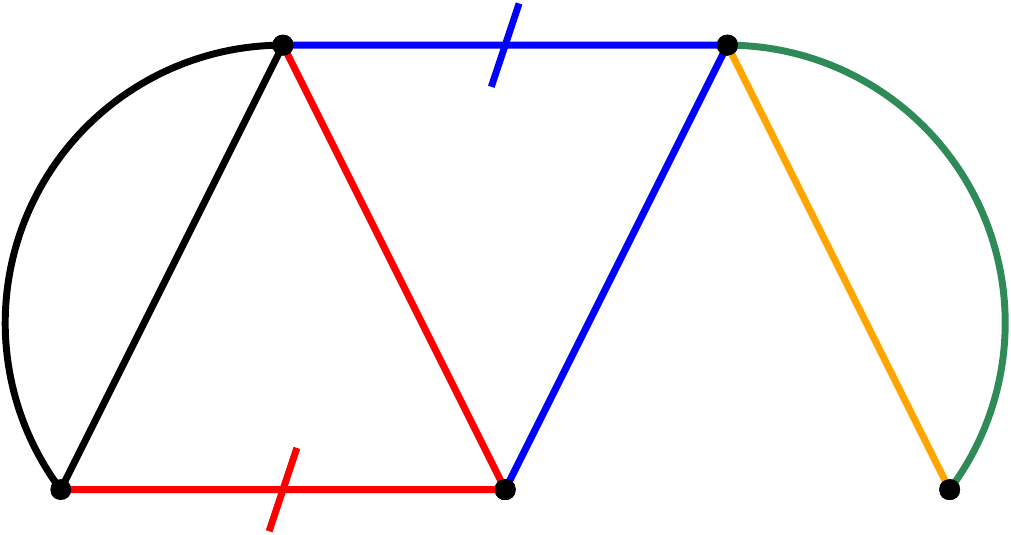}}}}_{\raisebox{5pt}{$\sst #3$}}}

\usepackage[usenames]{color}

\usepackage{tikz-lace-perc}
\usetikzlibrary{external}
\tikzsetexternalprefix{tikz-img/}
\tikzsetfigurename{lace-}
\usepackage{mathtools}
\renewcommand{\refeq}[1]{(\ref{eq:#1})}
\usepackage{extarrows}
\usepackage{centernot}
\usepackage[subrefformat=parens]{subcaption}
\graphicspath{{./figs/}}

\newcommand{\LaceCoefficient}[1]{\mathop{\pi_p^{\scriptscriptstyle (#1)}}\nolimits}
\newcommand{\LaceCoefficientFT}[1]{\mathop{\hat{\pi}_p^{\scriptscriptstyle (#1)}}\nolimits}
\newcommand{\SumOfLaceCoefficientsFT}[2][p]{%
  \ifthenelse{\equal{#2}{}}{
    \mathop{\hat{\Pi}_{#1}}\nolimits
  }{
    \mathop{\hat{\Pi}_{#1}^{\scriptscriptstyle\mathrm{#2}}}\nolimits
  }
}
\newcommand{\RemainderTerm}[1]{\mathop{R_p^{\scriptscriptstyle (#1)}}\nolimits}
\newcommand{\KroneckerDelta}[2]{\delta_{#1, #2}}
\newcommand{\Bond}[2]{\{#1, #2\}}
\newcommand{\DirectedBond}[2]{(#1, #2)}
\newcommand{\ClusterWithout}[3][]{%
  \ifthenelse{\equal{#1}{}}{
    \tilde{\mathcal{C}}^{#2}
  }{
    \tilde{\mathcal{C}}_{\scriptscriptstyle #1}^{#2}
  }%
  (#3)%
}
\newcommand{\PivotalBonds}[2]{\mathord{\mathtt{piv}}(#1, #2)}
\newcommand{\stcomp}[1]{#1^\mathrm{c}}
\newcommand{\event}[1]{\{#1\}}
\newcommand{\Event}[1]{\left\{#1\right\}}
\newcommand{\Connected}{\mathrel{\longleftrightarrow}}
\newcommand{\ConnectedThrough}[1]{\mathrel{\xlongleftrightarrow{#1}}}
\newcommand{\DoublyConnected}{\mathrel{\Longleftrightarrow}}

\newcommand{\DoublyConnectedThrough}[1]{\mathrel{\xLongleftrightarrow{#1}}}
\newcommand{\Probability}[1][]{%
  \mathop{%
    \ifthenelse{\equal{#1}{}}{
      \mathbb{P}_p
    }{
      \mathbb{P}_p^{\scriptscriptstyle #1}
    }%
  }\nolimits%
}
\newcommand{\Expectation}[1][]{%
  \mathop{%
    \ifthenelse{\equal{#1}{}}{
      \mathbb{E}_p
    }{
      \mathbb{E}_p^{\scriptscriptstyle #1}
    }%
  }\nolimits%
}
\newcommand{\DisjointPath}{\circ}
\newcommand{\supNorm}[1]{\lVert#1\rVert_{\infty}}
\newcommand{\SupNorm}[1]{\left\lVert#1\right\rVert_{\infty}}
\newcommand{\abs}[1]{\lvert#1\rvert}

\newcommand{\Laplacian}[1]{\hat{\Delta}_{#1}}

\theoremstyle{definition}

\newtheorem{definition}[thm]{Definition}

\title{A survey on the lace expansion
for the nearest-neighbor models on the BCC lattice}

\author{
Satoshi Handa\thanks{\tt handa@math.sci.hokudai.ac.jp}, \quad
Yoshinori Kamijima\thanks{\tt ykami@eis.hokudai.ac.jp}, \quad
Akira~Sakai\thanks{\tt sakai@math.sci.hokudai.ac.jp}\\
Department of Mathematics\\ Hokkaido University}

\begin{document}
\maketitle

\begin{abstract}
The aim of this survey is to explain, in a self-contained and relatively
beginner-friendly manner, the lace expansion for the nearest-neighbor models
of self-avoiding walk and percolation that converges in all dimensions above
6 and 9, respectively.  To achieve this, we consider a $d$-dimensional version
of the body-centered cubic (BCC) lattice, on which it is extremely easy to
enumerate various random-walk quantities.  Also, we choose a particular set
of bootstrapping functions, by which a notoriously complicated part of the
lace-expansion analysis becomes rather transparent.
\end{abstract}

\tableofcontents

\section{Introduction}\label{s:intro}
The lace expansion is one of the few mathematically rigorous methods to prove
critical behavior for various statistical-mechanical models in high dimensions.
It can show that the two-point function for the concerned model, up to the
critical point, is bounded by the Green function for the underlying random
walk in high dimensions.  During the course of learning this method, it also
provides good exercises in various mathematical skills from graph theory and
algebraic identities to Fourier analysis and probability theory.

First, we explain background, some historical facts and the purposes of
this survey.

\subsection{Background}\label{ss:background}
Corporation of infinitely many particles results in various intriguing and
challenging problems.  One of those is to understand phase transitions
and critical behavior of statistical-mechanical models, such as
percolation and the ferromagnetic Ising model.  For percolation, for example,
it exhibits a phase transition when the bond-occupation parameter $p$ crosses
its critical value $\pc$.  If $p$ is far below $\pc$, each cluster of
occupied vertices is so small that we may use standard probabilistic
techniques for i.i.d.~random variables to predict what happens in the
subcritical phase.  If $p$ is far above $\pc$, on the other hand,
vacant vertices can only form tiny islands and most of the other vertices
are connected to form a single gigantic cluster.  
However, when $p$ is close to $\pc$, the cluster of connected vertices
from the origin may be extremely large but porous in a nontrivial way,
and therefore naive perturbation methods fail.

A similar phenomenon occurs for self-avoiding walk (SAW), a century-old
statistical-mechanical model for linear polymers.  Consider a locally finite,
amenable and transitive graph as space.  A standard example is the
$d$-dimensional integer lattice $\Zd$.  The main observable to be investigated
is the SAW two-point function, which is the following generating function with
fugacity $p\ge0$:
\begin{align}\lbeq{Gsaw-intro}
G_p(x)=\sum_{\omega:o\to x}p^{|\omega|}\prod_{j=1}^{|\omega|}D(\omega_j
 -\omega_{j-1})\prod_{0\le s<t\le|\omega|}(1-\lambda\delta_{\omega_s,
 \omega_t}),
\end{align}
where the sum is over the nearest-neighbor paths $\omega$ on the concerned
lattice from the origin $o$ to $x$, $|\omega|$ is the number of steps along
$\omega$, and $D$ is the 1-step distribution of simple random walk (RW):
$D(x)=(2d)^{-1}\delta_{|x|,1}$ on $\Zd$.  The parameter $\lambda\in[0,1]$
is the intensity of self-avoidance; the model with $\lambda=1$ is
called strictly SAW, while the one with $\lambda\in(0,1)$ is called weakly
SAW.  The two-point function with $\lambda=0$ is equivalent to the RW Green
function $S_p(x)\equiv\sum_{n=0}^\infty p^nD^{*n}(x)$, where $D^{*n}$ is the
$n$-fold convolution of $D$.  The critical point (= the radius of convergence)
for RW is $p=1$.  For SAW, because of subadditivity, there is a critical point
$\pc\ge1$ such that the susceptibility $\chi_p\equiv\sum_xG_p(x)$ is finite if
and only if $p<\pc$ and diverges as $p\uparrow\pc$ (see, e.g., \cite{ms93}).

The way $\chi_p$ diverges is intriguing, as it shows power-law behavior
as $(\pc-p)^{-\gamma}$ with the critical exponent $\gamma$.  It is
considered to be universal in the sense that the value of $\gamma$ depends
only on $d$ and is insensitive to $\lambda\in(0,1]$ and the detail lattice
structure.  For example, the value of $\gamma$ for strictly SAW on $\mZ^2$
is believed to be $\frac{43}{32}$ and equal to that for weakly SAW on the
2-dimensional triangular lattice.  This is not the case for the critical
point $\pc$, as its value may vary depending on $\lambda\in(0,1]$ and the
detail lattice structure.  Other statistical-mechanical models that exhibit
divergence of the susceptibility are also characterized by the critical
exponent $\gamma$, and many physicists as well as mathematicians have been
trying hard to identify the value of $\gamma$ and classify the models into
different universality classes since last century.

\subsection{The mean-field theory}\label{ss:MFtheory}
Because of the nonlocal self-avoidance constraint
$\prod_{0\le s<t\le|\omega|}(1-\lambda\delta_{\omega_s,\omega_t})$ in
\refeq{Gsaw-intro}, SAW  does not enjoy the Markovian property, which holds
only when $\lambda=0$.  If there is a way to average out the self-avoidance
effect and absorb it into the fugacity $p$, then $G_p(x)$ may be approximated
by the RW Green function $S_\mu(x)$ with a mean-field fugacity
$\mu=\mu(\Zd,\lambda,p)$, and therefore $\chi_p$ may be approximated by $\sum_xS_\mu(x)=(1-\mu)^{-1}$.  Presumably, $\mu(\pc)=1$.  If $\mu$ is
left-differentiable at $\pc$, then this implies $\chi_p\asymp(\pc-p)^{-1}$
(i.e., $\chi_p$ is bounded above and below by positive multiples of
$(\pc-p)^{-1}$) as $p\uparrow\pc$.  In this respect, the mean-field value
for the critical exponent $\gamma$ is 1.

However, realizing the above idea is highly nontrivial.  As a first step,
one may want to use perturbation theory from the mean-field model
(i.e., $\lambda=0$).  The expansion of  the self-avoidance constraint in
powers of $\lambda>0$ yields
\begin{align}\lbeq{naive-exp}
\prod_{0\le s<t\le|\omega|}(1-\lambda\delta_{\omega_s,\omega_t})
 =\sum_{\Gamma\in\cG[0,|\omega|]}(-\lambda)^{|\Gamma|}\prod_{\{s,t\}\in
 \Gamma}\delta_{\omega_s,\omega_t},
\end{align}
where $\Gamma$, which is called a graph, is a set of pairs of indices on
$[0,|\omega|]\equiv\{0,1,\dots,|\omega|\}$, $\cG[0,|\omega|]$ is a set of
such graphs, and $|\Gamma|$ is the cardinality of $\Gamma$.  The trivial
contribution from $\Gamma\equiv\vno$ is the unperturbed solution $S_p(x)$,
which is already bad because its radius of convergence is 1, while $\pc\ge1$.
The first correction term proportional to $\lambda$ is
\begin{align}
-\lambda\sum_{\omega:o\to x}p^{|\omega|}\prod_{j=1}^{|\omega|}D(\omega_j
 -\omega_{j-1})\sum_{0\le s<t\le|\omega|}\delta_{\omega_s,\omega_t}=
-\lambda\big(S_p(o)-1\big)S_p^{*2}(x).
\end{align}
The higher-order correction terms are more involved, but the radius of
convergence of each term is always $p=1$.  What is worse, the alternating
series of those terms is absolutely convergent only when $p$ is close to zero,
because the sum over $\Gamma\in\cG[0,|\omega|]$ is potentially huge as long
as $\lambda>0$.  As a result, this naive expansion cannot be applied near
$\pc$ in order to justify the mean-field behavior.

\subsection{The infrared bound}\label{ss:IRbd}
Instead of deriving the exact solution for $\chi_p$, one may seek bounds
on $\chi_p$ or its derivative.  Indeed, it is not so difficult to show that
\cite{ms93}
\begin{align}
\frac{\chi_p^2}{1+\lambda\pc^2G_{\pc}^{*2}(o)}\le\frac{\mathrm{d}(p\chi_p)}
 {\mathrm{d}p}\le\chi_p^2.
\end{align}
The second inequality implies that $\chi_p$ is always bounded below by
$(1-p/\pc)^{-1}$.  Moreover, the first inequality implies that
$\chi_p$ is also bounded above by a multiple of $(1-p/\pc)^{-1}$, hence
$\gamma=1$, if
\begin{align}
G_{\pc}^{*2}(o)=\lim_{p\uparrow\pc}\int_{\Td}\hat G_p(k)^2
 \frac{\mathrm{d}^dk}{(2\pi)^d}<\infty,
\end{align}
where $\hat G_p(k)$ is the Fourier transform of the SAW two-point function
and $\Td\equiv[-\pi,\pi]^d$ is the $d$-dimensional torus of side length
$2\pi$ in the Fourier space.  It is a sufficient condition for the mean-field
behavior for $\chi_p$ and is called the bubble condition, named after the shape of
the diagram consisting of two line segments.  Whether or not the bubble
condition holds depends on the behavior of $\hat G_p(k)$ in the infrared
regime (i.e., around $k=0$).

For percolation, there is a similar condition to the bubble condition
under which $\gamma$ and other critical exponents take on their mean-field
values.  It is the cubic integrability of $\hat G_p(k)$ and is called the
triangle condition \cite{an84}.  Again, whether or not the triangle condition
holds depends on the infrared behavior of $\hat G_p(k)$.

Usually, there is no a priori bounds on $\hat G_p(k)$.  However, for some
spin models with a strong symmetry condition called reflection positivity
(e.g., the ferromagnetic Ising model with symmetric nearest-neighbor couplings
satisfies this condition), the two-point function enjoys the following
infrared bound \cite{fss76}: for any $d>2$, there is a
constant $K<\infty$ such that
\begin{align}\lbeq{IRbd-gen}
\|(1-\hat D)\hat G_p\|_\infty\equiv\sup_{k\in\Td}\big(1-\hat D(k)\big)|\hat
 G_p(k)|\le K\qquad\text{uniformly in $p$ close to }\pc.
\end{align}
If $D$ is a symmetric, non-degenerate and finite-range distribution with
variance $\sigma^2$, then $1-\hat D(k)\sim\frac{\sigma^2}{2d}|k|^2$ as
$|k|\to0$.  Suppose that the infrared bound holds for SAW and percolation.
Then
\begin{align}
G_{\pc}^{*n}(o)\le\int_{\Td}\bigg(\frac{K}{1-\hat D(k)}\bigg)^n
 \frac{\mathrm{d}^dk}{(2\pi)^d}\asymp\int_{\Td}
 \frac{\mathrm{d}^dk}{|k|^{2n}},
\end{align}
which implies that the bubble condition holds in all dimensions $d>4$ and
the triangle condition holds in all dimensions $d>6$.

On the other hand, there is some evidence (from hyperscaling inequalities,
numerical simulations, conformal field theory and so on) to suggest that
the critical exponents (if they exist) cannot take on their mean-field values
simultaneously if $d<4$ for SAW and $d<6$ for percolation.  In this respect,
the critical dimension $\dc$ is said to be 4 for SAW and 6 for percolation.

To complete the mean-field picture in high dimensions, it thus remains
to show that the infrared bound \refeq{IRbd-gen} holds for all dimensions
$d>\dc$.  Here, the lace expansion comes into play.

\subsection{The lace expansion}
In 1985, Brydges and Spencer \cite{bs85} came up to a fascinating idea.
First, they looked at the naive expansion \refeq{naive-exp}.  Next,
from each $\Gamma\in\cG[0,|\omega|]$, they isolated a connected graph
$\Gamma_0\subset\Gamma$ of the origin.  Then, they extracted a minimally
connected graph $\cL\subset\Gamma_0$ called a lace, and resummed all the
other edges in $\Gamma\setminus\cL$ to partially restore the self-avoidance
constraint.  This is what we nowadays call the algebraic lace expansion, named
after the shape of the aforesaid minimally connected graph.  Since then, the
algebraic lace expansion has been successfully applied to other models, such
as oriented percolation \cite{ny93}, lattice trees and lattice animals
\cite{hs90l}.

Later in 1990s, Hara and Slade (e.g., \cite{hs92a}) came up to a more
intuitively understandable way of deriving the lace expansion.  To
distinguish it from the algebraic lace expansion, we sometimes call it
the inclusion-exclusion lace expansion.  This opened up the possibility
of applying the lace expansion to a wider class of models, including
(unoriented) percolation \cite{hs90p}, the contact process \cite{s01},
the Ising model \cite{s07} and the (one-component) $\varphi^4$ model
\cite{s15}.

From now on, we simply call the latter the lace expansion.  We will show
its derivation for strictly SAW in Section~\ref{ss:derivation-SAW} and for
percolation in Section~\ref{ss:derivation-perc}.

The result of the lace expansion is formally explained by the following
recursion equation similar to that for the RW Green function: for any
$p<\pc$, there are functions $I_p$ and $J_p$ such that
\begin{align}
G_p(x)=I_p(x)+(J_p*G_p)(x).
\end{align}
If $I_p$ and $J_p$ satisfy certain regularity conditions, then it is natural
to believe that the global behavior of $G_p$ is also similar to that of the
RW Green function and therefore the infrared bound \refeq{IRbd-gen} holds.

However, since $I_p$ and $J_p$ are described by an alternating series
of the lace-expansion coefficients $\{\pi_p^{\sss(n)}\}_{n=0}^N$, each of which
involves complicated local interaction ($n$ represents the degree of
complexity), it is certainly not true that the aforesaid regularity
conditions always hold.  In fact, the regularity conditions require the
critical bubble $(D^{*2}*G_{\pc}^{*2})(o)$ for SAW and the critical triangle
$(D^{*2}*G_{\pc}^{*3})(o)$ for percolation to be small, not to be merely
finite.  This seemingly tautological statement (i.e., the critical
bubble/triangle have to be small in order to prove them to be finite) is
taken care of by the so-called bootstrapping argument, which will be
explained later in this survey.

During the course of the bootstrapping argument, we often assume that the
number of neighbors per vertex is sufficiently large.  Since each vertex has
$2d$ neighbors on $\Zd$, it means that $d$ is assumed to be large.  For SAW,
Hara and Slade \cite{hs92a,hs92b} succeeded in showing that $d\ge5$ is large
enough to prove mean-field results.  For percolation, however, the situation is
not as good as for SAW.  The best results so far were obtained by Fitzner and
van der Hofstad \cite{fh15b}, in which they proved mean-field results for
$d\ge11$ by using NoBLE, a perturbation method from non-backtracking
random walk (= memory-2 SAW).

There is another way to increase the number of neighbors per vertex.
Instead of taking $d$ large, we may enlarge the range $L$ of neighbors.
One such example is the spread-out lattice $\bar{\mZ}_L^d$, in which two
distinct vertices $x,y\in\Zd$ satisfying $\|x-y\|_\infty\le L$ are defined
to be neighbors, hence $(2L+1)^d-1$ neighbors per vertex.  By taking $L$
sufficiently large, all the models for which the lace expansion was obtained
are proven to exhibit mean-field behavior for all $d$ above the predicted
upper-critical dimensions \cite{hs90l,hs90p,ms93,ny93,s01,s07,s15}.

\subsection{The purposes of this survey}
Since we believe in universality, the mean-field results on the spread-out
lattice $\bar{\mZ}_L^d$, as long as $L<\infty$, are believed to hold on $\Zd$
as well.  This is proven to be true for SAW, but not yet for percolation.
We want to get rid of the artificial parameter $L$ and come up to a decent
nearest-neighbor lattice, on which 7-dimensional percolation is proven to
exhibit the mean-field behavior.  In an ongoing project with Lung-Chi Chen
and Markus Heydenreich \cite{chhks??}, we analyze the lace expansion for
percolation on a $d$-dimensional version of the body-centered cubic (BCC)
lattice, which has better features than the standard $\Zd$, as explained in
the next section.  Thanks to those features, enumeration of RW quantities
relevant to the lace-expansion analysis becomes extremely simple.  Also,
since those RW quantities are much smaller\footnote{A $d$-dimensional version
of the face-centered cubic (FCC) lattice has $d2^{d-1}$ neighbors per vertex,
more neighbors than on the BCC lattice, and therefore the RW quantities should
be much smaller on the FCC lattice.  However, since enumeration of those
quantities on the FCC lattice is not so simple (in fact, it is rather
complicated!), we decided to use the more charming BCC lattice.}
than the $\Zd$-counterparts, it is easy to get closer to the predicted
upper-critical dimension without introducing too much technical complexity.
One of the purposes of this survey is to explain the current status of the BCC
work and reveal the potential problems to overcome for completion of the
mean-field picture in high dimensions.

Another purpose of this survey is to provide a relatively short, self-contained
note on the lace expansion for the nearest-neighbor models.  Currently, the
best references on $\Zd$ are \cite{hs92a,hs92b} for SAW and \cite{fh15a,fh15b}
for percolation.  However, they are not necessarily accessible to beginners,
due to their length ($36 + 93$ pages for SAW and $79 + 92$ pages for
percolation) and complexity.  This is really unfortunate because, as mentioned
earlier, the lace expansion can provide a good playground for, e.g., graduate
students who may want to apply mathematical concepts and skills they learned
to interesting and important problems.  Considering this situation, we will
keep the material as simple as possible, instead of making all-out efforts to
go down to the predicted upper-critical dimensions.  That will be the final
goal of \cite{chhks??}.

\section{The models and the main result}
First, we provide precise definitions of the BCC lattice, self-avoiding walk
and percolation.  Then, we show the main result and explain its proof
assuming key propositions.

\subsection{The body-centered cubic (BCC) lattice}\label{ss:bcc}
\begin{figure}[t]
\begin{center}
\includegraphics[scale=0.7]{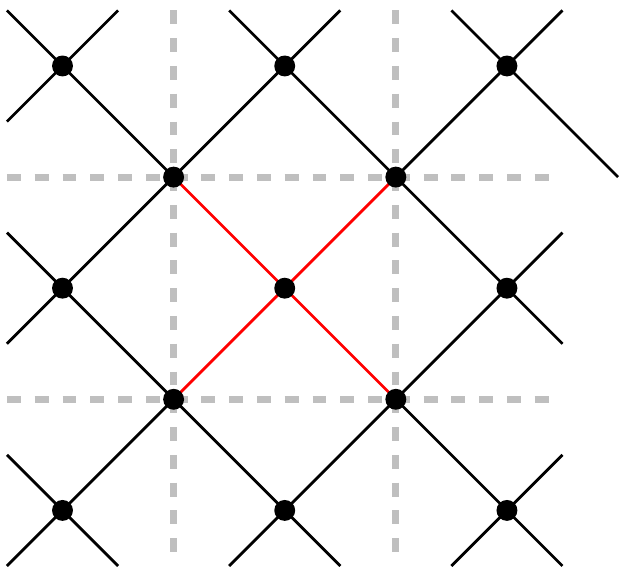}\hskip4pc
\includegraphics[scale=0.7]{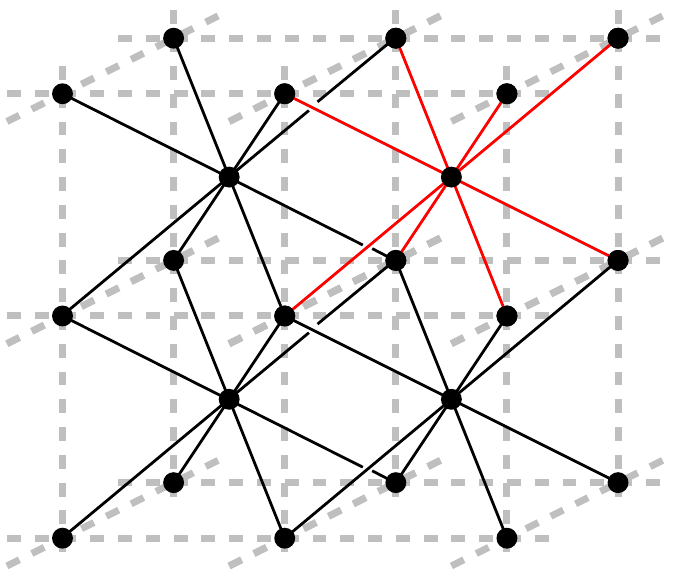}
\end{center}
\caption{\label{fig:bcc}
The basic structure (in red) of the BCC lattice $\Ld$ for $d=2,3$.}
\end{figure}

The $d$-dimensional BCC lattice $\Ld$ is a graph that contains the origin
$o=(0,\dots,0)$ and is generated by the set of neighbors
$\{x=(x_1,\dots,x_d):\prod_{j=1}^d|x_j|=1\}$.  It is equivalent to $\Zd$ when
$d=1$ and 2 (modulo rotation by $\pi/4$) but is more crowded in higher
dimensions in the sense that the degree of each vertex is $2^d$ on $\Ld$,
while it is $2d$ on $\Zd$.  We write $x\sim y$ if $x,y\in\Ld$ are neighbors,
i.e., $\prod_{j=1}^d|x_j-y_j|=1$.  It is a natural extension of the standard
3-dimensional BCC structure (see Figure~\ref{fig:bcc}).
The $d$-dimensional Brownian motion with the identity covariance matrix can
be constructed as the scaling limit of random walk (RW) on $\Ld$ generated
by the 1-step distribution
\begin{align}
D(x)=\frac1{2^d}\ind{x\sim o}=\prod_{j=1}^d\frac12\delta_{|x_j|,1}.
\end{align}
Due to this factorization and Stirling's formula\footnote{The two-sided bound~
$\dpst\frac1{12n+1}\le\log\frac{n!}{\sqrt{2\pi n}\big(\frac{n}e\big)^n}\le
\frac1{12n}$~ holds for all $n\in\mN$ \cite[Section~II.9]{f68}.}, we can obtain
a rather sharp bound on the $2n$-step return probability for all $n\in\mN$, as
\begin{align}\lbeq{stirling}
0\le(\pi n)^{-d/2}-D^{*2n}(o)\le\Big(1-e^{d(\frac1{24n+1}-\frac1{6n})}
 \Big)(\pi n)^{-d/2}\le\frac{2d}{15n}(\pi n)^{-d/2}.
\end{align}
Using this, we can easily evaluate various RW quantities, such as the RW loop
$\vep_1$, the RW bubble $\vep_2$ and the RW triangle $\vep_3$, defined as
\begin{align}\lbeq{RWquantities}
\vep_j&=(D^{*2}*S_1^{*j})(o)=\sum_{n=1}^\infty D^{*2n}(o)\times
 \begin{cases}
 1&[j=1],\\
 (2n-1)&[j=2],\\
 (2n-1)n&[j=3].
 \end{cases}
\end{align}
For example, if we split the sum into two at $n=N$, then the RW bubble
$\vep_2$ in dimensions $d>4$ can be estimated as
\begin{align}
0\le\vep_2-\sum_{n=1}^N(2n-1)D^{*2n}(o)\le2\pi^{-d/2}\int_N^\infty t^{1-d/2}\,
 \mathrm{d}t=\frac{4\pi^{-d/2}}{d-4}N^{(4-d)/2}.
\end{align}
If we choose $d=5$ and $N=100$ and use a calculator to evaluate the sum over
$n\le N$, then we obtain $\vep_2\le0.178465$.  Table~\ref{table:vep123}
summarizes the bounds on those RW quantities in different dimensions by
choosing $N=500$ (so that, by \refeq{stirling}, we can show that the RW
triangle $\vep_3$ for $d=7$ takes a value around the indicated number within
$10^{-6}$).
\begin{table}[b]
\caption{\label{table:vep123}Upper bounds on the RW loop, bubble and triangle
for $3\le d\le9$.}
\begin{center}
\begin{tabular}{r|lllllll}
         &    $d=3$ &    $d=4$ & $d=5$ & $d=6$ & $d=7$ & $d=8$ & $d=9$\\
\hline
$\vep_1$ & 0.393216 & 0.118637 & 0.046826 & 0.020461 & 0.009406
 & 0.004451 & 0.002144\\
$\vep_2$ & $\infty$ & $\infty$ & 0.178332 & 0.044004 & 0.015302
 & 0.006156 & 0.002678\\
$\vep_3$ & $\infty$ & $\infty$ & $\infty$ & $\infty$ & 0.052689
 & 0.012354 & 0.004148
\end{tabular}
\end{center}
\end{table}

\subsection{Self-avoiding walk}
As declared at the end of Section~\ref{s:intro}, we restrict our attention to
strictly SAW, which we simply call SAW from now on.  Let $\Omega(x,y)$ be the
set of self-avoiding paths on $\Ld$ from $x$ to $y$.  By convention,
$\Omega(x,x)$ is considered to be a singleton: a zero-step SAW at $x$.
Then, the SAW two-point function defined in the previous section can be
simplified as
\begin{align}\lbeq{Gsaw-def}
G_p(x)=\sum_{\omega\in\Omega(o,x)}p^{|\omega|}\prod_{j=1}^{|\omega|}D(\omega_j
 -\omega_{j-1}),
\end{align}
where the empty product is regarded as 1.  Recall that the susceptibility and
its critical point are defined as
\begin{align}\lbeq{chi&pc-def}
\chi_p=\sum_{x\in\Ld}G_p(x),&&
\pc=\sup\{p\ge0:\chi_p<\infty\}.
\end{align}

For more background and related results before 1993, we refer to the ``green"
book by Madras and Slade \cite{ms93}.  For recent progress in various
important problems, we refer to the monograph by Bauerschmidt et
al.~\cite{bdcgs12}.

\subsection{Percolation}
Here, we introduce bond percolation on $\Ld$.  Each bond $\{u,v\}\subset\Ld$
randomly takes either one of the two states, occupied or vacant, independently
of the other bonds.  We define the bond-occupation probability of a bond
$\{u,v\}$ as $pD(v-u)$, where $p\in[0,2^d]$ is the percolation
parameter, which is equal to the expected number of occupied bonds per
vertex.  Let $\mP_p$ be the associated probability measure, and denote its
expectation by $\mE_p$.

Next, we define the percolation two-point function.  In order to do so, we
first introduce the notion of connectivity.  We say that a self-avoiding path
$\omega=(\omega_0,\dots,\omega_{|\omega|})\in\Omega(x,y)$ is occupied if
either $x=y$ or every $b_j(\omega)\equiv\{\omega_{j-1},\omega_j\}$ for
$j=1,\dots,|\omega|$ is occupied.  We say that $x$ is connected to $y$,
denoted by $x\conn y$, if there is an occupied self-avoiding path
$\omega\in\Omega(x,y)$.  Then, we define the percolation two-point function as
\begin{align}\lbeq{perc2pt}
G_p(x)=\mP_p(o\conn x)
&=\mP_p\bigg(\bigcup_{\omega\in\Omega(o,x)}\{\omega\text{ is occupied}\}\bigg).
\end{align}
The susceptibility $\chi_p$ and its critical point $\pc$ are defined as in
\refeq{chi&pc-def}.  Menshikov \cite{m86} and Aizenman and Barsky \cite{ab87} independently proved that $\pc$ is unique in the sense that it can also be
characterized by the emergence of an infinite cluster of the origin:
\begin{align}
\pc=\inf\big\{p\in[0,2^d]:\mP_p(o\conn\infty)>0\big\}.
\end{align}
Recently, Duminil-Copin and Tassion \cite{dct16} came up to a particularly
simple proof of the uniqueness.  They also extended the idea to the Ising
model and dramatically simplified the proof of the uniqueness of the
critical temperature, first proven by Aizenman, Barsky and Fern\'andez \cite{abf87}.

For more background and related results before 1999, we refer to the excellent
book by Grimmett \cite{g99}.  The book by Bollob\'as and Riordan \cite{br06}
also contains progress after publication of Grimmett's book.

\subsection{The main result}
On the BCC lattice $\Ld$, we can prove the following result without introducing
too much technical complexity.

\begin{thm}[Infrared bound]\label{thm:main}
For SAW on $\mL^{\!d\ge6}$ and percolation on $\mL^{\!d\ge9}$,
there exists a model-dependent constant $K\in(0,\infty)$ such that
\begin{align}
\|(1-\hat D)\hat G_p\|_\infty\le K\qquad
 \text{uniformly in }p\in[1,\pc),
\end{align}
which implies the mean-field behavior, e.g., $\gamma=1$.
\end{thm}

In the proof of a key proposition necessary for the above theorem,
we will also show that
$\chi_1<\infty$.  This automatically implies the infrared bound for
$p\in[0,1)$, since
\begin{align}
\|(1-\hat D)\hat G_p\|_\infty\le2\chi_1<\infty.
\end{align}

The above result for SAW is not as sharp as the result in \cite{hs92a,hs92b},
where Hara and Slade proved the infrared bound on $\mZ^{d\ge5}$.  If we
simply follow their analysis with the same amount of work, then we should be
able to extend the above result to $\mL^{\!d\ge5}$.  However, as is mentioned
earlier, this is not our intention.  We include the result for SAW as an
example, just to show how easy to prove the infrared bound in such low
dimensions with relatively small effort.  Going down from 9 to 7 for
percolation will require more serious effort.  This will be the pursuit of
the joint work \cite{chhks??}.

The proof of the above theorem is rather straightforward, assuming the
following three propositions.  To state those propositions, we first
define
\begin{align}\lbeq{g12-def}
g_1(p)=p,&&
g_2(p)=\|(1-\hat D)\hat G_p\|_\infty.
\end{align}
Obviously, what we want to do is to show that $g_2(p)$ is bounded uniformly
in $p\in[1,\pc)$.  To define one more relevant function $g_3(p)$, we introduce
the notation for a sort of second derivative in the Fourier space, in a
particular direction.  For a function $\hat f$ on $\Td$ and
$k,l\in\Td$, we let
\begin{align}\lbeq{hatDelta}
\hat\Delta_k\hat f(l)=\frac{\hat f(l+k)+\hat f(l-k)}2-\hat f(l).
\end{align}
By simple trigonometric calculation, it is shown in
\cite[(5.17)]{s06}\footnote{It is shown in \cite[Lemma~5.7]{s06} that a
function $\hat A(k)=(1-\hat a(k))^{-1}$, where $\hat a$ is the Fourier
transform of a symmetric function $a(x)=a(-x)$ for all $x\in\Zd$, satisfies
the identity
\begin{align}\lbeq{trig-identity}
\hat\Delta_k\hat A(l)&=\frac{\hat A(l+k)+\hat A(l-k)}2\hat A(l)\,\hat\Delta_k
 \hat a(l)\nn\\
&\quad+\hat A(l+k)\hat A(l-k)\hat A(l)\bigg(\sum_xa(x)(\sin l\cdot x)(\sin k
 \cdot x)\bigg)^2.
\end{align}
The inequality \refeq{U-def} is obtained by applying the Schwarz inequality
to the sum in the above expression.} that the Fourier transform of the RW
Green function $\hat S_1(k)\equiv(1-\hat D(k))^{-1}$, which is well-defined in
a proper limit when $d>2$, obeys the inequality
\begin{align}\lbeq{U-def}
&|\hat\Delta_k\hat S_1(l)|\nn\\
&\le\hat U(k,l)\equiv\big(1-\hat D(k)\big)\bigg(\frac{\hat S_1(l+k)+\hat S_1
 (l-k)}2\hat S_1(l)+4\hat S_1(l+k)\hat S_1(l-k)\bigg).
\end{align}
Finally, we define
\begin{align}\lbeq{g3-def}
g_3(p)=\sup_{k,l}\frac1{\hat U(k,l)}\times
 \begin{cases}
 |\hat\Delta_k\hat G_p(l)|&[\mathrm{SAW}],\\
 |\hat\Delta_k(\hat G_p(l)p\hat D(l))|&[\mathrm{percolation}],
 \end{cases}
\end{align}
where the supremum near $k=0$ should be interpreted as the supremum over
the limit as $|k|\to0$.  It will be clear that $g_3$ is defined in slightly
different ways between the two models, due to the difference in the recursion
equations obtained by the lace expansion.

Now, we state the aforementioned three propositions and show that
they indeed imply Theorem~\ref{thm:main}.

\begin{prp}[Continuity]\label{prp:f-cont}
The functions $\{g_i(p)\}_{i=1}^3$ are continuous in $p\in[1,\pc)$.
\end{prp}

\begin{prp}[Initial conditions]\label{prp:f-initial}
For SAW on $\mL^{\!d\ge6}$ and percolation on $\mL^{\!d\ge8}$, 
there are model-dependent finite constants $\{K_i\}_{i=1}^3$ such that 
$g_i(1)<K_i$ for $i=1,2,3$.
\end{prp}

\begin{prp}[Bootstrapping argument]\label{prp:f-bootstrapping}
For SAW on $\mL^{\!d\ge6}$ and percolation on $\mL^{\!d\ge9}$, we fix
$p\in(1,\pc)$ and assume $g_i(p)\le K_i$, $i=1,2,3$, where
$\{K_i\}_{i=1}^3$ are the same constants as in Proposition~\ref{prp:f-initial}.
Then, the stronger inequalities $g_i(p)<K_i$, $i=1,2,3$, hold.
\end{prp}

\begin{figure}[t]
\begin{center}
\includegraphics[scale=0.5]{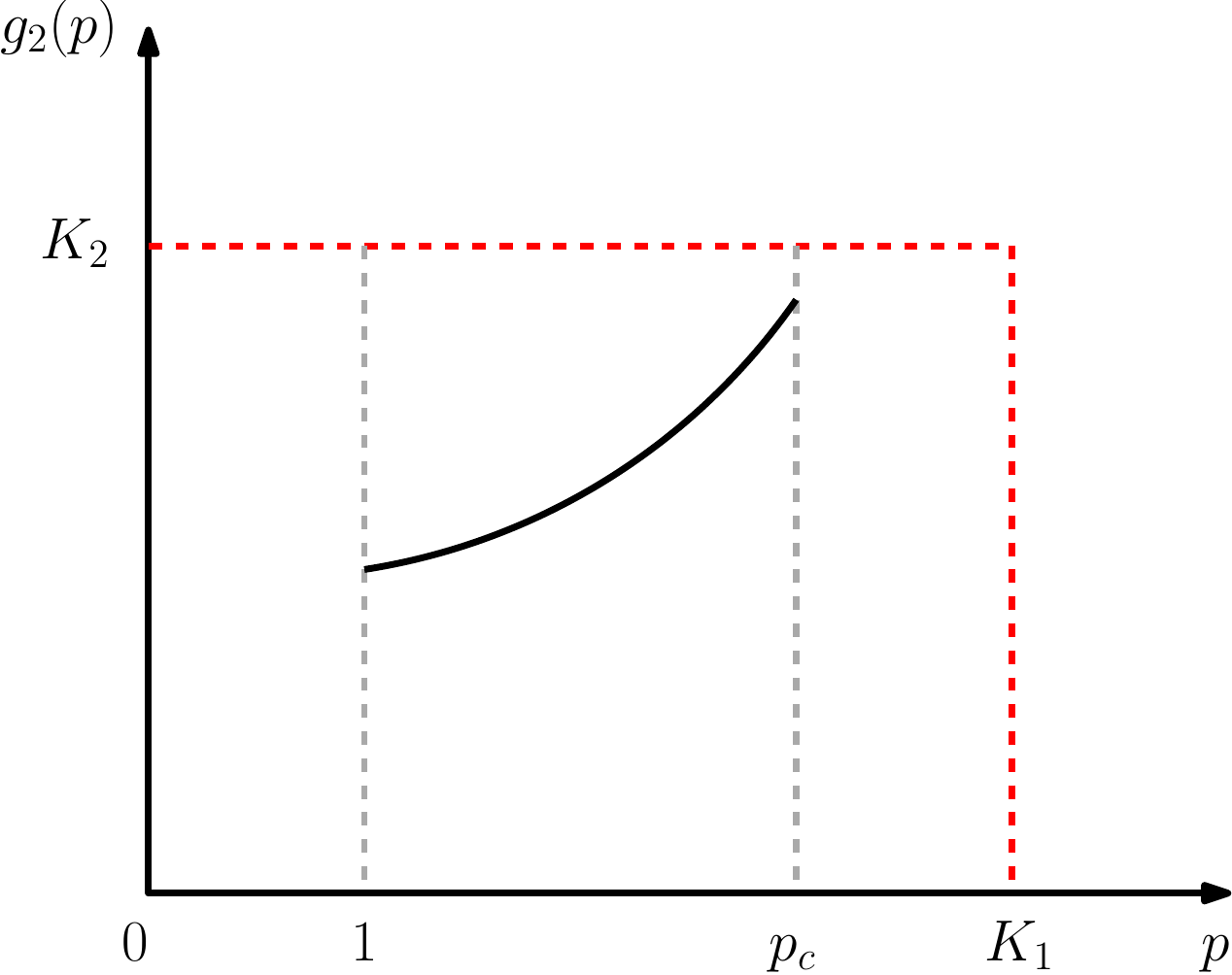}
\end{center}
\caption{\label{fig:g2}
Depiction of the proof of Theorem~\ref{thm:main} assuming
Propositions~\ref{prp:f-cont}--\ref{prp:f-bootstrapping}.}
\end{figure}

Since $g_2(p)$ is continuous in $p\in[1,\pc)$, with the initial value
$g_2(1)<K_2$, and cannot be equal to $K_2$ for $p\in(1,\pc\wedge K_1)$,
we can say that the strict inequality $g_2(p)<K_2$ holds for all
$p\in[1,\pc\wedge K_1)$.  Since the same argument applies to
$g_1(p)$, we can conclude $\pc\le K_1$, hence $g_2(p)<K_2$ for
all $p\in[1,\pc)$ (see Figure~\ref{fig:g2}).  This completes the proof of
Theorem~\ref{thm:main} assuming
Propositions~\ref{prp:f-cont}--\ref{prp:f-bootstrapping}.
\QED

\subsection{Where and how to use the lace expansion}\label{ss:whereandhow}
It remains to prove Propositions~\ref{prp:f-cont}--\ref{prp:f-bootstrapping}.
The proof of Proposition~\ref{prp:f-cont} is elementary, though cumbersome for
$g_3(p)$, and is explained in the next section.  To prove the other two
propositions, we will use the following lace expansion.

\begin{prp}[Lace expansion]\label{prp:lace-gen}
For any $p<\pc$ and $N\in\Zp\equiv\{0\}\cup\mN$, there exist
model-dependent nonnegative functions $\{\pi_p^{\sss(n)}\}_{n=0}^N$ on
$\Ld$ ($\pi_p^{\sss(0)}\equiv0$ for SAW) such that, if we define
$I_p^{\sss(N)}$ and $J_p^{\sss(N)}$ as
\begin{align}
I_p^{\sss(N)}(x)&=\delta_{o,x}+
 \begin{cases}
 0&[\text{SAW}],\\
 \sum_{n=0}^N(-1)^n\pi_p^{\sss(n)}(x)&[\text{percolation}],
 \end{cases}\\
J_p^{\sss(N)}(x)&=pD(x)+
 \begin{cases}
 \sum_{n=1}^N(-1)^n\pi_p^{\sss(n)}(x)&[\text{SAW}],\\
 \sum_{n=0}^N(-1)^n(\pi_p^{\sss(n)}*pD)(x)&[\text{percolation}],
 \end{cases}
\end{align}
then we obtain the recursion equation
\begin{align}\lbeq{le-gen}
G_p(x)=I_p^{\sss(N)}(x)+(J_p^{\sss(N)}*G_p)(x)+
 (-1)^{N+1}R_p^{\sss(N+1)}(x),
\end{align}
where the remainder $R_p^{\sss(N)}$ obeys the bound
\begin{align}\lbeq{rem-gen}
0\le R_p^{\sss(N)}(x)\le(\pi_p^{\sss(N)}*G_p)(x).
\end{align}
\end{prp}

The derivation of the lace expansion is model-dependent and is explained for
SAW in Section~\ref{ss:derivation-SAW} and for percolation in
Section~\ref{ss:derivation-perc}.

Here, we briefly explain where and how to use the lace expansion
to prove Propositions~\ref{prp:f-initial}--\ref{prp:f-bootstrapping}.
The details will be given in later sections.

\Proof{Step 1.}
First, we evaluate $\{g_i(p)\}_{i=1}^3$ in terms of sums of
$\hat\pi_p^{\sss(n)}(k)\equiv\sum_xe^{ik\cdot x}\pi_p^{\sss(n)}(x)$.

(i) Let $p\in[1,\pc)$ and 
suppose $\sum_{n=0}^\infty\hat\pi_p^{\sss(n)}(0)$ is small enough 
to ensure that
\begin{align}
\lim_{N\to\infty}\hat\pi_p^{\sss(N)}(0)=0,&&
\hat I_p(k)\equiv\lim_{N\to\infty}\hat I_p^{\sss(N)}(k)>0~~
 \text{uniformly in }p\text{ and }k.
\end{align}
The latter is always true for SAW since $\hat I_p(k)\equiv1$.  The former 
implies that
\begin{align}
0\le\sum_{x\in\Ld}R_p^{\sss(N)}(x)\le\hat\pi_p^{\sss(N)}(0)\chi_p
 \xrightarrow[N\to\infty]{}0.
\end{align}
Let (n.b.~$\pi_p^{\sss(0)}\equiv0$ for SAW)
\begin{align}
\hat\Pi_p(k)=\sum_{n=0}^\infty(-1)^n\hat\pi_p^{\sss(n)}(k),&&
\hat J_p(k)=p\hat D(k)+
 \begin{cases}
 \hat\Pi_p(k)&[\text{SAW}],\\
 \hat\Pi_p(k)p\hat D(k)&[\text{percolation}].
 \end{cases}
\end{align}
Then, by using \refeq{le-gen}, we obtain
\begin{align}\lbeq{chi-gen}
\chi_p\equiv\hat G_p(0)=\hat I_p(0)+\hat J_p(0)\chi_p=\frac{\hat I_p(0)}
 {1-\hat J_p(0)}.
\end{align}
Since $\chi_p\ge0$ and 
$\hat I_p(0)>0$, we can conclude $\hat J_p(0)\le1$, which implies
\begin{align}\lbeq{g1bd-byPi}
g_1(p)\le
 \begin{cases}
 1-\hat\Pi_p(0)&[\text{SAW}],\\
 (1+\hat\Pi_p(0))^{-1}&[\text{percolation}].
 \end{cases}
\end{align}

(ii) Next, by \refeq{le-gen} and \refeq{chi-gen}, we obtain
\begin{align}\lbeq{hatG-gen}
\hat G_p(k)=\frac{\hat I_p(k)}{1-\hat J_p(k)}=\frac{\hat I_p(k)}
 {-\hat\Delta_k\hat J_p(0)+\hat I_p(0)/\chi_p},
\end{align}
where we have used the symmetry $\hat J_p(k)=\hat J_p(-k)$ to obtain
$-\hat\Delta_k\hat J_p(0)=\hat J_p(0)-\hat J_p(k)$.  Suppose 
$-\sum_{n=0}^\infty\hat\Delta_k\hat\pi_p^{\sss(n)}(0)\equiv\sum_{n=0}^\infty
\sum_x(1-\cos k\cdot x)\pi_p^{\sss(n)}(x)$ is smaller than $1-\hat D(k)$
in order to ensure $-\hat\Delta_k\hat J_p(0)\ge0$.  Then, $\hat G_p(k)$ is
bounded as\footnote{For percolation, the non-negativity of $\hat G_p(k)$ is
elementary and proven in \cite[Lemma~3.3]{an84}.  The
actual proof goes as follows.  First, by translation-invariance, we can use any
vertex $y$ to rewrite $\hat G_p(k)$ as
\begin{align}
\hat G_p(k)=\sum_xe^{ik\cdot x}\mP_p(o\conn x)=\sum_xe^{ik\cdot x}\mP_p(y
 \conn x+y)=\mE_p\bigg[\sum_ze^{ik\cdot(z-y)}\ind{y\conn z}\bigg].
\end{align}
Then, by using the identity $1=\sum_y\ind{y\in\cC(o)}/|\cC(o)|$, where $\cC(o)$ is
the set of vertices connected from $o$, we can rewrite the rightmost expression
as
\begin{align}
\mE_p\bigg[\frac1{|\cC(o)|}\sum_{y\in\cC(o)}\sum_ze^{ik\cdot(z-y)}\ind{y\conn
 z}\bigg]
 =\mE_p\bigg[\frac1{|\cC(o)|}\sum_{y,z\in\cC(o)}e^{ik\cdot(z-y)}\bigg]
 =\mE_p\bigg[\bigg|\frac1{\sqrt{|\cC(o)|}}\sum_{z\in\cC(o)}e^{ik\cdot z}\bigg|^2
 \bigg]\ge0.
\end{align}}
\begin{align}\lbeq{irbd-gen}
0\le\hat G_p(k)\le\frac{\hat I_p(k)}{-\hat\Delta_k\hat J_p(0)}.
\end{align}
Since $p\ge1$, this implies
\begin{align}\lbeq{g2bd-byPi}
g_2(p)\le
 \begin{cases}
 \dpst\sup_k\bigg(1+\frac{-\hat\Delta_k\hat\Pi_p(0)}{1-\hat D(k)}
  \bigg)^{-1}&[\text{SAW}],\\
 \dpst\sup_k\bigg(1+\frac1{\hat I_p(k)}\frac{-\hat\Delta_k\hat\Pi_p(0)}
  {1-\hat D(k)}\bigg)^{-1}&[\text{percolation}],
 \end{cases}
\end{align}
where the supremum near $k=0$ should be interpreted as the supremum over the
limit as $|k|\to0$.

(iii) To evaluate $g_3(p)$, we want to use the identity \refeq{trig-identity}.
To do so for percolation, 
we first notice that, by using $\hat I_p(k)p\hat D(k)=\hat J_p(k)$ and 
\refeq{hatG-gen}, we obtain
\begin{align}
\hat G_p(k)p\hat D(k)=\frac{\hat J_p(k)}{1-\hat J_p(k)}=\frac1{1-\hat J_p(k)}
 -1\equiv\hat A_p(k)-1,
\end{align}
hence $\hat\Delta_k(\hat G_p(l)p\hat D(l))=\hat\Delta_k\hat A_p(l)$.
As a result, $g_3(p)$ for both models can be written as
\begin{align}
g_3(p)=\sup_{k,l}\frac{|\hat\Delta_k\hat A_p(l)|}{\hat U(k,l)}.
\end{align}
Then, by using \refeq{trig-identity} with $a(x)=J_p(x)$, 
noting $\hat A_p(k)=(-\hat\Delta_k\hat J_p(0)+\hat I_p(0)/\chi_p)^{-1}\ge0$ 
and applying the Schwarz inequality as in \cite[Lemma~5.7]{s06}, we obtain
\begin{align}\lbeq{g3bd-byPi}
g_3(p)\le\sup_{k,l}\frac{1-\hat D(k)}{\hat U(k,l)}\bigg(
&\frac{\hat A_p(l+k)+\hat A_p(l-k)}2\hat A_p(l)\frac{|\hat\Delta_k\hat J_p(l)|}
 {1-\hat D(k)}\nn\\
&+4\hat A_p(l+k)\hat A_p(l-k)\frac{-\hat\Delta_l\widehat{|J_p|}(0)}{1-\hat J_p
 (l)}\,\frac{-\hat\Delta_k\widehat{|J_p|}(0)}{1-\hat D(k)}\bigg),
\end{align}
where
\begin{align}
\widehat{|J_p|}(k)=\sum_{x\in\Ld}e^{ik\cdot x}|J_p(x)|.
\end{align}
We can further bound $|\hat\Delta_k\hat J_p(l)|$ and $-\hat\Delta_k
\widehat{|J_p|}(0)\equiv\widehat{|J_p|}(0)-\widehat{|J_p|}(k)\ge0$ in terms
of sums of $|\hat\Delta_k\hat\pi_p^{\sss(n)}(0)|$.  However, to simplify the
exposition, we refrain from doing so for now and postpone it to later sections.

So far, we have assumed that $\sum_{n=0}^\infty\hat\pi_p^{\sss(n)}(0)$ and 
$-\sum_{n=0}^\infty\hat\Delta_k\hat\pi_p^{\sss(n)}(0)$ are small enough to 
carry out the above computations.  Sufficient conditions to this assumption 
are
\begin{align}\lbeq{suffcond-saw}
\sum_{n=1}^\infty\hat\pi_p^{\sss(n)}(0)<\infty,&&
\sup_k\sum_{n=1}^\infty\frac{-\hat\Delta_k\hat\pi_p^{\sss(n)}(0)}{1-\hat
 D(k)}<1
\end{align}
for SAW, and
\begin{align}\lbeq{suffcond-perc}
\sum_{n=0}^\infty\hat\pi_p^{\sss(n)}(0)+\sup_k\sum_{n=0}^\infty\frac{-\hat
 \Delta_k\hat\pi_p^{\sss(n)}(0)}{1-\hat D(k)}<1
\end{align}
for percolation (cf., \refeq{g2bd-prebyBWr}).  
These conditions are to be verified eventually.
\QED

\Proof{Step 2.}
As shown in \refeq{g1bd-byPi}, \refeq{g2bd-byPi} and \refeq{g3bd-byPi}, the
bootstrapping functions $\{g_i(p)\}_{i=1}^3$ are bounded in terms of sums of
$\hat\pi_p^{\sss(n)}(0)$ and sums of $|\hat\Delta_k\hat\pi_p^{\sss(n)}(0)|$.
In the second step, we evaluate those lace-expansion coefficients in terms of
smaller diagrams, such as
\begin{align}\lbeq{LBT-def}
L_p=\|(pD)^{*2}*G_p\|_\infty,&&
B_p=\|(pD)^{*2}*G_p^{*2}\|_\infty,&&
T_p=\|(pD)^{*2}*G_p^{*3}\|_\infty.
\end{align}
For example, we can bound $\hat\pi_p^{\sss(n)}$ for $n\ge2$ as
\begin{align}
0\le\hat\pi_p^{\sss(n)}(0)\le
 \begin{cases}
 B_p(p\|D\|_\infty+L_p)r^{n-2}&[\text{SAW}],\\
 (1+\frac12B_p+T_p)^2r\rho^{n-1}&[\text{percolation}],
 \end{cases}
\end{align}
where
\begin{align}\lbeq{common-ratios}
r=p\|D\|_\infty+L_p+B_p,&&
\rho=\bigg(1+\frac12B_p+T_p\bigg)(r+T_p)+T_p(2r+T_p)
\end{align}
See Sections~\ref{s:LESAW}--\ref{s:LEperc} for the proof of the above
inequality and the bounds on $\hat\pi_p^{\sss(0)}(0)$ and $\hat\pi_p^{\sss(1)}(0)$.
It will also be shown that the amplitude of
$|\hat\Delta_k\hat\pi_p^{\sss(n)}(0)|/(1-\hat D(k))$ is bounded in a similar
fashion, with the common ratio $r$ for SAW and $\rho$ for percolation.
Therefore, the assumptions made in Step~1 hold if $L_p,B_p,T_p$ and other
diagrams in the bounds are small enough.
\QED

\Proof{Step 3.}
In the final step, we investigate the aforesaid diagrams and prove that,
by choosing appropriate values for $\{K_i\}_{i=1}^3$, those diagrams are
indeed small enough for SAW on $\mL^{\!d\ge6}$ and for percolation on
$\mL^{\!d\ge9}$.

(i) For $p=1$, we only need to use the trivial inequality $G_1(x)\le S_1(x)$,
$x\in\Ld$, for both models to obtain that, for $d>2$ (as mentioned earlier,
$\hat S_1(k)\equiv(1-\hat D(k))^{-1}$ is well-defined in a proper limit when
$d>2$),
\begin{align}\lbeq{L1bd-gen}
L_1\le\|D^{*2}*S_1\|_\infty=\int_{\Td}\frac{\hat D(k)^2}{1-\hat
 D(k)}\,\frac{\text{d}^dk}{(2\pi)^d}=(D^{*2}*S_1)(o)\equiv\vep_1.
\end{align}
Similarly, we obtain
\begin{align}\lbeq{B1T1bd-gen}
B_1\le\vep_2,&&
T_1\le\vep_3.
\end{align}
Consulting with Table~\ref{table:vep123} in Section~\ref{ss:bcc}, we can see
that, even in $\dc+1$ dimensions, $r$ and $\rho$ in \refeq{common-ratios} are
small enough for the bootstrapping functions $\{g_i(p)\}_{i=1}^3$ to be
convergent.

(ii) The strategy for $p\in(1,\pc)$ is different from that for $p=1$, because
there is no \emph{a priori} bound on $G_p$ in terms of $S_1$.  Here, we use
the assumptions $g_i(p)\le K_i$, $i=1,2,3$, to evaluate the diagrams.  For
example,
\begin{align}\lbeq{Lpbd-gen}
L_p\le p^2\int_{\Td}\hat D(k)^2|\hat G_p(k)|\,\frac{\text{d}^dk}
 {(2\pi)^d}\le K_1^2K_2\underbrace{\int_{\Td}\frac{\hat D(k)^2}
 {1-\hat D(k)}\,\frac{\text{d}^dk}{(2\pi)^d}}_{=\vep_1}.
\end{align}
Similarly,
\begin{align}\lbeq{BpTpbd-gen}
B_p\le K_1^2K_2^2\vep_2,&&
T_p\le K_1^2K_2^3\vep_3.
\end{align}
As a result, $r$ and $\rho$ in \refeq{common-ratios} become functions of
$\{K_i\}_{i=1,2}$.  If we choose their values appropriately, then we can derive
the improved bound $g_1(p)<K_1$ for all $d\ge\dc+1$.  To improve the bounds on
$\{g_i(p)\}_{i=2,3}$, we also have to control $K_3$.  This is the worst enemy
that keeps us from going down to $\dc+1$ dimensions.  In \cite{chhks??},
we will make all-out efforts to overcome this problem.
\QED

\subsection{Organization}
In the rest of this survey, we prove the above propositions in detail.  In
Section~\ref{s:cont}, we prove Proposition~\ref{prp:f-cont} for both models.

In Section~\ref{s:LESAW}, we prove
Propositions~\ref{prp:f-initial}--\ref{prp:lace-gen} for SAW as follows.
In Section~\ref{ss:derivation-SAW}, we explain the derivation of the lace
expansion (Proposition~\ref{prp:lace-gen}) for SAW.
In Section~\ref{ss:Pibds-saw}, we prove bounds on the lace-expansion
coefficients in terms of basic diagrams, as briefly explained in \emph{Step~2}
in Section~\ref{ss:whereandhow}.
In Section~\ref{ss:LBWbd}, we prove bounds on those basic diagrams in terms of
RW quantities, as explained in \emph{Step~3} in Section~\ref{ss:whereandhow}.  
Applying them to the bounds on the bootstrapping functions
$\{g_i(p)\}_{i=1}^3$ obtained in \emph{Step~1} in
Section~\ref{ss:whereandhow}, we prove 
Propositions~\ref{prp:f-initial}--\ref{prp:f-bootstrapping} on 
$\mL^{\!d\ge6}$.  Finally, in Section~\ref{ss:discussion-saw}, 
we provide further discussion to potentially improve our results.

In Section~\ref{s:LEperc}, we prove
Propositions~\ref{prp:f-initial}--\ref{prp:lace-gen} for percolation as
follows. In Section~\ref{ss:derivation-perc}, we derive the lace
expansion (Proposition~\ref{prp:lace-gen}) for percolation.
In Section~\ref{ss:Pibds-perc}, we prove bounds on the lace-expansion
coefficients in terms of basic diagrams, as briefly explained in \emph{Step~2}
in Section~\ref{ss:whereandhow}.
In Section~\ref{ss:TVOHbd}, we prove bounds on those basic diagrams in terms of
RW quantities, as explained in \emph{Step~3} in Section~\ref{ss:whereandhow}.  
Applying them to the bounds on the bootstrapping functions
$\{g_i(p)\}_{i=1}^3$ obtained in \emph{Step~1} in
Section~\ref{ss:whereandhow}, we prove Proposition~\ref{prp:f-initial} on 
$\mL^{\!d\ge8}$ and Proposition~\ref{prp:f-bootstrapping} on 
$\mL^{\!d\ge9}$.  In Section~\ref{ss:discussion-perc}, 
we provide further discussion to potentially improve our results.

\section{Continuity of the bootstrapping functions}\label{s:cont}
In this section, we prove Proposition~\ref{prp:f-cont}.  First, we recall
\refeq{g12-def} and \refeq{g3-def} for the bootstrapping functions
$\{g_i(p)\}_{i=1}^3$.  Obviously, $g_1(p)\equiv p$ is continuous.
To prove continuity of the other two, we introduce
\begin{align}
\tilde g_{2,k}(p)&=\big(1-\hat D(k)\big)\hat G_p(k),\\
\tilde g_{3,k,l}(p)&=\frac1{\hat U(k,l)}\times
 \begin{cases}
 \hat\Delta_k\hat G_p(l)&[\text{SAW}],\\
 \hat\Delta_k(\hat G_p(l)\hat D(l))&[\text{percolation}],
 \end{cases}
\end{align}
and show that they are continuous in $p\in[1,\pc)$ for every
$k,l\in\Td$.  However, since
\begin{align}
g_2(p)=\sup_{k\in\Td}|\tilde g_{2,k}(p)|,&&
g_3(p)=\sup_{k,l\in\Td}|\tilde g_{3,k,l}(p)|,
\end{align}
and the supremum of continuous functions is not necessarily continuous,
we must be a bit more cautious here.  The following elementary lemma provides
a sufficient condition for the supremum to be continuous.

\begin{lmm}[Lemma~5.13 of \cite{s06}, in our language]
Fix $p_0\in[1,\pc)$ and let $\{\hat f_k(p)\}_{k\in\Td}$ be an
equicontinuous family of functions in $p\in[1,p_0]$.  Suppose that
$\sup_{k\in\Td}\hat f_k(p)<\infty$ for every $p\in[1,p_0]$.
Then, $\sup_{k\in\Td}\hat f_k(p)<\infty$ is continuous in
$p\in[1,p_0]$.
\end{lmm}

Therefore, in order to prove continuity of $\{g_i(p)\}_{i=2,3}$ in
$p\in[1,\pc)$, we want to show that
$\{\tilde g_{2,k}(p)\}_{k\in\Td}$ and
$\{\tilde g_{3,k,l}(p)\}_{k,l\in\Td}$ are equicontinuous
families of functions in $p\in[1,p_0]$ for each $p_0\in[1,\pc)$.
To prove this, it then suffices to show that the following (i) and (ii) hold.
\begin{enumerate}[(i)]
\item
$\tilde g_{2,k}(p)$ and $\partial_p\tilde g_{2,k}(p)$
are finite uniformly in $k\in\Td$ and $p\in[1,p_0]$.
\item
$\tilde g_{3,k,l}(p)$ and $\partial_p\tilde g_{3,k,l}(p)$
are finite uniformly in $k,l\in\Td$ and $p\in[1,p_0]$.
\end{enumerate}

To prove (i) is not so hard.  By $0\le1-\hat D(k)\le2$,
$|\hat G_p(k)|\le\chi_p$ and the monotonicity of $\chi_p$ in $p$,
we obtain $|\tilde g_{2,k}(p)|\le2\chi_{p_0}<\infty$ uniformly in
$k\in\Td$ and $p\in[1,p_0]$.  Moreover, by subadditivity
for SAW, Russo's formula and the BK inequality for percolation
(see, e.g., \cite{g99}), and then using translation-invariance, we obtain
\begin{align}\lbeq{subaddBKbd}
0\le\partial_pG_p(x)\le(D*G_p^{*2})(x),
\end{align}
hence
\begin{align}
|\partial_p\tilde g_{2,k}(p)|\le2\sum_x(D*G_p^{*2})(x)\le2\chi_{p_0}^2
 <\infty,
\end{align}
uniformly in $k\in\Td$ and $p\in[1,p_0]$, as required.

To prove (ii) needs extra care, especially near $k=0$, because of the factor
$1-\hat D(k)$ in $\hat U(k,l)$.  From here, we prove (ii) for SAW and for
percolation separately.

\Proof{Proof of (ii) for SAW.}
First, by using the telescopic inequality in
\cite[Appendix~A]{fh15a}\footnote{Although \refeq{telescope} is a result of
simple trigonometric computation, it is not so easy to come up to the actual
proof.  The actual proof of \cite[Appendix~A]{fh15a} goes as follows.  First,
take the real part of the telescopic identity
$1-\exp(i\sum_{j=1}^Jt_j)=\sum_{j=1}^J(1-e^{it_j})\exp(i\sum_{h=1}^{j-1}t_h)$,
where the empty sum for $j=1$ is regarded as zero.  Then, use the
inequalities $|\sin\sum_{h=1}^{j-1}t_h|\le\sum_{h=1}^{j-1}|\sin t_h|$,
$|\sin t_j||\sin t_h|\le(\sin^2t_j+\sin^2t_h)/2$ and
$\sin^2t_j\le2(1-\cos t_j)$ to obtain
\begin{align}
1-\cos\bigg(\sum_{j=1}^Jt_j\bigg)-\sum_{j=1}^J(1-\cos t_j)&=-\sum_{j=1}^J
 \underbrace{(1-\cos t_j)\bigg(1-\cos\sum_{h=1}^{j-1}t_h\bigg)}_{\ge0}+\sum_{j
 =1}^J(\sin t_j)\sin\sum_{h=1}^{j-1}t_h\nn\\
&\le\sum_{j=1}^J\sum_{h=1}^{j-1}\frac{\sin^2t_j+\sin^2t_h}2
 \le(J-1)\sum_{j=1}^J(1-\cos t_j),
\end{align}
which implies \refeq{telescope}.}
\begin{align}\lbeq{telescope}
0\le1-\cos\sum_{j=1}^Jt_j\le J\sum_{j=1}^J(1-\cos t_j),
\end{align}
we obtain
\begin{align}\lbeq{SAW-DeltaGbd1}
&|\hat\Delta_k\hat G_p(l)|\le\sum_x(1-\cos k\cdot x)G_p(x)\nn\\
&\quad=\sum_x\sum_{\omega\in\Omega(o,x)}\bigg(1-\cos\sum_{i=1}^{|\omega|}k\cdot
 (\omega_i-\omega_{i-1})\bigg)p^{|\omega|}\prod_{j=1}^{|\omega|}D(\omega_j
 -\omega_{j-1})\nn\\
&\quad\le\sum_{u,v,x}\big(1-\cos k\cdot(v-u)\big)\sum_{\omega\in\Omega(o,x)}
 |\omega|\sum_{i=1}^{|\omega|}\ind{b_i(\omega)=(u,v)}\,p^{|\omega|}\prod_{j
 =1}^{|\omega|}D(\omega_j-\omega_{j-1}).
\end{align}
Ignoring the self-avoidance constraint between
$\eta\equiv(\omega_0,\dots,\omega_{i-1})$ and
$\xi\equiv(\omega_i,\dots,\omega_{|\omega|})$ and using
translation-invariance, we can further bound $|\hat\Delta_k\hat G_p(l)|$ as
\begin{align}\lbeq{SAW-DeltaGbd2}
|\hat\Delta_k\hat G_p(l)|&\le\sum_{u,v,x}\big(1-\cos k\cdot(v-u)\big)
 pD(v-u)\sum_{\substack{\eta\in\Omega(o,u)\\ \xi\in\Omega(v,x)}}(|\eta|+|\xi|
 +1)\nn\\
&\hskip3pc\times p^{|\eta|}\prod_{i=1}^{|\eta|}D(\eta_i-\eta_{i-1})~p^{|\xi|}
 \prod_{j=1}^{|\xi|}D(\xi_j-\xi_{j-1})\nn\\
&\le2p\big(1-\hat D(k)\big)\chi_p\sum_x\sum_{\omega\in\Omega(o,x)}
 (|\omega|+1)p^{|\omega|}\prod_{j=1}^{|\omega|}D(\omega_j-\omega_{j-1}).
\end{align}
However, by the identity $|\omega|+1=\sum_y\ind{y\in\omega}$ for a
self-avoiding path $\omega$, subadditivity and translation-invariance,
the sum in the last line is bounded as
\begin{align}
\sum_x\sum_{\omega\in\Omega(o,x)}(|\omega|+1)p^{|\omega|}\prod_{j=1}^{|\omega|}
 D(\omega_j-\omega_{j-1})\le\sum_{x,y}G_p(y)G_p(x-y)=\chi_p^2.
\end{align}
As a result, we arrive at
\begin{align}\lbeq{SAW-DeltaGbd3}
|\hat\Delta_k\hat G_p(l)|\le2p_0\big(1-\hat D(k)\big)\chi_{p_0}^3,
\end{align}
which implies that $\tilde g_{3,k,l}(p)$ is finite uniformly in
$k,l\in\Td$ and $p\in[1,p_0]$.

For the derivative $\partial_p\tilde g_{3,k,l}(p)\equiv\hat U(k,l)^{-1}
\hat\Delta_k\partial_p\hat G_p(l)$, we note that
\begin{eqnarray}\lbeq{SAW-DeltaGbd4}
|\hat\Delta_k\partial_p\hat G_p(l)|&\stackrel{\refeq{subaddBKbd}}\le&\sum_x
 (1-\cos k\cdot x)(D*G_p^{*2})(x)\nn\\
&\stackrel{\refeq{telescope}}\le&3\bigg(\big(1-\hat D(k)\big)\chi_p^2+2\chi_p
 \underbrace{\sum_v(1-\cos k\cdot v)G_p(v)}_{=\hat\Delta_k\hat G_p
 (0)}\bigg)\nn\\
&\stackrel{\refeq{SAW-DeltaGbd3}}\le&3(1-\hat D(k))\chi_{p_0}^2(1+4p_0
 \chi_{p_0}^2).
\end{eqnarray}
Therefore, $\partial_p\tilde g_{3,k,l}(p)$ is also finite uniformly in
$k,l\in\Td$ and $p\in[1,p_0]$.
\QED

\Proof{Proof of (ii) for percolation.}
First, we note that
\begin{eqnarray}
\big|\hat\Delta_k\big(\hat G_p(l)\hat D(l)\big)\big|
&\le&\sum_x(1-\cos k\cdot x)(G_p*D)(x)\nn\\
&\stackrel{\refeq{telescope}}\le&2\Big(|\hat\Delta_k\hat G_p(0)|
 +\chi_p\big(1-\hat D(k)\big)\Big),
\end{eqnarray}
and that
\begin{eqnarray}
\big|\hat\Delta_k\big(\partial_p\hat G_p(l)\hat D(l)\big)\big|&\le&\sum_x
 (1-\cos k\cdot x)(\partial_pG_p*D)(x)\nn\\
&\stackrel{\refeq{subaddBKbd}\,\&\,\refeq{telescope}}\le&2\Big(|\hat\Delta_k
 \partial_p\hat G_p(0)|+\chi_p^2\big(1-\hat D(k)\big)\Big)\nn\\
&\stackrel{\refeq{SAW-DeltaGbd4}}\le&2\Big(6\chi_p|\hat\Delta_k\hat G_p(0)|
 +4\chi_p^2\big(1-\hat D(k)\big)\Big).
\end{eqnarray}
Therefore, to evaluate $\tilde g_{3,k,l}(p)$ and
$\partial_p\tilde g_{3,k,l}(p)$, it suffices to evaluate
$|\hat\Delta_k\hat G_p(0)|$.

To evaluate $|\hat\Delta_k\hat G_p(0)|$ by using \refeq{telescope}, as we did
for SAW, we first rewrite the expression \refeq{perc2pt} for $G_p(x)$.  To
do so, we introduce ordering among self-avoiding paths from $o$ to $x$ as
follows.  For each vertex $x$, let $B(x)$ be the set of bonds incident on $x$.
Order the elements in $B(x)$ in an arbitrary but fixed manner.  For a pair of
bonds $b,b'\in B(x)$, we write $b\prec b'$ if $b$ is lower than $b'$ in that
ordering.  For a pair of self-avoiding paths $\omega,\omega'\in\Omega(x,y)$,
we write $\omega\prec\omega'$ if at the first time $\tau$ when $\omega$
becomes incompatible with $\omega'$ (therefore $b_j(\omega)=b_j(\omega')$
for all $j<\tau$) we have $b_\tau(\omega)\prec b_\tau(\omega')$.  We say that
$\omega$ is occupied if all $b_1(\omega),\dots,b_{|\omega|}(\omega)$ are
occupied.  Let $E_{x,y}(\omega)$ be the event that $\omega\in\Omega(x,y)$ is
the lowest occupied path from $x$ to $y$:
\begin{align}\lbeq{lowestoccupied}
E_{x,y}(\omega)=\{\omega\text{ is occupied}\}\setminus\bigcup_{\substack{
 \omega'\in\Omega(x,y)\\ (\omega'\prec\omega)}}\{\omega'\text{ is occupied}\}.
\end{align}
Then, we can rewrite the expression \refeq{perc2pt} for $G_p(x)$ as
\begin{align}
G_p(x)=\sum_{\omega\in\Omega(o,x)}\mP_p\big(E_{o,x}(\omega)\big).
\end{align}
Similarly to \refeq{SAW-DeltaGbd1}, we can bound $|\hat\Delta_k\hat G_p(0)|$ as
\begin{align}
|\hat\Delta_k\hat G_p(0)|
&\le\sum_x\sum_{\omega\in\Omega(o,x)}|\omega|\sum_{i=1}^{|\omega|}\big(1-\cos
 k\cdot(\omega_i-\omega_{i-1})\big)\,\mP_p\big(E_{o,x}(\omega)\big)\nn\\
&=\sum_{u,v,x}\big(1-\cos k\cdot(v-u)\big)\sum_{\omega\in\Omega(o,x)}|\omega|
 \sum_{i=1}^{|\omega|}\ind{b_i(\omega)=(u,v)}\,\mP_p\big(E_{o,x}(\omega)\big).
\end{align}
Let $\eta=(\omega_0,\dots,\omega_i)$ and
$\xi=(\omega_i,\dots,\omega_{|\omega|})$ and denote their concatenation in that
order by $\eta\circ\xi$.  Then, the above inequality is equivalent to
\begin{align}
|\hat\Delta_k\hat G_p(0)|\le\sum_{u,v,x}\big(1-\cos k\cdot(v-u)\big)\sum_{\eta
 \in\Omega(o,v)}\ind{b_{|\eta|}(\eta)=(u,v)}\nn\\
\times\sum_{\substack{\xi\in\Omega(v,x)\\ (\eta\circ\xi\in\Omega(o,x))}}
 (|\eta|+|\xi|)\,\mP_p\big(E_{o,x}(\eta\circ\xi)\big).
\end{align}

Next, we rewrite $\mP_p(E_{o,x}(\eta\circ\xi))$.  To do so, we introduce a
peculiar cluster of $\eta$ as follows.  Given a vertex $y$ and a bond
$b\in B(y)$, we define $\cC_{\prec b}(y)$ to be the set of vertices that are
connected from $y$ via an occupied bond $b'\in B(y)$ with $b'\prec b$;
if there are no such occupied bonds, then we define $\cC_{\prec b}(y)=\{y\}$.
Given a self-avoiding path $\eta$, we let
\begin{align}
\cC_{\prec\eta}=\bigcup_{j=1}^{|\eta|}\cC_{\prec b_j(\eta)}(\eta_{j-1}).
\end{align}
Notice that ther terminal point $\eta_{|\eta|}$ is not in $\cC_{\prec\eta}$.  Using this notation and
recalling \refeq{lowestoccupied}, we can rewrite the event
$E_{o,x}(\eta\circ\xi)$ for $\eta\in\Omega(o,v)$ and $\xi\in\Omega(v,x)$
with $\eta\circ\xi\in\Omega(o,x)$ as
\begin{align}
E_{o,x}(\eta\circ\xi)=E_{o,v}(\eta)\cap\big\{E_{v,x}(\xi)\text{ occurs on }\Ld
 \setminus\cC_{\prec\eta}\big\}.
\end{align}
For a $V\subset\Ld$, we let $\mP_p^V$ be the percolation measure defined by
making all bonds $b$ with $b\cap(\Ld\setminus V)\ne\vno$ vacant.  Then, we
obtain the rewrite
\begin{align}
\mP_p(E_{o,x}(\eta\circ\xi))=\mE_p\Big[\indic_{E_{o,v}(\eta)}\,\mP_p^{\Ld
 \setminus\cC_{\prec\eta}}\big(E_{v,x}(\xi)\big)\Big],
\end{align}
hence
\begin{align}\lbeq{perc-DeltaGbd1}
|\hat\Delta_k\hat G_p(0)|&\le\sum_{u,v,x}\big(1-\cos k\cdot(v-u)\big)\sum_{\eta
 \in\Omega(o,v)}\ind{b_{|\eta|}(\eta)=(u,v)}\nn\\
&\qquad\times\Bigg(|\eta|\,\mE_p\bigg[\indic_{E_{o,v}(\eta)}\sum_{\substack{\xi
 \in\Omega(v,x)\\ (\eta\circ\xi\in\Omega(o,x))}}\mP_p^{\Ld\setminus\cC_{\prec
 \eta}}\big(E_{v,x}(\xi)\big)\bigg]\nn\\
&\qquad\qquad+\mE_p\bigg[\indic_{E_{o,v}(\eta)}\sum_{\substack{\xi\in\Omega(v,
 x)\\ (\eta\circ\xi\in\Omega(o,x))}}|\xi|\,\mP_p^{\Ld\setminus\cC_{\prec\eta}}
 \big(E_{v,x}(\xi)\big)\bigg]\Bigg).
\end{align}

The contribution from the first expectation is evaluated as follows.  First,
we note that the sum over $\xi$ can be replaced by the sum over
$\xi\in\Omega(v,x)$ that are restricted in $\Ld\setminus\cC_{\prec\eta}$,
or $\mP_p^{\Ld\setminus\cC_{\prec\eta}}(E_{v,x}(\xi))=0$ otherwise.
Then, the resulting sum equals the restricted two-point function on
$\Ld\setminus\cC_{\prec\eta}$ and is bounded by the full two-point function
$G_p(x-v)$.  Therefore,
\begin{align}\lbeq{perc-DeltaGbd2}
&\sum_{\eta\in\Omega(o,v)}\ind{b_{|\eta|}(\eta)=(u,v)}|\eta|\,\mE_p\bigg[
 \indic_{E_{o,v}(\eta)}\sum_{\substack{\xi\in\Omega(v,x)\\ (\eta\circ\xi\in
 \Omega(o,x))}}\mP_p^{\Ld\setminus\cC_{\prec\eta}}\big(E_{v,x}(\xi)\big)
 \bigg]\nn\\
&\le\sum_{\eta\in\Omega(o,v)}\ind{b_{|\eta|}(\eta)=(u,v)}|\eta|\,\mP_p\big(
 E_{o,v}(\eta)\big)\,G_p(x-v).
\end{align}
We apply the same analysis to $\eta=\zeta\circ(u,v)$, where
$\zeta=(\eta_0,\dots,\eta_{|\eta|-1})$, and obtain
\begin{align}\lbeq{perc-DeltaGbd3}
\refeq{perc-DeltaGbd2}&\le pD(v-u)\,G_p(x-v)\sum_{\zeta\in\Omega(o,u)}(|\zeta|
 +1)\,\mP_p\big(E_{o,u}(\zeta)\big)\nn\\
&=pD(v-u)\,G_p(x-v)\sum_y\sum_{\substack{\zeta'\in\Omega(o,y)\\ \zeta''\in
 \Omega(y,u)\\ (\zeta'\circ\zeta''\in\Omega(o,u))}}\mP_p\big(E_{o,u}(\zeta'
 \circ\zeta'')\big),
\end{align}
where the equality is due to the identity $|\zeta|+1=\sum_y\ind{y\in\zeta}$.
Again, by the same analysis as discussed above, we finally obtain
\begin{align}\lbeq{perc-DeltaGbd4}
\refeq{perc-DeltaGbd3}\le G_p^{*2}(u)\,pD(v-u)\,G_p(x-v).
\end{align}

The contribution from the second expectation in \refeq{perc-DeltaGbd1} can be
evaluated in a similar way, and the result is
\begin{align}\lbeq{perc-DeltaGbd5}
&\sum_{\eta\in\Omega(o,v)}\ind{b_{|\eta|}(\eta)=(u,v)}\mE_p\bigg[
 \indic_{E_{o,v}(\eta)}\sum_{\substack{\xi\in\Omega(v,x)\\ (\eta\circ\xi\in
 \Omega(o,x))}}|\xi|\,\mP_p^{\Ld\setminus\cC_{\prec\eta}}\big(E_{v,x}(\xi)\big)
 \bigg]\nn\\
&\le G_p(u)\,pD(v-u)\,G_p^{*2}(x-v).
\end{align}
Substituting \refeq{perc-DeltaGbd4}--\refeq{perc-DeltaGbd5} back into
\refeq{perc-DeltaGbd1}, we obtain the same bound as
\refeq{SAW-DeltaGbd3}:
\begin{align}
|\hat\Delta_k\hat G_p(0)|\le2p_0\big(1-\hat D(k)\big)\chi_{p_0}^3,
\end{align}
which implies finiteness of $\tilde g_{3,k,l}(p)$ and
$\partial_p\tilde g_{3,k,l}(p)$ uniformly in $k,l\in\Td$ and $p\in[1,p_0]$, as
required.  This completes the proof of Proposition~\ref{prp:f-cont}.
\QED

\section{Lace-expansion analysis for self-avoiding walk}\label{s:LESAW}
In this section, we prove Propositions~\ref{prp:f-initial}--\ref{prp:lace-gen}
for SAW.  First, in Section~\ref{ss:derivation-SAW}, we explain the derivation
of the lace expansion, Proposition~\ref{prp:lace-gen}, for SAW.
In Section~\ref{ss:Pibds-saw}, we prove bounds on the lace-expansion
coefficients in terms of basic diagrams, such as $L_p$ and $B_p$.  Finally,
in Section~\ref{ss:LBWbd}, we prove bounds on those basic diagrams in terms of
RW loops and RW bubbles and use them to prove 
Propositions~\ref{prp:f-initial}--\ref{prp:f-bootstrapping} on
$\mL^{\!d\ge6}$.  We close this section by addressing potential elements for
extending the result to 5 dimensions, in Section~\ref{ss:discussion-saw}.

\subsection{Derivation of the lace expansion}\label{ss:derivation-SAW}
Proposition~\ref{prp:lace-gen} for SAW is restated as follows.

\begin{prp}[Lace expansion for SAW]\label{prp:LESAW}
For any $p<\pc$ and $N\in\mN$, there are nonnegative
functions $\{\pi_p^{\sss(n)}\}_{n=1}^N$ on $\Ld$ such that, if we define
$\Pi_p^{\sss(N)}$ as
\begin{align}
\Pi_p^{\sss(N)}(x)=\sum_{n=1}^N(-1)^n\pi_p^{\sss(n)}(x),
\end{align}
then we obtain the recursion equation
\begin{align}\lbeq{lesaw}
G_p(x)=\delta_{o,x}+\big((pD+\Pi_p^{\sss(N)})*G_p\big)(x)+
 (-1)^{N+1}R_p^{\sss(N+1)}(x),
\end{align}
where the remainder $R_p^{\sss(N)}$ obeys the bound
\begin{align}\lbeq{remsaw}
0\le R_p^{\sss(N)}(x)\le(\pi_p^{\sss(N)}*G_p)(x).
\end{align}
\end{prp}

\Proof{Sketch proof of Proposition~\ref{prp:LESAW}.}
First, we derive the first expansion, i.e., \refeq{lesaw} for $N=1$.  
For notational convenience, we use
\begin{align}
P(\omega)=p^{|\omega|}\prod_{j=1}^{|\omega|}D(\omega_j -\omega_{j-1}).
\end{align}
Then, by splitting the sum in \refeq{Gsaw-def} into two depending on whether
$|\omega|$ is zero or positive, we obtain
\begin{align}\lbeq{lesaw-1st}
G_p(x)=\delta_{o,x}+\sum_{\substack{\omega\in\Omega(o,x)\\ (|\omega|\ge1)}}
 P(\omega)=\delta_{o,x}+\sum_ypD(y)\sum_{\omega\in\Omega(y,x)}P(\omega)\ind{o
 \notin\omega}.
\end{align}
This is depicted as
\begin{align}
\includegraphics[scale=0.5]{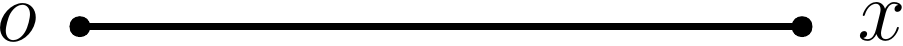}~=\delta_{o,x}
 +~\raisebox{-4pt}{\includegraphics[scale=0.5]{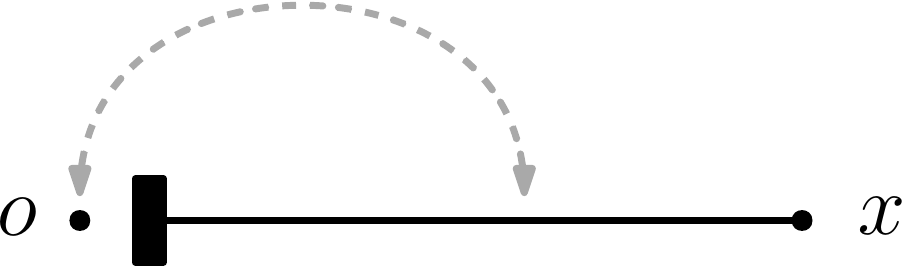}}
\end{align}
where the rectangle next to the origin represents that there is a bond from
$o$ to a neighboring vertex $y$, which is summed over $\Ld$ and unlabeled
in the picture, and the dashed two-sided arrow represents mutual avoidance
between $o$ and SAWs from $y$ to $x$, which corresponds to the indicator
$\ind{o\notin\omega}$ in \refeq{lesaw-1st}.  Using the identity 
$\ind{o\notin\omega}=1-\ind{o\in\omega}$ due to the inclusion-exclusion
relation, we complete the first expansion as
\begin{align}
G_p(x)&=\delta_{o,x}+~\raisebox{-4pt}{\includegraphics[scale=0.5]{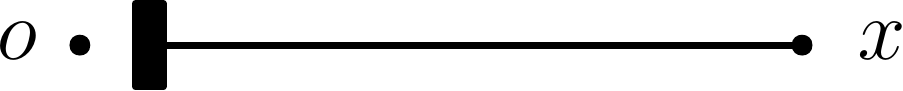}}~-~
 \underbrace{\raisebox{-2.9pc}{\includegraphics[scale=0.5]{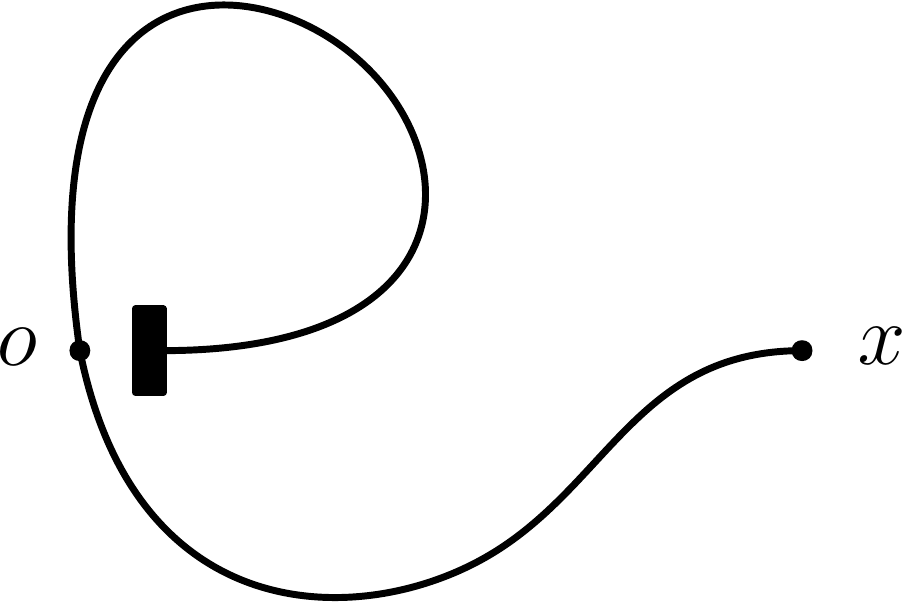}}}_{\equiv
 R_p^{\sss(1)}(x)}\nn\\
&=\delta_{o,x}+(pD*G_p)(x)-R_p^{\sss(1)}(x).
\end{align}

Next, we expand the remainder $R_p^{\sss(1)}(x)$ to complete the first 
expansion.  Splitting each SAW from $y$ (summed over $\Ld$ and unlabeled
in the picture) to $x$ through $o$ into two SAWs, $\omega_1\in\Omega(y,o)$
and $\omega_2\in\Omega(o,x)$ (in red), we can rewrite $R_p^{\sss(1)}(x)$ as
\begin{align}
R_p^{\sss(1)}(x)=~\raisebox{-2.9pc}{\includegraphics[scale=0.5]{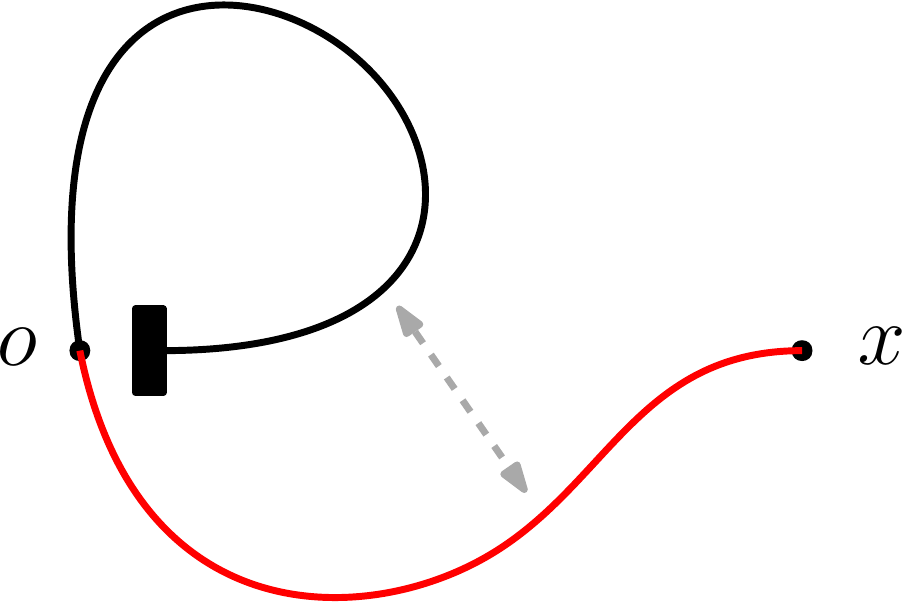}}
\end{align}
where the dashed two-sided arrow implies that the concatenation of $\omega_1$
and $\omega_2$ in this order, denoted $\omega_1\circ\omega_2$, is SAW.  Using
the identity $\ind{\omega_1\circ\omega_2\text{ is SAW}}
=1-\ind{\omega_1\circ\omega_2\text{ is not SAW}}$, we obtain
\begin{align}
R_p^{\sss(1)}(x)&=~\raisebox{-2.8pc}{\includegraphics[scale=0.5]{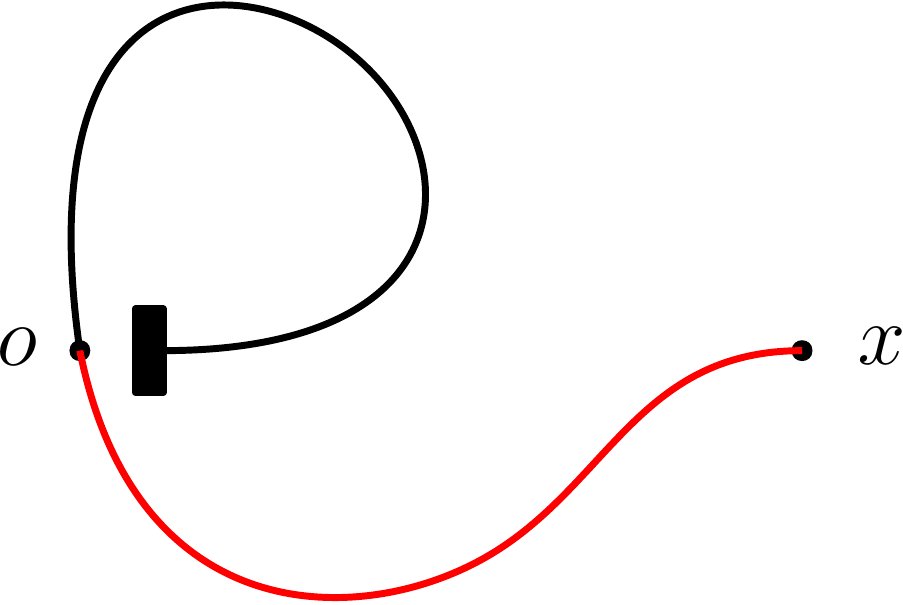}}~-~
 \underbrace{\raisebox{-3.1pc}{\includegraphics[scale=0.5]{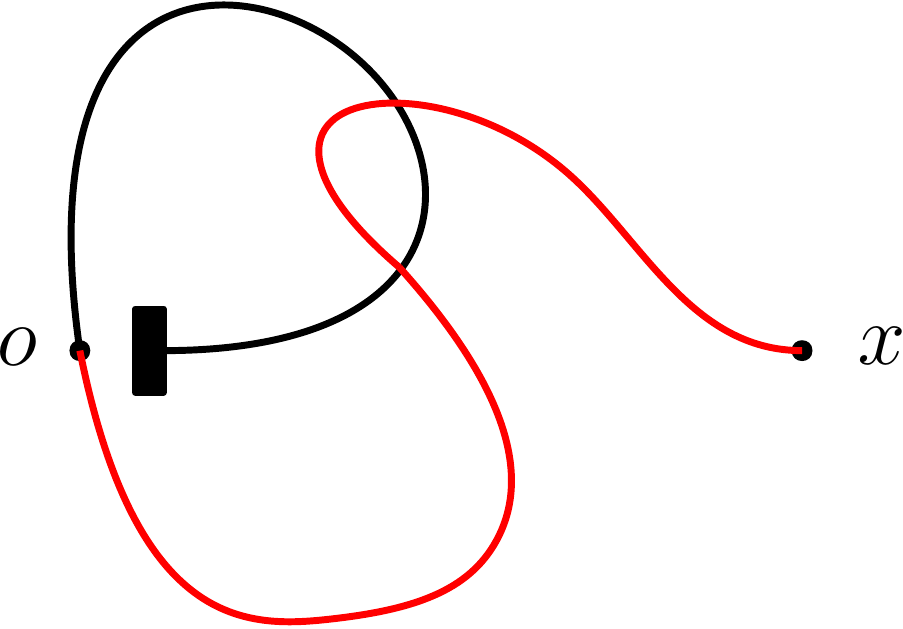}}}_{\equiv
 R_p^{\sss(2)}(x)}\nn\\
&=\sum_y\underbrace{\pisawone{o=y}}_{\equiv\pi_p^{\sss(1)}(y)}G_p(x-y)
 -R_p^{\sss(2)}(x),
\end{align}
where the precise definition of $\pi_p^{\sss(1)}(x)$ is the following:
\begin{align}
\pi_p^{\sss(1)}(x)=(pD*G_p)(o)\,\delta_{o,x}.
\end{align}
Since $R_p^{\sss(2)}(x)$ is nonnegative, this also implies \refeq{remsaw} for
$N=1$.  This completes the first expansion.

To show how to derive the higher-order expansion coefficients, we further
demonstrate the expansion of the remainder $R_p^{\sss(2)}(x)$.  Since
$\omega_1\circ\omega_2$ is not SAW, there must be at least one vertex other
than $o$ where $\omega_2$ hits $\omega_1$.  Take the first such vertex, say,
$z\ne o$, which is summed over $\Ld$ and unlabeled in the following picture,
and split $\omega_2\in\Omega(o,x)$ into two SAWs, $\omega_{21}\in\Omega(o,z)$
and $\omega_{22}\in\Omega(z,x)$ (in blue), so that
$\omega_1\cap\omega_{21}=\{o,z\}$.  Then, we can rewrite $R_p^{\sss(2)}(x)$ as
\begin{align}
R_p^{\sss(2)}(x)=~\raisebox{-3.1pc}{\includegraphics[scale=0.5]{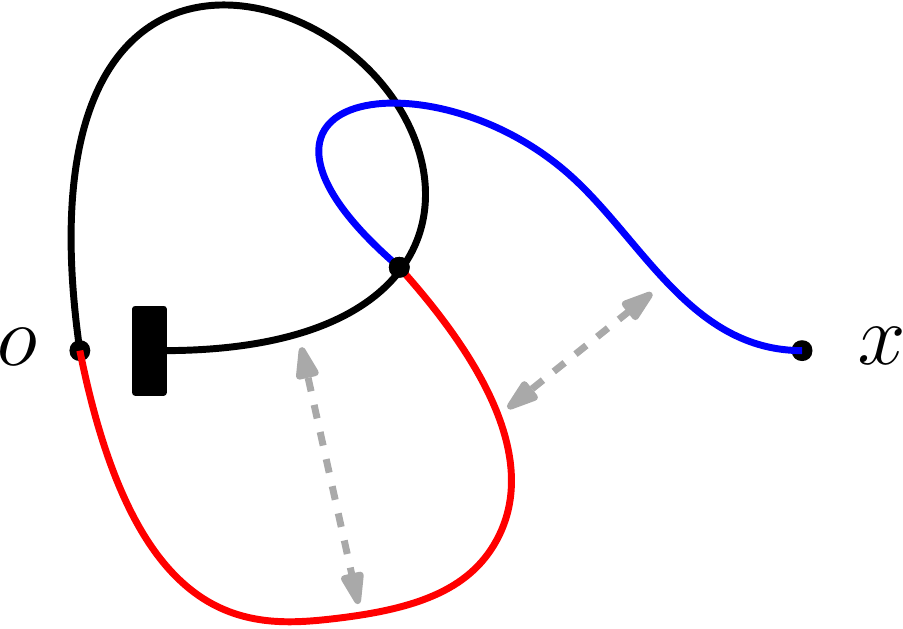}}
\end{align}
where the dashed two-sided arrow between the red $\omega_{21}$ and the blue
$\omega_{22}$ implies that the concatenation $\omega_{21}\circ\omega_{22}$ is
SAW.  Using the identity $\ind{\omega_{21}\circ\omega_{22}
\text{ is SAW}}=1-\ind{\omega_{21}\circ\omega_{22}\text{ is not SAW}}$,
we obtain
\begin{align}
R_p^{\sss(2)}(x)&=~\raisebox{-3.1pc}{\includegraphics[scale=0.5]{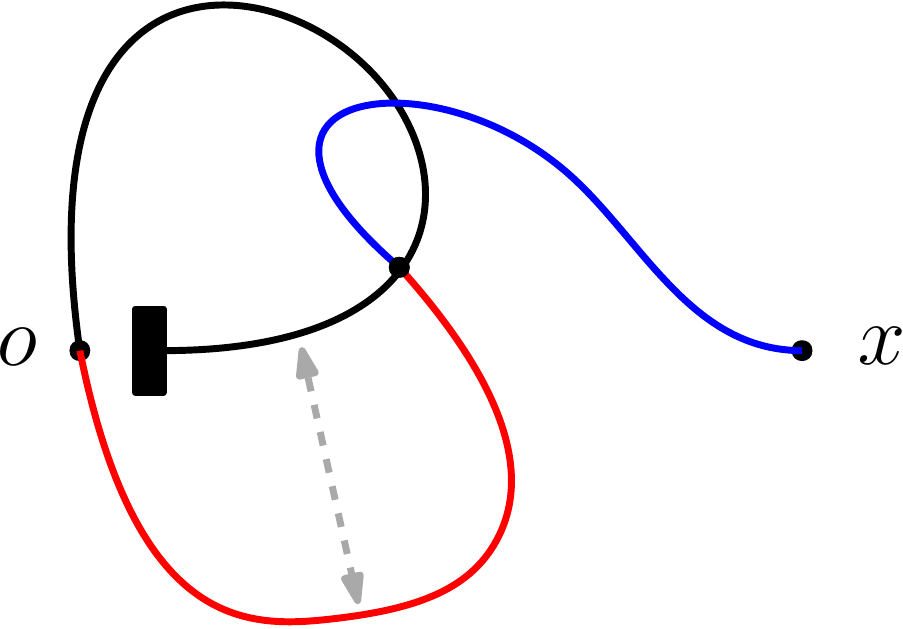}}~-~
 \underbrace{\raisebox{-3.1pc}{\includegraphics[scale=0.5]{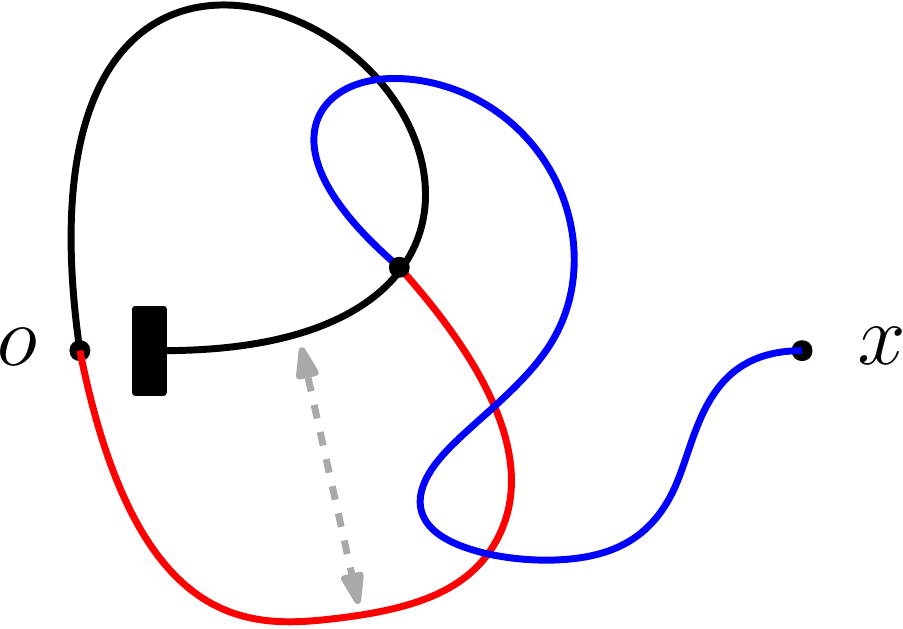}}}_{\equiv
 R_p^{\sss(3)}(x)}\nn\\
&=\sum_y\underbrace{\pisawtwo{o}{y}}_{\equiv\pi_p^{\sss(2)}(y)}G_p(x-y)
 -R_p^{\sss(3)}(x),
\end{align}
where the precise definition of $\pi_p^{\sss(2)}(x)$ is the following:
\begin{align}\lbeq{pisaw2-def}
\pi_p^{\sss(2)}(x)&=(1-\delta_{o,x})\sum_{\omega_1,\omega_2,\omega_3\in\Omega
 (o,x)}P(\omega_1)P(\omega_2)P(\omega_3)\prod_{i\ne j}\ind{\omega_i\cap\omega_j
 =\{o,x\}}.
\end{align}
Since $R_p^{\sss(3)}(x)$ is nonnegative, this implies \refeq{remsaw} for
$N=2$, as required.

By repeated application of inclusion-exclusion relations, we obtain the lace
expansion \refeq{lesaw}, with the lace-expansion coefficients depicted as
\begin{align}\lbeq{pisaw-diagr}
\pi_p^{\sss(3)}(x)=\pisawthree{o}{x},&&
\pi_p^{\sss(4)}(x)=\pisawfour{o}{x}{},&&
\pi_p^{\sss(5)}(x)=\pisawfive{o}{x}{},~\dots
\end{align}
where the slashed line segments represent SAWs with length $\ge0$, while the
others represent SAWs with length $\ge1$.  The unlabeled vertices are summed
over $\Ld$.  Due to the construction explained above, the red line segments
avoid the black ones, the blue ones avoid the red ones, the yellow ones avoid
the blue ones, and so on.  We complete the sketch proof of
Proposition~\ref{prp:LESAW}.
\QED

\subsection{Diagrammatic bounds on the expansion coefficients}
 \label{ss:Pibds-saw}
As explained in \emph{Step~1} in Section~\ref{ss:whereandhow}, the
bootstrapping functions $\{g_i(p)\}_{i=1}^3$ are bounded in terms of sums of
$\hat\pi_p^{\sss(n)}(0)$ and $|\hat\Delta_k\hat\pi_p^{\sss(n)}(0)|$.  In this
subsection, we prove bounds on those quantities in terms of basic diagrams,
such as $L_p$ and $B_p$ in \refeq{LBT-def}, as briefly explained in
\emph{Step~2} in Section~\ref{ss:whereandhow}.  Recall that
\begin{align}\lbeq{LBr-redef}
L_p=\|(pD)^{*2}*G_p\|_\infty,&&
B_p=\|(pD)^{*2}*G_p^{*2}\|_\infty,&&
r=p\|D\|_\infty+L_p+B_p.
\end{align}
We also define
\begin{align}\lbeq{B'W-def}
B_p'=\|(pD)^{*4}*G_p^{*2}\|_\infty,&&
\hat W_p(k)=\sup_x(1-\cos k\cdot x)G_p(x).
\end{align}

\begin{lmm}[Diagrammatic bounds on the expansion coefficients]\label{lmm:DBSAW}
The expansion coefficients
$\hat\pi_p^{\sss(n)}(0)\equiv\sum_x\pi_p^{\sss(n)}(x)$ and $|\hat\Delta_k
\hat\pi_p^{\sss(n)}(0)|\equiv\sum_x(1-\cos k\cdot x)\pi_p^{\sss(n)}(x)$, both
nonnegative, obey the following bounds:
\begin{align}
\hat\pi_p^{\sss(n)}(0)&\le
 \begin{cases}
 L_p&[n=1],\\
 B_p(p\|D\|_\infty+L_p)r^{n-2}&[n\ge2],
 \end{cases}\lbeq{pisaw-sumbd}\\
|\hat\Delta_k\hat\pi_p^{\sss(n)}(0)|&\le
 \begin{cases}
 B_p^2\hat W_p(k)m^2r^{2m-2}&[n=2m+1],\\
 B_p^2\hat W_p(k)m(m-1)r^{2m-3}+B_p\hat W_p(k)mr^{2m-2}&[n=2m].
 \end{cases}\lbeq{pisaw-diffbd}
\end{align}
For $|\hat\Delta_k\hat\pi_p^{\sss(2)}(0)|$, in particular,
the following bound also holds:
\begin{align}\lbeq{imp-2-bd}
|\hat\Delta_k\hat\pi_p^{\sss(2)}(0)|\le\frac{3B_p}{2^d}p\big(1-\hat
 D(k)\big)+B_p'\hat W_p(k).
\end{align}
\end{lmm}

\paragraph{Remark.}
As shown in \refeq{W1bd-saw} and \refeq{Wpbd-saw} in the next subsection,
$\|\hat W_p/(1-\hat D)\|_\infty$ could be relatively large, compared to
$L_p,B_p$ and $r$.  Therefore, if we want to have a good bound on
$|\hat\Delta_k\hat\pi_p^{\sss(n)}(0)|/(1-\hat D(k))$, we should have a small
multiplicative factor to $\hat W_p(k)$.  By \refeq{pisaw-diffbd}, that
multiplicative factor is at most $\lfloor\frac{n}2\rfloor^2B_pr^{n-2}$ for
$n\ge2$ (it is zero for $n=1$, due to the definition of $\pi_p^{\sss(1)}$)
and the dominant contribution comes from the case of $n=2$, i.e., $B_p$.
In \refeq{imp-2-bd}, on the other hand, the multiplicative factor to
$\hat W_p(k)$ is $B'_p$, which is potentially much smaller than $B_p$.
This can be seen by comparing the RW versions of $B_p$ and $B'_p$,
which are the RW bubble $\vep_2$ and
\begin{align}
\vep'_2=(D^{*4}*S_1^2)(o)=\sum_{n=2}^\infty(2n-3)D^{*2n}(o).
\end{align}
Table~\ref{table:vep'2} summarizes the bounds on those RW bubbles
that are evaluated as explained in Section~\ref{ss:bcc}.
\begin{table}[b]
\caption{\label{table:vep'2}Comparison of upper bounds on the RW bubbles
for $4\le d\le 9$.}
\begin{center}
\begin{tabular}{r|llllll}
          &    $d=4$ & $d=5$ & $d=6$ & $d=7$ & $d=8$ & $d=9$\\
\hline
$\vep_2$  & $\infty$ & 0.178332 & 0.044004 & 0.015302 & 0.006156 & 0.002678\\
$\vep_2'$ & $\infty$ & 0.115931 & 0.018708 & 0.004302 & 0.001161 & 0.000344\\
\end{tabular}
\end{center}
\end{table}

The amount of extra work caused by the use of \refeq{imp-2-bd} instead of using
only \refeq{pisaw-diffbd} is quite small.  However, this is the key to be able
to go down to 6 dimensions.  We will get back to this point in
Section~\ref{ss:discussion-saw}.

\bigskip

\Proof{Sketch proof of Lemma~\ref{lmm:DBSAW}.}
In the following, we repeatedly use the trivial inequality
\begin{align}\lbeq{sbad-eq}
G_p(x)\ind{x\ne o}\le(pD*G_p)(x).
\end{align}
For example,
\begin{align}
\hat\pi_p^{\sss(1)}(0)=\sum_{x\ne o}pD(x)\,G_p(x)\le\big((pD)^{*2}*G_p\big)(o)
 \le L_p.
\end{align}
For $n\ge2$, we first decompose $\hat\pi_p^{\sss(n)}(0)$ by using
subadditivity and then repeatedly apply \refeq{sbad-eq} to obtain
\refeq{pisaw-sumbd}.  For example,
\begin{eqnarray}\lbeq{pisaw2bd}
\hat\pi_p^{\sss(2)}(0)=\pisawtwo{o}{}
&\le&\bigg(\sum_{x\ne o}G_p(x)^2\bigg)\bigg(\sup_{x\ne o}G_p(x)\bigg)\nn\\
&\stackrel{\refeq{sbad-eq}}\le&\underbrace{\big((pD)^{*2}*G_p^{*2}\big)(o)}_{\le
 B_p}\bigg(\sup_{x\ne o}(pD*G_p)(x)\bigg),
\end{eqnarray}
and
\begin{eqnarray}\lbeq{pisaw5bd}
\hat\pi_p^{\sss(5)}(0)&=&\pisawfive{o}{}{}\nn\\
&\le&\bigg(\sum_{x\ne o}G_p(x)^2\bigg)\bigg(\sup_{x\ne o}\sum_{y\ne x}G_p(y)
 G_p(x-y)\bigg)^3\bigg(\sup_{x\ne o}G_p(x)\bigg)\nn\\
&\stackrel{\refeq{sbad-eq}}\le&B_p\bigg(\sup_x(pD*G_p^{*2})(x)\bigg)^3
 \bigg(\sup_x(pD*G_p)(x)\bigg).
\end{eqnarray}
In general, $\hat\pi_p^{\sss(n)}(0)$ for $n\ge2$ is bounded by the right-most
expression with the power 3 replaced by $n-2$.
Notice that, by omitting the spatial variables, we have
\begin{align}
pD*G_p=pD*\big(\delta+(1-\delta)G_p\big)\stackrel{\refeq{sbad-eq}}\le
 pD+(pD)^{*2}*G_p,
\end{align}
where $\delta$ is the Kronecker delta, hence
\begin{align}
\sup_x(pD*G_p)(x)&\le p\|D\|_\infty+\|(pD)^{*2}*G_p\|_\infty
 =p\|D\|_\infty+L_p.
\end{align}
Similarly, we have
\begin{eqnarray}
pD*G_p^{*2}&=&pD*G_p*\big(\delta+(1-\delta)G_p\big)\nn\\
&\stackrel{\refeq{sbad-eq}}\le&pD*G_p+(pD)^{*2}*G_p^{*2}\nn\\
&=&pD*\big(\delta+(1-\delta)G_p\big)+(pD)^{*2}*G_p^{*2}\nn\\
&\stackrel{\refeq{sbad-eq}}\le&pD+(pD)^{*2}*G_p+(pD)^{*2}*G_p^{*2},
\end{eqnarray}
hence
\begin{align}\lbeq{boundbyr}
\sup_x(pD*G_p^{*2})(x)&\le p\|D\|_\infty+\|(pD)^{*2}*G_p\|_\infty
 +\|(pD)^{*2}*G_p^{*2}\|_\infty\nn\\
&=p\|D\|_\infty+L_p+B_p\equiv r.
\end{align}
This completes the proof of \refeq{pisaw-sumbd}.

Next, we prove \refeq{pisaw-diffbd} for $n=2m+1$.  Since $\pi_p^{\sss(1)}(x)$
is proportional to $\delta_{o,x}$ and therefore
$\hat\Delta_k\hat\pi_p^{\sss(1)}(0)\equiv0$, we can assume $m\ge1$.  To bound
$|\hat\Delta_k\hat\pi_p^{\sss(2m+1)}(0)|\equiv\sum_x(1-\cos k\cdot x)
\pi_p^{\sss(2m+1)}(x)$ for $m\ge1$, we first identify the diagram vertices
along the lowest diagram path from $o$ to $x$, say, $y_1,\dots,y_{m-1}$, and
then split $x$ into $\{y_j-y_{j-1}\}_{j=1}^m$, where $y_0=o$ and $y_m=x$.
For example,
\begin{align}
|\hat\Delta_k\hat\pi_p^{\sss(5)}(0)|=\sum_{y_1,y_2}\bigg(1-\cos\sum_{j=1,2}k
 \cdot(y_j-y_{j-1})\bigg)\pisawfive{y_0=o}{y_2}{~y_1}
\end{align}
Then, by using \refeq{telescope} and subadditivity, we obtain
\begin{align}
|\hat\Delta_k\hat\pi_p^{\sss(5)}(0)|&\le2\sum_{y_1,y_2}\Big((1-\cos k\cdot y_1)
 +\big(1-\cos k\cdot(y_2-y_1)\big)\Big)\nn\\
&\hskip4pc\times\bigg(G_p(y_1)\pisawfivew{o}{y_2}{y_1}+~G_p(y_2-y_1)
 \pisawfiveww{o}{y_2}{~y_1}\bigg)\nn\\
&\le2\hat W_p(k)\bigg(\pisawfivew{o}{}{}~+\pisawfiveww{o}{}{}\bigg).
\end{align}
Each remaining diagram is bounded, by following similar decomposition to
\refeq{pisaw2bd}--\refeq{pisaw5bd} and then using \refeq{boundbyr}, by
$B_p^2r^2$, yielding the desired bound on
$|\hat\Delta_k\hat\pi_p^{\sss(5)}(0)|$.  In general,
\begin{align}
|\hat\Delta_k\hat\pi_p^{\sss(2m+1)}(0)|&\le m\hat W_p(k)\times\Big(m
 \text{ diagrams, each bounded by }B_p^2r^{2m-2}\Big)\nn\\
&\le B_p^2\hat W_p(k)m^2r^{2m-2},
\end{align}
as required.

To prove \refeq{pisaw-diffbd} for $n=2m$, we follow the same line as above for
$n=2m+1$.  To bound $|\hat\Delta_k\hat\pi_p^{\sss(2m)}(0)|\equiv\sum_x(1-\cos
 k\cdot x)\pi_p^{\sss(2m)}(x)$, we first identify the diagram vertices along
the lowest diagram path from $o$ to $x$, say, $y_1,\dots,y_{m-1}$, and then
split $x$ into $\{y_j-y_{j-1}\}_{j=1}^m$, where $y_0=o$ and $y_m=x$.  For
example,
\begin{align}
|\hat\Delta_k\hat\pi_p^{\sss(4)}(0)|=\sum_{y_1,y_2}\bigg(1-\cos\sum_{j=1,2}k
 \cdot(y_j-y_{j-1})\bigg)\pisawfour{y_0=o}{y_2}{\quad~y_1}
\end{align}
Then, by using \refeq{telescope} and subadditivity, we obtain
\begin{align}
|\hat\Delta_k\hat\pi_p^{\sss(4)}(0)|&\le2\sum_{y_1,y_2}\Big((1-\cos k\cdot y_1)
 +\big(1-\cos k\cdot(y_2-y_1)\big)\Big)\nn\\
&\hskip4pc\times\bigg(G_p(y_1)\pisawfourw{o}{y_2}{\quad~y_1}+~G_p(y_2-y_1)
 \pisawfourww{o}{\;y_2}{\quad~~y_1}\bigg)\nn\\
&\le2\hat W_p(k)\bigg(\pisawfourw{o}{}{}~+\pisawfourww{o}{}{}~\bigg).
\end{align}
Following similar decomposition to \refeq{pisaw2bd}--\refeq{pisaw5bd} and
using \refeq{boundbyr}, we can bound the first diagram by $B_p^2r$, while
the second diagram is bounded by $B_pr^2$, yielding the desired bound on
$|\hat\Delta_k\hat\pi_p^{\sss(4)}(0)|$.  In general,
\begin{align}
|\hat\Delta_k\hat\pi_p^{\sss(2m)}(0)|&\le m\hat W_p(k)\times\bigg(\Big((m-1)
 \text{ diagrams, each bounded by }B_p^2r^{2m-3}\Big)\nn\\
&\hskip7pc+\Big(1\text{ diagram, bounded by }B_pr^{2m-2}\Big)\bigg)\nn\\
&\le B_p^2\hat W_p(k)m(m-1)r^{2m-3}+B_p\hat W_p(k)mr^{2m-2},
\end{align}
as required.

To prove the bound \refeq{imp-2-bd} on $|\hat\Delta_k\hat\pi_p^{\sss(2)}(0)|$,
we recall the definition \refeq{pisaw2-def} and divide $\pi_p^{\sss(2)}(x)$
into $\pi_p^{\sss(2),=1}(x)$ and $\pi_p^{\sss(2),\ge2}(x)$, where
\begin{align}
\pi_p^{\sss(2),=1}(x)&=(1-\delta_{o,x})\sum_{\substack{\omega_1,\omega_2,
 \omega_3\in\Omega(o,x)\\ (\exists i:|\omega_i|=1)}}P(\omega_1)P(\omega_2)
 P(\omega_3)\prod_{i\ne j}\ind{\omega_i\cap\omega_j=\{o,x\}},\\
\pi_p^{\sss(2),\ge2}(x)&=(1-\delta_{o,x})\sum_{\substack{\omega_1,\omega_2,
 \omega_3\in\Omega(o,x)\\ (\forall i:|\omega_i|\ge2)}}P(\omega_1)P(\omega_2)
 P(\omega_3)\prod_{i\ne j}\ind{\omega_i\cap\omega_j=\{o,x\}}.
\end{align}
Then, by symmetry, the contribution from $\pi_p^{\sss(2),=1}(x)$ is bounded as
\begin{align}\lbeq{symm-improve}
|\hat\Delta_k\hat\pi_p^{\sss(2),=1}(0)|&\le3\sum_{x\sim o}(1-\cos k\cdot x)
 pD(x)\bigg(\underbrace{\sum_{\omega\in\Omega(o,x)}P(\omega)}_{\le(pD*G_p)
 (x)}\bigg)^2\nn\\
&\le3\bigg(\sup_{x\sim o}(pD*G_p)(x)^2\bigg)p\sum_x(1-\cos k\cdot x)D(x)\nn\\
&=3\bigg(\frac1{2^d}\underbrace{\sum_{x\sim o}(pD*G_p)(x)^2}_{\le B_p}\bigg)
 p\big(1-\hat D(k)\big),
\end{align}
while the contribution from $\pi_p^{\sss(2),\ge2}(x)$ is easily bounded as
\begin{align}
|\hat\Delta_k\hat\pi_p^{\sss(2),\ge2}(0)|&\le\sum_x(1-\cos k\cdot x)\bigg(
 \sum_{\substack{\omega\in\Omega(o,x)\\ (|\omega|\ge2)}}P(\omega)\bigg)^3\nn\\
&\le\underbrace{\sum_x\big((pD)^{*2}*G_p\big)(x)^2}_{\le B'_p}\bigg(
 \underbrace{\sup_x(1-\cos k\cdot x)G_p(x)}_{=\hat W_p(k)}\bigg),
\end{align}
This completes the proof of Lemma~\ref{lmm:DBSAW}.
\QED

\subsection{Diagrammatic bounds on the bootstrapping functions}
Let
\begin{align}
\hat\Pi_p\odd(k)=\sum_{m=0}^\infty\hat\pi_p\ssc{2m+1}(k),&&
\hat\Pi_p\even(k)=\sum_{m=1}^\infty\hat\pi_p\ssc{2m}(k).
\end{align}
Suppose that $r\equiv p\|D\|_\infty+L_p+B_p<1$.  Then, 
by Lemma~\ref{lmm:DBSAW}, we obtain 
\begin{gather}
0\le\hat\Pi_p\odd(0)\le L_p+B_p(p\|D\|_\infty+L_p)\frac{r}{1-r^2},
 \lbeq{odd-bd}\\
0\le\hat\Pi_p\even(0)\le B_p(p\|D\|_\infty+L_p)\frac1{1-r^2},\lbeq{even-bd}\\
\sup_k\frac{|\hat\Delta_k\hat\Pi_p\odd(0)|}{1-\hat D(k)}\le\frac{B_p^2(1+r^2)}
 {(1-r^2)^3}\bigg\|\frac{\hat W_p}{1-\hat D}\bigg\|_\infty,
 \lbeq{sum-diff-bd-odd}\\
\sup_k\frac{|\hat\Delta_k\hat\Pi_p\even(0)|}{1-\hat D(k)}\le\frac{3B_p}{2^d}
 p+\bigg(B'_p+B_p^2\frac{2r}{(1-r^2)^3}+B_p\frac{r^2(2-r^2)}{(1-r^2)^2}\bigg)
 \bigg\|\frac{\hat W_p}{1-\hat D}\bigg\|_\infty.\lbeq{sum-diff-bd-even}
\end{gather}
Applying these bounds to \refeq{g1bd-byPi}, \refeq{g2bd-byPi} and 
\refeq{g3bd-byPi}, we obtain the following bounds on the bootstrapping 
functions $\{g_i(p)\}_{i=1}^3$.

\begin{lmm}\label{lmm:gbds-saw}
Suppose $r<1$ and that $L_p,B_p,B'_p,\|\hat W_p/(1-\hat D)\|_\infty$ are so 
small that the two inequalities in \refeq{suffcond-saw} hold.  Then, we have
\begin{align}
g_1(p)&\le1+L_p+\frac{B_p(p\|D\|_\infty+L_p)r}{1-r^2},\lbeq{g1bd-byBWr}\\
g_2(p)&\le\bigg(1-\frac{B_p^2(1+r^2)}{(1-r^2)^3}\bigg\|\frac{\hat W_p}{1-\hat
 D}\bigg\|_\infty\bigg)^{-1},\lbeq{g2bd-byBWr}\\
g_3(p)&\le\max\{g_2(p),1\}^3\nn\\
&\quad\times\Bigg(\bigg(1+\frac{3B_p}{2^d}\bigg)p+\bigg(B'_p
 +\frac{B_p^2}{(1-r^2)(1-r)^2}+\frac{B_pr^2(2-r^2)}{(1-r^2)^2}\bigg)\bigg\|\frac{\hat W_p}{1-\hat D}
 \bigg\|_\infty\Bigg)^2.\lbeq{g3bd-byBWr}
\end{align}
\end{lmm}

\Proof{Proof.}
The bounds on $g_1(p)$ and $g_2(p)$ are easy; since
$\hat\Pi_p(0)=\hat\Pi_p\even(0)-\hat\Pi_p\odd(0)$ and 
$-\hat\Delta_k\hat\Pi_p(0)=|\hat\Delta_k\hat\Pi_p\even(0)|
-|\hat\Delta_k\hat\Pi_p\odd(0)|$, we obtain
\begin{align}
g_1(p)&\stackrel{\refeq{g1bd-byPi}}\le1+\hat\Pi_p\odd(0)
 \stackrel{\refeq{odd-bd}}\le1+L_p+\frac{B_p(p\|D\|_\infty+L_p)r}{1-r^2},\\
g_2(p)&\stackrel{\refeq{g2bd-byPi}}\le\sup_k\bigg(1-\frac{|\hat\Delta_k\hat
 \Pi_p\odd(0)|}{1-\hat D(k)}\bigg)^{-1}\stackrel{\refeq{sum-diff-bd-odd}}
 \le\bigg(1-\frac{B_p^2(1+r^2)}{(1-r^2)^3}\bigg\|\frac{\hat W_p}{1-\hat D}
 \bigg\|_\infty\bigg)^{-1}.
\end{align}
For $g_3(p)$, since $\hat G_p(k)=\hat A_p(k)\equiv1/(1-\hat J_p(k))$ for SAW
and $|\hat G_p(k)|\le g_2(p)\hat S_1(k)\equiv g_2(p)/(1-\hat D(k))$, we obtain
\begin{align}\lbeq{g3bd-prebyBWr}
g_3(p)&\stackrel{\refeq{g3bd-byPi}}\le\sup_{k,l}\frac{1-\hat D(k)}{\hat
 U(k,l)}\bigg(\frac{\hat S_1(l+k)+\hat S_1(l-k)}2\hat S_1(l)g_2(p)^2
 \frac{|\hat\Delta_k\hat J_p(l)|}{1-\hat D(k)}\nn\\
&\qquad+4\hat S_1(l+k)\hat S_1(l-k)g_2(p)^3\frac{-\hat\Delta_l
 \widehat{|J_p|}(0)}{1-\hat D(l)}\,\frac{-\hat\Delta_k\widehat{|J_p|}(0)}
 {1-\hat D(k)}\bigg)\nn\\
&\stackrel{\refeq{U-def}}\le\max\{g_2(p),1\}^3\max\bigg\{\sup_{k,l}
 \frac{|\hat\Delta_k\hat J_p(l)|}{1-\hat D(k)},~\bigg(\sup_k\frac{-\hat
 \Delta_k\widehat{|J_p|}(0)}{1-\hat D(k)}\bigg)^2~\bigg\}.
\end{align}
Since $J_p=pD+\Pi_p$ for SAW, we have
\begin{align}
\frac{|\hat\Delta_k\hat J_p(l)|}{1-\hat D(k)}&=\frac1{1-\hat D(k)}\bigg|\sum_x
 (1-\cos k\cdot x)e^{il\cdot x}\big(pD(x)+\Pi_p(x)\big)\bigg|\nn\\
&\le\frac1{1-\hat D(k)}\sum_x(1-\cos k\cdot x)\big(pD(x)+\Pi_p\odd(x)+\Pi_p
 \even(x)\big)\nn\\
&\le p+\frac{|\hat\Delta_k\hat\Pi_p\even(0)|}{1-\hat D(k)}+\frac{|\hat\Delta_k
 \hat\Pi_p\odd(0)|}{1-\hat D(k)},
\end{align}
which is larger than 1, since $p\ge1$.  It is easy to check that
$-\hat\Delta_k\widehat{|J_p|}(0)/(1-\hat D(k))$ obeys the same bound.  
Therefore, by using \refeq{sum-diff-bd-odd}--\refeq{sum-diff-bd-even}, 
we obtain
\begin{align}
g_3(p)&\le\max\{g_2(p),1\}^3\bigg(p+\sup_k\frac{|\hat\Delta_k\hat\Pi_p
 \even(0)|}{1-\hat D(k)}+\sup_k\frac{|\hat\Delta_k\hat\Pi_p\odd(0)|}{1-\hat
 D(k)}\bigg)^2\nn\\
&\le\max\{g_2(p),1\}^3\nn\\
&\quad\times\Bigg(\bigg(1+\frac{3B_p}{2^d}\bigg)p+\bigg(B'_p+\frac{B_p^2}
 {(1-r^2)(1-r)^2}+\frac{B_pr^2(2-r^2)}{(1-r^2)^2}\bigg)\bigg\|\frac{\hat W_p}
 {1-\hat D}\bigg\|_\infty\Bigg)^2,
\end{align}
as required.
\QED

\subsection{Bounds on diagrams in terms of random-walk quantities}
 \label{ss:LBWbd}
In this subsection, we evaluate the diagrams 
for $p\in[1,\pc)$ and complete the proof of 
Propositions~\ref{prp:f-initial}--\ref{prp:f-bootstrapping}.  

First, we evaluate the diagrams for $p\in(1,\pc)$ under the bootstrapping 
assumptions.

\begin{lmm}\label{lmm:LpBpWpbd-saw}
Let $d\ge5$ and $p\in(1,\pc)$ and suppose that $g_i(p)\le K_i$, $i=1,2,3$,
for some constants $\{K_i\}_{i=1}^3$.  Then, we have
\begin{align}\lbeq{LpBpB'pbd-saw}
L_p\le K_1^2K_2\vep_1,&&
B_p\le K_1^2K_2^2\vep_2,&&
B'_p\le K_1^4K_2^2\vep'_2,
\end{align}
\begin{align}\lbeq{Wpbd-saw}
\bigg\|\frac{\hat W_p}{1-\hat D}\bigg\|_\infty\le5K_3(1+2\vep_1+\vep_2).
\end{align}
\end{lmm}

\Proof{Proof.}
The first two inequalities in \refeq{LpBpB'pbd-saw} have already been
explained in \refeq{Lpbd-gen}--\refeq{BpTpbd-gen}.  Similarly, by using
$g_i(p)\le K_i$, $i=1,2$, we have
\begin{align}
B'_p\le p^4\int_{\Td}\hat D(k)^4\hat G_p(k)^2\frac{\text{d}^dk}{(2\pi)^d}
 \le K_1^4K_2^2\underbrace{\int_{\Td}\frac{\hat D(k)^4}{(1-\hat D(k))^2}\,
 \frac{\text{d}^dk}{(2\pi)^d}}_{=(D^{*4}*S_1^{*2})(o)}=K_1^4K_2^2\vep'_2.
\end{align}
For \refeq{Wpbd-saw}, we use $g_3(p)\le K_3$ to obtain
\begin{align}\lbeq{W1prebd-saw1}
0\le(1-\cos k\cdot x)G_p(x)=\int_{\Td}\big(-\hat\Delta_k\hat G_p(l)\big)
 e^{il\cdot x}\,\frac{\text{d}^dl}{(2\pi)^d}\le K_3\int_{\Td}\hat U(k,l)\,
 \frac{\text{d}^dl}{(2\pi)^d},
\end{align}
uniformly in $x$ and $k$.  Then, by \refeq{U-def} and using the Schwarz 
inequality, the right-hand side is further bounded by
\begin{align}\lbeq{W1prebd-saw2}
5K_3\big(1-\hat D(k)\big)\int_{\Td}\hat S_1(l)^2\frac{\text{d}^dl}{(2\pi)^d}
 =5K_3\big(1-\hat D(k)\big)S_1^{*2}(o).
\end{align}
Since $S_1^{*2}(o)=\sum_{n=0}^\infty(2n+1)D^{*2n}(o)=1+2\vep_1+\vep_2$ 
(see \refeq{RWquantities}), this completes the proof of 
Lemma~\ref{lmm:LpBpWpbd-saw}.
\QED

Next, we evaluate the diagrams at $p=1$ by using the trivial inequality 
$G_1(x)\le S_1(x)$.  Here, we do not need the bootstrapping assumptions.

\begin{lmm}\label{lmm:L1B1W1bd-saw}
Let $d\ge5$ and $p=1$.  Then, we have
\begin{align}\lbeq{L1B1B'1bd-saw}
L_1\le\vep_1,&&
B_1\le\vep_2,&&
B'_1\le\vep'_2,
\end{align}
\begin{align}\lbeq{W1bd-saw}
\bigg\|\frac{\hat W_1}{1-\hat D}\bigg\|_\infty\le5(1+2\vep_1+\vep_2).
\end{align}
\end{lmm}

\Proof{Proof.}
The first two inequalities in \refeq{L1B1B'1bd-saw} have already been
explained in \refeq{L1bd-gen}--\refeq{B1T1bd-gen}.  Similarly, by the
trivial inequality $G_1(x)\le S_1(x)$, we have
\begin{align}
B'_1\le\|D^{*4}*S_1^{*2}\|_\infty=\int_{\Td}\frac{\hat D(k)^4}{(1
 -\hat D(k))^2}\,\frac{\text{d}^dk}{(2\pi)^d}=\vep'_2.
\end{align}
Also, by following the same line as \refeq{W1prebd-saw1}--\refeq{W1prebd-saw2}, 
we obtain
\begin{align}
(1-\cos k\cdot x)G_p(x)\le(1-\cos k\cdot x)S_1(x)
&=\int_{\Td}\big(-\hat\Delta_k\hat S_1(l)\big)e^{il\cdot x}\,\frac{\text{d}^dl}
 {(2\pi)^d}\nn\\
&\le\int_{\Td}\hat U(k,l)\,\frac{\text{d}^dl}{(2\pi)^d}\nn\\
&\le5\big(1-\hat D(k)\big)(1+2\vep_1+\vep_2).
\end{align}
This completes the proof of \refeq{W1bd-saw}.
\QED

\Proof{Proof of Proposition~\ref{prp:f-initial}.}
Since $\vep_1$ and $\vep_2$ are finite for $d\ge5$ (see
Table~\ref{table:vep123} in Section~\ref{ss:bcc}) and decreasing in $d$ 
(because $D^{*2n}(o)\equiv\big(\binom{2n}n2^{-2n}\big)^d$ on $\Ld$ 
is decreasing in $d$), we have
\begin{align}
r=\|D\|_\infty+L_1+B_1\stackrel{\refeq{L1B1B'1bd-saw}}\le2^{-d}+\vep_1+\vep_2
 \le
 \begin{cases}
 0.257&[d=5],\\
 0.081&[d\ge6].
 \end{cases}
\end{align}
In addition, by \refeq{odd-bd}--\refeq{sum-diff-bd-even} and 
Lemma~\ref{lmm:L1B1W1bd-saw} (see also Table~\ref{table:vep'2} in 
Section~\ref{ss:Pibds-saw}), we have
\begin{align}
\sum_{n=1}^\infty\hat\pi_1^{\sss(n)}(0)\le
 \begin{cases}
 0.066&[d=5],\\
 0.023&[d\ge6],
 \end{cases}&&
\sup_k\sum_{n=1}^\infty\frac{-\hat\Delta_k\hat\pi_1^{\sss(n)}(0)}{1-\hat D(k)}\le
 \begin{cases}
 1.331&[d=5],\\
 0.120&[d\ge6],
 \end{cases}
\end{align}
which imply that the inequalities in \refeq{suffcond-saw} hold for all $d\ge6$ 
(but not for $d=5$).  Then, by Lemma~\ref{lmm:gbds-saw}, we obtain
\begin{align}
g_1(1)&\le1+\vep_1+\frac{\vep_2(2^{-d}+\vep_1)r}{1-r^2}\le1.021,
 \lbeq{g1(1)bd-saw}\\
g_2(1)&\le\bigg(1-5(1+2\vep_1+\vep_2)\frac{\vep_2^2(1+r^2)}{(1-r^2)^3}
 \bigg)^{-1}\le1.012,\lbeq{g2(1)bd-saw}\\
g_3(1)&\le(1.011)^3\Bigg(1+\frac{3\vep_2}{2^d}+5(1+2\vep_1+\vep_2)\bigg(\vep'_2
 +\frac{\vep_2^2}{(1-r^2)(1-r)^2}+\frac{\vep_2r^2(2-r^2)}{(1-r^2)^2}\bigg)
 \Bigg)^2\nn\\
&\le1.301.\lbeq{g3(1)bd-saw}
\end{align}
Proposition~\ref{prp:f-initial} holds as long as $K_1>1.021$, $K_2>1.012$  and 
$K_3>1.301$.
\QED

\Proof{Proof of Proposition~\ref{prp:f-bootstrapping}.}
Let
\begin{align}
K_1=K_2=1.03,&&
K_3=1.79,
\end{align}
so that Proposition~\ref{prp:f-initial} holds for $d\ge6$.
Using Table~\ref{table:vep123} in Section~\ref{ss:bcc}, we have
\begin{align}\lbeq{rpbd}
r\stackrel{\refeq{LpBpB'pbd-saw}}\le K_12^{-d}+K_1^2K_2\vep_1+K_1^2K_2^2\vep_2
 \le0.088.
\end{align}
In addition, by \refeq{odd-bd}--\refeq{sum-diff-bd-even} and 
Lemma~\ref{lmm:LpBpWpbd-saw} (see also Table~\ref{table:vep'2} in 
Section~\ref{ss:Pibds-saw}), we have
\begin{align}
\sum_{n=1}^\infty\hat\pi_p^{\sss(n)}(0)\le0.025,&&
\sup_k\sum_{n=1}^\infty\frac{-\hat\Delta_k\hat\pi_p^{\sss(n)}(0)}{1-\hat D(k)}
 \le0.257,
\end{align}
which imply that the inequalities in \refeq{suffcond-saw} hold.  Then, similarly 
to \refeq{g1(1)bd-saw}--\refeq{g3(1)bd-saw}, we obtain
\begin{align}
g_1(p)&\le1+K_1^2K_2\vep_1+\frac{K_1^2K_2^2\vep_2(K_12^{-d}+K_1^2K_2\vep_1)
 r}{1-r^2}\le1.023<K_1,\\
g_2(p)&\le\bigg(1-5K_3(1+2\vep_1+\vep_2)\frac{K_1^4K_2^4\vep_2^2(1+r^2)}
 {(1-r^2)^3}\bigg)^{-1}\le1.026<K_2,\\
g_3(p)&\le(1.025)^3\Bigg(\bigg(1+\frac{3K_1^2K_2^2\vep_2}{2^d}\bigg)K_1+5K_3(1
 +2\vep_1+\vep_2)\nn\\
&\hskip8pc\times\bigg(K_1^4K_2^2\vep'_2+\frac{K_1^4K_2^4\vep_2^2}{(1-r^2)
 (1-r)^2}+\frac{K_1^2K_2^2\vep_2r^2(2-r^2)}{(1-r^2)^2}\bigg)\Bigg)^2\nn\\
&\le1.789<K_3.
\end{align}
This completes the proof of Proposition~\ref{prp:f-bootstrapping}.
\QED

\subsection{Further discussion}\label{ss:discussion-saw}
We have been able to prove convergence of the lace expansion for SAW on
$\mL^{\!d\ge6}$ in full detail, in such a small number of pages, rather easily.
This is due to the simple structure of the BCC lattice $\Ld$ and the choice of
the bootstrapping functions $\{g_i(p)\}_{i=1}^3$ (and thanks to the extra
effort explained in the remark after Lemma~\ref{lmm:DBSAW}).
Of course, if we follow the same analysis as Hara and Slade \cite{hs92a,hs92b},
we should be able to extend the result to 5 dimensions.  But, then, the amount
of work and the level of technicality would be almost the same, and it would
not make this survey attractive or accessible to beginners.
Instead of following the analysis of \cite{hs92a,hs92b},
we keep the material as simple as possible and just summarize elements by
which we could improve our analysis.  Those elements are the following.
\begin{enumerate}
\item
Apparently, the largest contribution comes from 
$|\hat\Delta_k\hat\pi_p^{\sss(2)}(0)|$.  To improve its bound, we introduced an 
extra diagram, i.e., $B'_p\equiv\|(pD)^{*4}*G_p^{*2}\|_\infty$.  As a result, 
we were able to improve the applicable range from $d\ge7$ to $d\ge6$.  
It is natural to guess that the introduction of longer bubbles, like
$B^{\sss(n)}_p\equiv\|(pD)^{*2n}*G_p^{*2}\|_\infty$, could result in
the desired applicable range $d\ge5$.  Indeed, its RW counterpart
$(D^{*2n}*S_1^{*2})(o)$ gets smaller as $n$ increases.  However, since 
$B^{\sss(n)}_p$ has the exponentially growing factor $p^{2n}$, 
there must be an optimal $n_*\in\mN$ at which $B^{\sss(n)}_p$ attains its 
minimum.  So far, our naive computation failed to achieve 
convergence of the lace expansion in $\mL^{\!d\ge5}$ 
by merely introducing $B^{\sss(n)}_p$ up to $n=3$.
\item
The reason why we introduced $B'_p$ is because the current bound on
$\|\hat W_p/(1-\hat D)\|_\infty$ in \refeq{W1bd-saw} and \refeq{Wpbd-saw}
is not small.  In particular, the relatively large factor 5 in
\refeq{W1bd-saw} and \refeq{Wpbd-saw} is due to the use of the Schwarz
inequality, as explained in the third footnote.  Therefore, if we could
achieve a better bound on \refeq{trig-identity}, hopefully
without using the Schwarz inequality, it would be of great help.
\item
In \refeq{g1bd-byBWr}--\refeq{g2bd-byBWr}, we discarded the contributions 
from $\hat\Pi_p\even(0)$ and $|\hat\Delta_k\hat\Pi_p\even(0)|$.
By Lemma~\ref{lmm:DBSAW}, we can speculate
$\hat\Pi_p\even(0)\le\hat\Pi_p\odd(0)$ and
$|\hat\Delta_k\hat\Pi_p\even(0)|\ge|\hat\Delta_k\hat\Pi_p\odd(0)|$.  This
means that, if we include their effect into computation, then $g_1(p)$ could
be much closer to 1 (see \refeq{g1bd-byPi}) and $g_2(p)$ could be even smaller
than 1 (see \refeq{g2bd-byPi}), and as a result, we could achieve the desired
applicable limit $d\ge5$.  However, to make use of those even terms, we must
also control lower bounds on $g_1(p)$ and $g_2(p)$, and to do so, we need
nontrivial lower bounds on the lace-expansion coefficients.  Heading towards
this direction would significantly increase the amount of work and technical
details, as in \cite{hs92a,hs92b}, which is against our motivation of writing
this survey.
\item
We evaluated $\hat G_p(k)$ by $\hat S_1(k)\equiv(1-\hat D(k))^{-1}$
uniformly in $k\in\Td$, i.e., in both infrared and ultraviolet regimes.  
However, doing so in the ultraviolet regime (i.e., bounding $G_p(x)$ by
$S_1(x)$ for small $x$) is not efficient, and as a result, it requires 
$d$ to be relatively large.  To overcome this problem, we may want to 
incorporate the idea of ultraviolet regularization, first introduced in 
\cite{an84} for percolation.  This approach has never been investigated 
in the previous lace-expansion work, but it could provide a natural way to 
analyze in dimensions close to $\dc$.
\end{enumerate}

\section{Lace-expansion analysis for percolation}\label{s:LEperc}
In this section, we prove Propositions~\ref{prp:f-initial}--\ref{prp:lace-gen}
for percolation.  First, in Section~\ref{ss:derivation-perc}, we explain the derivation of the lace expansion, Proposition~\ref{prp:lace-gen}, for percolation.
In Section~\ref{ss:Pibds-perc}, we prove bounds on the lace-expansion
coefficients in terms of basic diagrams.  However, unlike SAW, we need more diagrams, such as $T_p$ and $\hat V_p^j(k)$ for $j=0,1,2,3$. Finally, in Section~\ref{ss:TVOHbd}, we prove bounds on those basic diagrams in terms of RW loops, bubbles and triangles and use them to prove Proposition~\ref{prp:f-initial}
on $\mL^{\!d\ge8}$ and Proposition~\ref{prp:f-bootstrapping} on
$\mL^{\!d\ge9}$.  We close this section by addressing potential elements for 
extending the result to 7 dimensions, in Section~\ref{ss:discussion-perc}.

\subsection{Derivation of the lace expansion}\label{ss:derivation-perc}
Proposition~\ref{prp:lace-gen} for percolation is restated as follows.
\begin{prp}[\cite{hs90p}]
  \label{prp:LEPERC}
  For any $p<\pc$ and $N\in\mathbb{Z}_{+}$, there are nonnegative
  functions $\{\pi_p^{\sss(n)}\}_{n=0}^N$ on $\Ld$ such that, if we define
  $\Pi_p^{\sss(N)}$ as
  \begin{equation}
    \Pi_p^{\sss(N)}(x)=\sum_{n=0}^N(-1)^n\pi_p^{\sss(n)}(x),
  \end{equation}
  then we obtain the recursion equation
  \begin{equation}
    \label{eq:perc-laceexp}
    G_p(x)=
      \delta_{o,x}+\Pi_p\ssc{N}(x)+\big((\delta+\Pi_p\ssc{N})*pD*G_p\big)(x)
      +(-1)^{N+1}R_p\ssc{N+1}(x),
  \end{equation}
  where the remainder $R_p\ssc{N+1}(x)$ obeys the bound
  \begin{equation}
    \label{eq:perc-Rbd}
    0\le R_p\ssc{N+1}(x)\le(\pi_p\ssc{N}*G_p)(x).
  \end{equation}
\end{prp}

To prove the above proposition, we first introduce some notions and notation.
\begin{definition}
  Fix a bond configuration
  and let $x,y,u,v\in\Ld$.
  \begin{enumerate}[(i)]
    \item
      Given a bond $b$, we define $\ClusterWithout{b}{x}$ to be the set of vertices
      connected to $x$ in the new configuration obtained by setting $b$ to be vacant.
    \item
      We say that a directed bond $\DirectedBond{u}{v}$ is pivotal
      for the connection from $x$ to $y$
      if $x\Connected u$ occurs in $\ClusterWithout{\Bond{u}{v}}{x}$
      (i.e., $x$ is connected to $u$ without using $\Bond{u}{v}$)
      and if $v\Connected y$ occurs in the complement of $\ClusterWithout{\Bond{u}{v}}{x}$,
      denoted by $\stcomp{\ClusterWithout{\Bond{u}{v}}{x}}$.
      Let $\PivotalBonds{x}{y}$ be the set of directed pivotal bonds for the connection
      from $x$ to $y$.
    \item
      We say that $x$ is doubly connected to $y$, denoted by $x\DoublyConnected y$,
      if either $x=y$ or $x\Connected y$ and $\PivotalBonds{x}{y}=\vno$.
    \item
      Given a set of vertices $A\subset\Ld$,
      we say that $x$ and $y$ are connected in $A$
      if either $x=y\in A$ or
      there is an occupied self-avoiding path from $x$ to $y$ 
consisting of vertices in $A$.
      We write this event as $\event{x\Connected y\text{ in } A}$.
    \item
      Given a set of vertices $A\subset\Ld$,
      we say that $x$ and $y$ are connected through $A$
      if either $x=y\in A$ or
      every occupied self-avoiding path from $x$ to $y$ contains vertices in $A$.
      We write this event as $\event{x\ConnectedThrough{A} y}$.
  \end{enumerate}
\end{definition}

\begin{proof}[Sketch proof of Proposition~\ref{prp:LEPERC}]
First, we derive the first expansion, i.e., \refeq{perc-laceexp} for $N=0$.
By splitting the event $\event{o\Connected x}$ into two depending on whether
or not there is a pivotal bond for the connection from $o$ to $x$, we first 
obtain
  \begin{align}
    \lbeq{first-expansion}
    G_p(x) = \mP_p\big(\underbrace{o\conn x,~\piv(o,x)=\vno}_{=\,\{o\db x\}}\big)
 +\mP_p\big( o\conn x,~\piv(o,x)\ne\vno\big).
  \end{align}
Let
  \begin{align}
    \pi_p^{\sss(0)}(x) = \mP_p( o\db x) - \delta_{o,x}.
  \end{align}
Then, by definition, the first term in \refeq{first-expansion} is 
$\delta_{o,x}+\pi_p^{\sss(0)}(x)$.  To expand the second term in 
\refeq{first-expansion}, we use the first pivotal bond 
$b\equiv\DirectedBond{\underline{b}}{\overline{b}}$ 
for the connection from $o$ to $x$, 
so that $o\DoublyConnected\underline{b}$ in $\ClusterWithout{b}{o}$
  and $\overline{b}\Connected x$ in $\stcomp{\ClusterWithout{b}{o}}$.
  Since those two events are independent of the occupation status
  of $b$, we obtain
  \begin{align}
&\mP_p\big( o\conn x,~\piv(o,x)\ne\vno\big)\nn\\
    &= \sum_{b} \mP_p\Big(
      o\DoublyConnected\underline{b} \text{ in } \ClusterWithout{b}{o},~
      b\text{ occupied},~
      \bar b\Connected x \text{ in } \stcomp{\ClusterWithout{b}{o}}\Big)\nn\\
    &= \sum_b pD(b)
      \Expectation[0]\Big[
        \ind{o\DoublyConnected\underline{b}}~
        \mP_p^{\sss1}\Big(\bar b\conn x \text{ in }\stcomp{\ClusterWithout[0]{b}{o}}
        \Big)\Big],
    \label{eq:nested-expectation}
  \end{align}
where $D(b)$ is the abbreviation for $D(\bar b-\underline b)$, and 
the extra indices\footnote{This rewrite is due to the tower property 
$\mathbb{E}\left[ X\right] = \mathbb{E}\left[\mathbb{E}\left[ X \mid 
\mathcal{G}\right]\right]$, where $\mathbb{E}\left[ X \mid \mathcal{G}\right]$ 
is the conditional expectation of a random variable $X$ with respect to a 
sub-$\sigma$-algebra $\mathcal{G}$.} represent that $\ClusterWithout[0]{b}{o}$ 
is random against $\Expectation[0]$ but deterministic against 
$\mP_p^{\sss1}$.  
In the last line, we have dropped ``in $\ClusterWithout[0]{b}{o}$" by 
using the fact that $\{\bar b\conn x$ in 
$\stcomp{\ClusterWithout[0]{b}{o}}\}=\vno$ when 
$\ind{o\DoublyConnected\bar b\text{ in } \ClusterWithout[0]{b}{o}} \ne 
\ind{o\DoublyConnected\underline{b}}$.

Now we introduce schematic drawings, such as
  \begin{align}
    \label{eq:perc-diagrammatic-representation}
\delta_{o,x}+\pi_p^{\sss(0)}(x)=\picPercCoefficient[colored=true]{0}{o}{x},&&
    \eqref{eq:nested-expectation}=
 \picPercRemainder[colored=true, mode=avoid, colored=true]{1}{o}{x}.
  \end{align}
In the second drawing, the parallel short line segments in the middle 
represents $pD(b)$, which is summed over all bonds $b$ and unlabeled.  
  The dashed two-sided arrow represents mutual avoidance between $\ClusterWithout[0]{b}{o}$
  (in black) and $\ClusterWithout[1]{b}{x}$ (in red).
  By the inclusion-exclusion relation 
$\{\bar b\conn x\text{ in }\stcomp{\ClusterWithout[0]{b}{o}}\}=\{\bar
 b\conn x\}\setminus\{\bar b\ConnectedThrough{\ClusterWithout[0]{b}{o}}x\}$,
we complete the first expansion as
  \begin{align}
    G_p(x)
    &= \picPercCoefficient[colored=true]{0}{o}{x} + \picPercRemainder[colored=true]{1}{o}{x} - \underbrace{\picPercRemainder[colored=true, mode=cross]{1}{o}{x}}_{\equiv \RemainderTerm{1}(x)} \notag\\
    &= \KroneckerDelta{o}{x} + \LaceCoefficient{0}(x)
      + \left(\left(\delta + \LaceCoefficient{0}\right) * pD * G_p\right) (x)
      - \RemainderTerm{1}(x).
    \label{eq:perc-inclusion-exclusion}
  \end{align}
The precise definition of the remainder $R_p^{\sss(1)}(x)$ is 
\begin{align}
R_p^{\sss(1)}(x)=\sum_bpD(b)\mE_p^{\sss0}\Big[
 \ind{o\DoublyConnected\underline{b}}~
 \mP_p^{\sss1}\Big(\bar b\ConnectedThrough{\stcomp{\ClusterWithout[0]{b}{o}}}
 x\Big)\Big].
\end{align}

  Next, we expand the remainder $\RemainderTerm{1}(x)$ to derive the second 
expansion, i.e., \refeq{perc-laceexp} for $N=1$.  To do so, and to derive 
the higher-order expansion later, we have to deal with the event 
$\{v\ConnectedThrough{A}x\}$ for some vertex $v$ and a vertex set $A$.  Let
  \begin{equation}
    E( v, x; A) = \big\{v\ConnectedThrough{A} x\big\} \setminus \bigcup_{b\in\PivotalBonds{v}{x}}\big\{v\ConnectedThrough{A} \underline{b}\big\}.
  \end{equation}
Intuitively, if we regard a percolation cluster of $v$ containing $x$ as 
a string of sausages from $v$ to $x$, then $E(v, x; A)$ is considered to be 
the event that the last sausage is the first one that goes through $A$. 
Then, we can split the event $\event{v\ConnectedThrough{A} x}$ into two 
disjoint events as
  \begin{equation}
    \label{eq:cutting-bond-partition}
    \big\{v\ConnectedThrough{A} x\big\}
    = E\left( v, x; A\right) \cup \Big\{\exists b\in\PivotalBonds{v}{x} \text{ occupied } \&\ v\ConnectedThrough{A}\underline{b}\Big\}.
  \end{equation}
Let 
  \begin{equation}
    \LaceCoefficient{1}(x) = \sum_{b} pD(b)
      \Expectation[0]\Big[
        \ind{o\DoublyConnected\underline{b}}\,
        \mP_p^{\sss1}\Big(E\big(\bar b, x; \ClusterWithout[0]{b}{o}\big)\Big)\Big],
  \end{equation}
so that we have
\begin{align}
\RemainderTerm{1}(x)= \LaceCoefficient{1}(x)+\sum_{b_1}pD(b_1)
&\Expectation[0]\Big[\ind{o\DoublyConnected\underline{b_1}}\nn\\
&\times\mP_p^{\sss1}\Big(\exists b_2\in\piv(\overline{b_1},x)\text{ occupied }
 \&\ \overline{b_1}\ConnectedThrough{\ClusterWithout[0]{b_1}{o}}
\underline{b_2}\Big)\Big].
\end{align}
Notice that the event $\{\exists b\in\PivotalBonds{v}{x}$ occupied 
\& $v\ConnectedThrough{A}\underline{b}\}$ in \refeq{cutting-bond-partition} 
can be rewritten by identifying the first element in 
$\{b\in\piv(v,x):v\ConnectedThrough{A}\underline b\}$ as
\begin{align}\label{eq:cutting-bond-lemma}
&\Event{\exists b\in\piv(v,x)\text{ occupied }\&~v\ConnectedThrough{A}
 \underline b}\nn\\
&=\bigcup_b\Event{E(v, \underline{b}; A) \text{ occurs in }
 \ClusterWithout{b}{v}}\cap \Event{b \text{ occupied}}\cap
 \Event{\overline{b}\Connected x \text{ in } \stcomp{\ClusterWithout{b}{v}}}.
  \end{align}
By this rewrite and using the fact that the first and third events on the 
right-hand side are independent of the occupation status of $b$, we obtain 
(cf., \eqref{eq:nested-expectation})
\begin{align}
\RemainderTerm{1}(x)= \LaceCoefficient{1}(x)+\sum_{b_1,b_2}pD(b_1)pD(b_2)
 \mE_p^{\sss0}\bigg[\ind{o\DoublyConnected\underline{b_1}}\,
&\mE_p^{\sss1}\Big[\indn{E(\overline{b_1}, \underline{b_2};\ClusterWithout[0]
 {b_1}{o})}\nn\\
&\times\mP_p^{\sss2}\Big(\overline{b_2}\Connected x \text{ in }
 \stcomp{\ClusterWithout[1]{b_2}{\overline{b_1}}}\Big)\Big]\bigg],
  \end{align}
where we have dropped ``occurs in $\ClusterWithout[1]{b_2}{\overline{b_1}}$" 
by using the fact that $\{\overline{b_2}\Connected x \text{ in }
 \stcomp{\ClusterWithout[1]{b_2}{\overline{b_1}}}\}=\vno$ when $\ind{E(\overline{b_1}, \underline{b_2}; \ClusterWithout[0]{b_1}{o}) \text{ occurs in } \ClusterWithout[1]{b_2}{\overline{b_1}}} \neq \indn{E(\overline{b_1}, \underline{b_2}; \ClusterWithout[0]{b_1}{o})}$.
By using similar schematic drawings to \eqref{eq:perc-diagrammatic-representation}, 
the above identity for $\RemainderTerm{1}$ is rewritten as
  \begin{equation}
    \label{eq:sausage-1}
    \RemainderTerm{1}(x) = \picPercCoefficient[colored=true]{1}{o}{x} + \picPercRemainder[colored=true, mode=avoid]{2}{o}{x},
  \end{equation}
where the dashed two-sided arrow represents mutual avoidance between $\ClusterWithout[1]{b_2}{\overline{b_1}}$
  (in red) and $\ClusterWithout[2]{b_2}{x}$ (in blue).
By the inclusion-exclusion relation
$\{\overline{b_2}\Connected x \text{ in }\stcomp{\ClusterWithout[1]{b_2}{\overline{b_1}}}\} = \{\overline{b_2}\Connected x\}\setminus \{\overline{b_2}\ConnectedThrough{\ClusterWithout[1]{b_2}{\overline{b_1}}} x\}$,
  we arrive at the second expansion
  \begin{align*}
    \RemainderTerm{1}(x)
    &= \picPercCoefficient[colored=true]{1}{o}{x}\\
      &\qquad + \picPercRemainder[colored=true]{2}{o}{x} - \underbrace{\picPercRemainder[colored=true, mode=cross]{2}{o}{x}}_{\equiv\RemainderTerm{2}(x)}\\
    &= \LaceCoefficient{1}(x) + \left(\LaceCoefficient{1} * pD * G_p\right) (x) - \RemainderTerm{2}(x).
  \end{align*}

  To show how to derive the higher-order expansion coefficients,
  we further demonstrate the expansion of the remainder $\RemainderTerm{2}(x)$ 
by using schematic drawings.
  Using \eqref{eq:cutting-bond-partition} and \eqref{eq:cutting-bond-lemma},
  we can rewrite $\RemainderTerm{2}(x)$ as
  \begin{equation}\lbeq{R2schematicsdrawings}
    \RemainderTerm{2}(x) = \underbrace{\picPercCoefficient[colored=true]{2}{o}{x}}_{\equiv\LaceCoefficient{2}(x)} + \picPercRemainder[colored=true, mode=avoid]{3}{o}{x}.
  \end{equation}
The precise definition of $\LaceCoefficient{2}(x)$ is 
\begin{align}
    \LaceCoefficient{2}(x) =
    \sum_{b_1,b_2} pD(b_1)pD(b_2)
    \Expectation[0]\bigg[
      \ind{o\DoublyConnected\underline{b_1}}\,
      \Expectation[1]\Big[
        \indn{E(\overline{b_1}, \underline{b_2}; \ClusterWithout[0]{b_1}{o})}\,
 \mP_p^{\sss2}\Big(E(\overline{b_2}, x; \ClusterWithout[1]{b_2}{\overline{b_1}})\Big)
 \Big]\bigg].
\end{align}
As in the previous stages of the expansion, the dashed two-sided arrow 
in \refeq{R2schematicsdrawings} represents mutual avoidance between
  $\ClusterWithout[2]{b_3}{\overline{b_2}}$ (in blue) and $\ClusterWithout[3]{b_3}{x}$ (in green).
Then, by the inclusion-exclusion relation
$\{\overline{b_3}\Connected x \text{ in }\stcomp{\ClusterWithout[2]{b_3}
{\overline{b_2}}}\} = \{\overline{b_3}\conn x\}\setminus\{\overline{b_3}
\ConnectedThrough{\ClusterWithout[2]{b_3}{\overline{b_2}}}x\}$,
  we obtain
  \begin{align}
    \RemainderTerm{2}(x)
    &= \picPercCoefficient[colored=true]{2}{o}{x} + \picPercRemainder[colored=true]{3}{o}{x}\\
      &\qquad - \underbrace{\picPercRemainder[colored=true, mode=cross]{3}{o}{x}}_{\equiv\RemainderTerm{3}(x)}\\
    &= \LaceCoefficient{2}(x) + \left(\LaceCoefficient{2} * pD * G_p\right) (x) - \RemainderTerm{3}(x).
  \end{align}
By repeated applications of inclusion-exclusion to the remainders, 
we can derive the higher-order expansion coefficients, such as
  \begin{align}
    \LaceCoefficient{3}(x) &= \picPercCoefficient[colored=true]{3}{o}{x},\\
    \LaceCoefficient{4}(x) &= \picPercCoefficient[colored=true]{4}{o}{x}.
  \end{align}
We complete the sketch proof of Proposition~\ref{prp:LEPERC}.
\end{proof}

\subsection{Diagrammatic bounds on the expansion coefficients} \label{ss:Pibds-perc}
As explained in \emph{Step~1} in Section~\ref{ss:whereandhow}, the
bootstrapping functions $\{g_i(p)\}_{i=1}^3$ are bounded in terms of sums of
$\hat\pi_p^{\sss(n)}(0)$ and $|\hat\Delta_k\hat\pi_p^{\sss(n)}(0)|$.  In this
subsection, we prove bounds on those quantities in terms of basic diagrams,
such as $T_p$ in \eqref{eq:LBT-def} and $\hat V_p^0$, $\hat V_p^1$, $\hat V_p^2$, 
$\hat V_p^3$, defined as
\begin{align}
T_p&=\supNorm{(pD)^{*2} * G_p^{*3}},\\[5pt]
\hat{V}_p^0(k)&= \sum_x \left( pD * G_p\right)(x)^2 \left( 1 - \cos k\cdot x\right),\\
\hat{V}_p^1(k)
&= \sup_x \sum_y \left( pD * G_p\right)(y) \left( 1 - \cos k\cdot y\right) G_p(x - y),\\
  \hat{V}_p^2(k)
&= \sup_x \sum_y \left( pD * G_p\right)(y) \left( 1 - \cos k\cdot y\right) \left( pD * G_p\right) (x - y),\\
\hat{V}_p^3(k)
&= \sup_{x,y} \sum_{\{v_j\}_{j=1}^5} G_p(v_1)\,( pD * G_p) (v_2 - v_1)\big( 1 - \cos k\cdot (v_2 - v_1)\big)\nn\\
&\qquad \times( pD * G_p) (v_3 - v_2)\,(pD*G_p)(v_4 - v_1)\,
 ( pD * G_p)(v_4 - v_2)\nn\\
&\qquad \times ( pD * G_p)(v_5-v_4)\,( pD * G_p)(x-v_5) G_p(y+v_3-v_5).
\end{align}
Recall $r = p \supNorm{D} + L_p + B_p$, and we also define
\begin{equation}
  \rho = T_p \left( 2r + T_p\right) + \left( r + T_p\right) \left( 1 + \frac{B_p}{2} + T_p\right).
\end{equation}

\begin{lmm}[Diagrammatic bounds on the expansion coefficients]
  \label{lmm:DBPERC}
  The expansion coefficients $\LaceCoefficientFT{n}(0)\equiv\sum_{x}\LaceCoefficient{n}(x)$
  and $\abs{\Laplacian{k}\LaceCoefficientFT{n}(0)}\equiv\sum_{x}\left( 1 - \cos k\cdot x\right)\LaceCoefficient{n}(x)$,
  both nonnegative, obey the following bounds:
\begin{align}
\LaceCoefficientFT{n}(0)&\le
 \begin{cases}
 B_p / 2 & [n=0],\\
 ( 1 + B_p / 2 + T_p)^2 r \rho^{n-1} & [n\geq 1],
 \end{cases}\label{eq:piperc-sumbd}\\[5pt]
\abs{\Laplacian{k}\LaceCoefficientFT{0}(0)}
&\le\frac12\hat{V}_p^0(k), \label{eq:piperc-diffbd-0}\\[5pt]
\abs{\Laplacian{k}\LaceCoefficientFT{1}(0)}
&\le\bigg(1 + 2 (B_p+r) + \frac{3}{4} B_p (B_p+2r) + 3rT_p\bigg)
 \hat{V}_p^0(k)\nn\\
&\quad+( 8 + 6B_p + 9T_p)T_p\hat{V}_p^2(k).\label{eq:piperc-diffbd-1}
\end{align}
  For $m\geq 1$,
\begin{align}\label{eq:piperc-diffbd-even}
\abs{\Laplacian{k}\LaceCoefficientFT{2m}(0)}\nn\\
\le( 4m+1)&\Bigg(\rho^{2m-1}\Big(2r \hat{V}_p^0(k) + T_p \hat{V}_p^2(k)
 + T_p \hat{V}_p^1(k)\Big) \bigg( 1 + \frac{B_p}{2} + T_p\bigg)\nn\\
&+ m\rho^{2m-1}\big(\hat{V}_p^2(k) + \hat{V}_p^1(k)\big) \left( 1
 + \frac{B_p}{2} + T_p\right)^2\nn\\
&+ \rho^{2m-2}\bigg(( 2m-1) \Big(r^2 \hat{V}_p^0(k) + r T_p \hat{V}_p^2(k)
 + r T_p \hat{V}_p^1(k) + T_p^2 \hat{V}_p^1(k)\Big)\nn\\
&\qquad\qquad+ ( m-1) \big(r T_p\hat{V}_p^0(k) + \hat{V}_p^3(k)\big)\bigg)
 \left( 1 + \frac{B_p}{2} + T_p\right)^2\Bigg),
\end{align}
  and
\begin{align}\label{eq:piperc-diffbd-odd}
\abs{\Laplacian{k}\LaceCoefficientFT{2m+1}(0)}\nn\\
\le( 4m+3)&\Bigg(2\rho^{2m} \Big(r\hat{V}_p^0(k)
 + T_p\hat{V}_p^2(k)\Big) \left( 1 + \frac{B_p}{2} + T_p\right)\nn\\
&+ \rho^{2m}\big(( m+1) \hat{V}_p^2(k) + m\hat{V}_p^1(k)\big)
 \left( 1 + \frac{B_p}{2} + T_p\right)^2\nn\\
&+ m\rho^{2m-1}\bigg(2\Big(r^2 \hat{V}_p^0(k) + r T_p \hat{V}_p^2(k)
 + r T_p \hat{V}_p^1(k) + T_p^2 \hat{V}_p^1(k)\Big)\nn\\
&\qquad\qquad+ rT_p\hat{V}_p^0(k) + \hat{V}_p^3(k)\bigg)
 \left( 1 + \frac{B_p}{2} + T_p\right)^2\Bigg).
\end{align}
\end{lmm}

The rest of this subsection is devoted to showing the above bounds 
on $\hat\pi_p^{\sss(n)}(0)$ and $|\hat\Delta_k\hat\pi_p^{\sss(n)}(0)|$ for 
$n=0$ in Section~\ref{sss:pi0}, for $n=1$ in Section~\ref{sss:pi1}, 
for $n=2$ in Section~\ref{sss:pi2}, and for $n\ge3$ in 
Section~\ref{sec:estimation-perc-lace-diff}.

\subsubsection{Bounds on $\LaceCoefficientFT{0}(0)$ and $\abs{\Laplacian{k}\LaceCoefficientFT{0}(0)}$}\label{sss:pi0}
By the Boolean and BK inequalities, we obtain
\begin{align}
  \LaceCoefficient{0}(x)
  &=\Probability( o\DoublyConnected x) - \KroneckerDelta{o}{x} \notag\\
  &=\mP_p\Bigg(\bigcup_{\substack{b_1,\, b_2\in B(o)\\ (b_1 \prec b_2)}}
    \Bigl\{
      \Event{b_1 \text{ occupied},\ \overline{b_1}\Connected x}
      \DisjointPath \Event{b_2 \text{ occupied},\ \overline{b_2}\Connected x}
    \Bigr\}\Bigg)\nn\\
  &\leq \sum_{\substack{b_1,\, b_2\in B(o)\\ (b_1 \prec b_2)}}
    pD(b_1) G_p(x - \overline{b_1}) pD(b_2) G_p(x - \overline{b_2}) \notag\\
  &\leq \frac{1}{2} \left( pD * G_p\right)(x)^2,
  \label{eq:bound-pi0}
\end{align}
where we have used the ordering $\prec$ introduced above \eqref{eq:lowestoccupied}.
The factor $1/2$ in the last line is due to ignoring the ordering.
Then, summing over $x$ yields \refeq{piperc-sumbd} for $n=0$.

The bound \refeq{piperc-diffbd-0} on $|\hat\Delta_k\hat\pi_p^{\sss(0)}(0)|$ is 
also achieved by multiplying both sides of \eqref{eq:bound-pi0} by 
$1-\cos k\cdot x$ and summing the resulting inequality over $x$.

\subsubsection{Bounds on $\LaceCoefficientFT{1}(0)$ and $\abs{\Laplacian{k}\LaceCoefficientFT{1}(0)}$}\label{sss:pi1}
First, we prove \refeq{piperc-sumbd} for $n=1$ and 
\refeq{piperc-diffbd-1} by assuming the following diagrammatic bound on 
$\LaceCoefficient{1}(x)$:
\begin{align}\label{eq:perc-diagrammatic-estimate-1}
&\LaceCoefficient{1}(x)\nn\\
&\le\picPercFirstLace[leftline=vanish, rightline=vanish]{o}{}{}{x}
    + \frac{1}{2} \picPercFirstLace[leftline=vanish, rightline=vertex]{o}{}{}{x}
    + \picPercFirstLace[leftline=vanish, rightline=side]{o}{}{}{x}
    + \frac{1}{2} \picPercFirstLace[leftline=vertex, rightline=vanish]{o}{}{}{x}
    + \frac{1}{4} \picPercFirstLace[leftline=vertex, rightline=vertex]{o}{}{}{x} \nn\\
&\quad
    + \frac{1}{2} \picPercFirstLace[leftline=vertex, rightline=side]{o}{}{}{x}
    + \picPercFirstLace[leftline=side, rightline=vanish]{o}{}{}{x}
    + \frac{1}{2} \picPercFirstLace[leftline=side, rightline=vertex]{o}{}{}{x}
    + \picPercFirstLace[leftline=side, rightline=side]{o}{}{}{x},
\end{align}
where we have used the following two types of line segments:
\begin{align}\label{eq:diagram-line} 
\picLine{o}{x}=G_p(x),&&
\picLine[onestep]{o}{x}=( pD * G_p)(x).
\end{align}
As in the case for SAW (cf., e.g., \refeq{pisaw-diagr}), the unlabeled 
vertices are summed over $\Ld$.
The proof of \refeq{perc-diagrammatic-estimate-1} is given at the end of 
Section~\ref{sss:pi1}.

\Proof{Proof of \eqref{eq:piperc-sumbd} for $n=1$ assuming 
\refeq{perc-diagrammatic-estimate-1}.}
  The bound on $\LaceCoefficientFT{1}(0)$ is obtained by summing both sides of
  \eqref{eq:perc-diagrammatic-estimate-1} over $x$
  and repeatedly using translation-invariance.  For example,
  \begin{align}
    \picPercFirstLace[leftline=side, rightline=vertex]{o}{}{}{}
    = \picPercFirstLace[leftline=side, rightline=vanish]{}{o}{}{}
\underbrace{\picBubble[toplinegap=true, bottomlinegap=true]{o}{}}_{\leq B_p}
    \le \underbrace{\picTriangle{}{o}{}}_{\leq T_p} \bigg(\underbrace{\sup_x
 \picOpenBubble[toplinegap=true]{o}{x}{}}_{\le r~(\because\refeq{boundbyr})}
 \bigg) B_p
 \le T_p r B_p,
  \end{align}
  and
  \begin{align}
    \picPercFirstLace[leftline=side, rightline=side]{o}{}{}{}
 \le\underbrace{\picTriangle[left]{}{o}{}}_{\leq T_p}
 \sup_{x} \sum_{y} \picInsRightTriangle{o}{x}{y}
&\leq T_p
      \bigg(\sup_{x,\, z}\sum_{y} \picTwoLines[toplinegap=true]{o}{x}{y}{y+z}\bigg)
      \underbrace{\picTriangle[right]{o}{}{}}_{\leq T_p}\nn\\
&\leq T_p \bigg(\underbrace{\sup_{x,\, z} \picOpenBubble[toplinegap=true]{o}{x-z}{}}_{\le r~(\because\refeq{boundbyr})}\bigg) T_p
  \le T_p r T_p.
  \end{align}
  Applying the same analysis to the other diagrams, we obtain
  \begin{align*}
    \LaceCoefficientFT{1}(0)
    &\leq \left( pD * G_p^{*2}\right)(o)
      + \frac{1}{2} \left( pD * G_p^{*2}\right)(o) B_p
      + r T_p
      + \frac{1}{2} B_p \left( pD * G_p^{*2}\right)(o)\\
      &\quad
      + \frac{1}{4} B_p \left( pD * G_p^{*2}\right)(o) B_p
      + \frac{1}{2} B_p r T_p
      + T_p r
      + \frac{1}{2} T_p r B_p
      + T_p r T_p\\
    &= \left( pD * G_p^{*2}\right)(o) \left( 1 + \frac{B_p}{2}\right)^2
      + r \Bigg( 2 \left( 1 + \frac{B_p}{2}\right) T_p + T_p^2\Bigg)\\
    &\stackrel{\mathclap{\eqref{eq:boundbyr}}}{\leq}~\left( 1 + \frac{B_p}{2} + T_p\right)^2 r,
  \end{align*}
  as required.
\QED

\Proof{Sketch proof of \eqref{eq:piperc-diffbd-1}.}
The bound on $\abs{\Laplacian{k}\LaceCoefficientFT{1}(0)}$ is obtained by 
multiplying $1 - \cos k\cdot x$ to both sides of
\eqref{eq:perc-diagrammatic-estimate-1} and summing the resulting expression 
over $x$.
To decompose the diagrams into the basic diagrams, we also use the telescopic 
inequality \eqref{eq:telescope}, translation-invariance and the trivial inequality
\begin{align}\lbeq{bound-at_least_one_line-perc}
  \picLine{o}{x} \ind{x\neq o} = G_p(x) \ind{x\neq o}
  \leq \left( pD * G_p\right)(x) = \picLine[onestep]{o}{x}.
\end{align}
For example,
  \begin{align}
&\sum_x\picPercFirstLace[leftline=vanish, rightline=vertex]{o}{}{}{x}(1-\cos k\cdot x)\nn\\
    &\stackrel{\mathclap{\eqref{eq:telescope}}}{\leq}
      2 \sum_{x,\, y} \picPercFirstLace[leftline=vanish, rightline=vertex]{o}{}{y}{x}
      \Big(\left( 1 - \cos k\cdot y\right) + \left( 1 - \cos k\cdot (x - y)\right)\Big)\nn\\
    &\stackrel{\mathclap{\eqref{eq:bound-at_least_one_line-perc}}}{\leq} 2 \underbrace{\sum_{y} \picBubble[toplinegap=true, bottomlinegap=true]{o}{y} ( 1 - \cos k\cdot y)}_{=\hat{V}_p^0(k)}~
        \underbrace{\picBubble[toplinegap=true, bottomlinegap=true]{o}{}}_{\le B_p}\nn\\
      &\quad
      + 2 \underbrace{\picBubble[toplinegap=true]{o}{}}_{\le r~(\because\refeq{boundbyr})}~
        \underbrace{\sup_{y}\sum_{x} \picBubble[toplinegap=true, bottomlinegap=true]{y}{x} \big( 1 - \cos k\cdot (x - y)\big))}_{=\hat{V}_p^0(k)}\nn\\
    &\le 2 \hat{V}_p^0(k) B_p + 2 r \hat{V}_p^0(k).
\end{align}
Another example is the following:
\begin{align}\label{eq:example-first_lace_9-diff}
&\sum_x\picPercFirstLace[leftline=side, rightline=side]{o}{}{}{x}(1-\cos k\cdot
 x)\nn\\
&\stackrel{\mathclap{\eqref{eq:telescope}}}\le3\sum_{\{x_j\}_{j=1}^3}
 \picPercFirstLace[leftline=side, rightline=side]{o}{x_1}{x_2}{x_3}\sum_{j=1}^3
 \big((1 - \cos k\cdot(x_j-x_{j-1})\big),
\end{align}
where $x_0=o$.  The contribution from $1 - \cos k\cdot x_1$ is 
bounded by
\begin{align}
3\bigg(\underbrace{\sup_x\sum_y\picOpenBubble[mode=rightopen,
 toplinegap=true]{y}{o}{x}( 1 - \cos k\cdot y)}_{=\hat V_p^2(k)}
 \bigg)\bigg(\underbrace{\sup_x
 \picOpenTriangle[mode=rightopen]{}{}{o}{x}}_{=T_p}\bigg)
 \picTriangle[right]{o}{}{}
 \stackrel{\mathclap{\eqref{eq:bound-at_least_one_line-perc}}}{\leq} 3 \hat{V}_p^2(k) T_p T_p.
  \end{align}
The contribution from $1 - \cos k\cdot(x_3-x_2)$ obeys the same bound, because
\begin{align}
3~\picTriangle[left]{o}{}{}~\bigg(\underbrace{\sup_x
 \picOpenTriangle[mode=leftopen]{o}{x}{}{}}_{=T_p}\bigg)
        \sup_{x}\sum_{y} \picOpenBubble[mode=leftopen, toplinegap=true]{o}{y}{x}( 1 - \cos k\cdot y)
&\stackrel{\mathclap{\eqref{eq:bound-at_least_one_line-perc}}}{\leq} 3 T_p T_p \hat{V}_p^2(k).
\end{align}
The contribution from $1 - \cos k\cdot(x_2-x_1)$ is bounded by
\begin{align}
3~\picTriangle[left]{o}{}{}~\bigg(\sup_{x,z}\sum_{y}
 \picTwoLines[toplinegap=true]{o}{x}{y}{y+z}( 1 - \cos k\cdot y)\bigg)
 \picTriangle[right]{o}{}{}\nn\\
\leq 3 T_p\bigg(\sup_{x,z}\sum_{y}\picOpenBubble[toplinegap=true]{o}{x-z}{y}
 ( 1 - \cos k\cdot y)\bigg)T_p
&\stackrel{\mathclap{\eqref{eq:bound-at_least_one_line-perc}}}\le
 3 T_p \hat{V}_p^2(k) T_p.
\end{align}
As a result, \eqref{eq:example-first_lace_9-diff} is bounded by 
$9T_p^2\hat{V}_p^2(k)$.  The other terms can be estimated similarly.  
We complete the proof of \eqref{eq:piperc-diffbd-1}.
\QED

\Proof{Proof of \refeq{perc-diagrammatic-estimate-1}.}
First, we recall the definition of $\LaceCoefficient{1}(x)$:
\begin{equation}
  \LaceCoefficient{1}(x) = \sum_{b} pD(b)
    \Expectation[0]\Big[
      \ind{o\DoublyConnected\underline{b}}\,
      \mP_p^{\sss1}\Big(E\big(\overline{b}, x; \ClusterWithout[0]{b}{o}\big)\Big)\Big].
\end{equation}

\begin{figure}[tp]
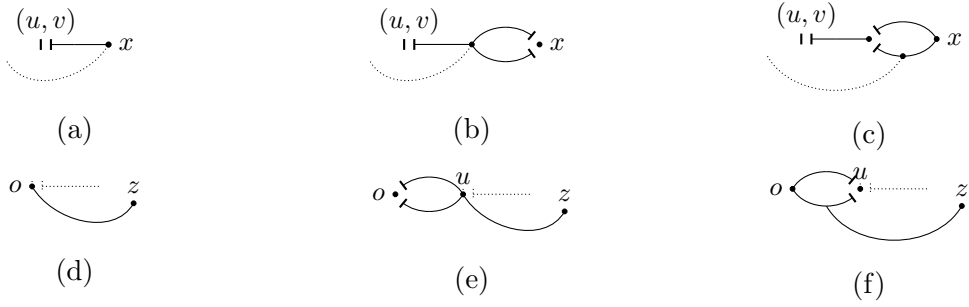

  \begin{subfigure}{0.33\linewidth}
    \centering
    \picLastSausage[mode=vanish, auxiliary=true]{\DirectedBond{u}{v}}{x}
    \caption{}
    \label{fig:last-sausage-vanish}
  \end{subfigure}%
  \begin{subfigure}{0.33\linewidth}
    \centering
    \picLastSausage[mode=vertex, auxiliary=true]{\DirectedBond{u}{v}}{x}
    \caption{}
    \label{fig:last-sausage-vertex}
  \end{subfigure}%
  \begin{subfigure}{0.33\linewidth}
    \centering
    \picLastSausage[mode=side, auxiliary=true]{\DirectedBond{u}{v}}{x}
    \caption{}
    \label{fig:last-sausage-side}
  \end{subfigure}
  \begin{subfigure}{0.33\linewidth}
    \centering
    \picFirstSausage[mode=vanish, auxiliary=true]{o}{u}{z}
    \caption{}
    \label{fig:first-sausage-vanish}
  \end{subfigure}%
  \begin{subfigure}{0.33\linewidth}
    \centering
    \picFirstSausage[mode=vertex, auxiliary=true]{o}{u}{z}
    \caption{}
    \label{fig:first-sausage-vertex}
  \end{subfigure}%
  \begin{subfigure}{0.33\linewidth}
    \centering
    \picFirstSausage[mode=side, auxiliary=true]{o}{u}{z}
    \caption{}
    \label{fig:first-sausage-side}
  \end{subfigure}
  \caption{Schematic drawings for the events in \eqref{eq:distinguish-sausage} 
and \eqref{eq:distinguish-first-sausage}.  The real and dotted line segments are on 
different probability spaces.  The arcs having short line segments at one of the 
two end vertices represent nonzero connections.}
  \label{fig:last-first-sausage}
\end{figure}

Let $x\DoublyConnectedThrough{A} y$ be the event that $x$ is doubly connected 
to $y$ through $A$, i.e., there are at least two occupied paths from $x$ to $y$ 
and every occupied path from $x$ to $y$ has vertices of $A$.  Then, by 
definition, we have
\begin{equation}
  \label{eq:distinguish-last-event}
  E(v, x; A)
  \subset \bigcup_y\Big\{\Event{v \Connected y}
    \DisjointPath \big\{y\DoublyConnectedThrough{A} x\big\}\Big\}.
\end{equation}
Splitting $\{y\DoublyConnectedThrough{A}x\}$ into three events 
depending on where the double connection traverses $A$, we obtain
\begin{align}
\big\{y\DoublyConnectedThrough{A} x\big\}&\subset\event{y=x\in A}\cup
 \{y\db x\ne y\in A\}\nn\\
&\quad\cup\,\bigcup_{z(\ne y)}\Big\{\{y\conn z\in A\}\circ\{z\conn x\}\circ\{y
 \conn x\ne y\}\Big\},
\end{align}
hence (see Figure~\ref{fig:last-first-sausage}(a)--(c))
\begin{align}\label{eq:distinguish-sausage}
E(v, x; A)&\subset\underbrace{\{v\conn x\in A\}}_\text{(a)}\,\cup\,
 \underbrace{\bigcup_y\big\{\{v\conn y\in A\}\circ\{y\db x\ne y\}
 \big\}}_\text{(b)}\nn\\
&\quad\cup\underbrace{\bigcup_{\substack{y,z\\ (y\ne z)}}\Big\{\{v\conn
 y\}\circ\{y\conn z\in A\}\circ\{z\conn x\}\circ\{y\conn x\ne y\}
 \Big\}}_\text{(c)}.
\end{align}
Then, by the BK inequality and the argument around \eqref{eq:bound-pi0} 
to derive the factor $1/2$, and applying the trivial inequality 
\eqref{eq:bound-at_least_one_line-perc} to the nonzero connections 
(e.g., $\mP_p(y\conn x\ne y)\le(pD*G_p)(x-y)$), we obtain
\begin{align}\label{eq:bound-last-sheet}
\sum_{v} pD(v - u)&\Probability\big( E(v, x; A)\big)\nn\\
\le\sum_{y,z} \ind{z\in A}&\bigg(( pD * G_p)(x - u) \KroneckerDelta{y}{x}
 \KroneckerDelta{z}{x}+ \frac12( pD * G_p)(y - u) ( pD * G_p)(x - y)^2
 \KroneckerDelta{z}{y}\nn\\
&+( pD * G_p)(y - u) ( pD * G_p)(z - y) G_p(x - z) ( pD * G_p)(x - y)\bigg),
\end{align}
and therefore, by using the diagrammatic representations in \refeq{diagram-line},
\begin{align}\label{eq:last-diagrams}
\LaceCoefficient{1}(x) \leq \sum_{u,z}\mP_p(o\db u,~o\conn z)\Bigg(
 \underbrace{\picLastDiagram[mode=vanish, position=top]{u}{}{}{x}\delta_{z,x}
 +\frac12\picLastDiagram[mode=bubble, position=top]{u}{z}{}{x}
 +\picLastDiagram[mode=triangle, position=top]{u}{z}{}{x}}_{\dpst\equiv\bR{z}ux}
 \Bigg).
\end{align}
We emphasize that each of the above line segments is a two-point function 
$G_p$, and not a connection event described in 
Figure~\ref{fig:last-first-sausage}.

It remains to investigate $\mP_p(o\db u,~o\conn z)$ in 
\eqref{eq:last-diagrams}.
Splitting the event into three depending on which vertex on the backbone 
from $o$ to $u$ a connection to $z$ comes out of, we have (see 
Figure~\ref{fig:last-first-sausage}(d)--(f))
\begin{align}\label{eq:distinguish-first-sausage}
&\{o\db u,~o\conn z\}\nn\\[5pt]
&\subset\underbrace{\{o=u\conn z\}}_\text{(d)}\,\cup\,\underbrace{\big\{\{o\db
 u\ne o\}\circ\{u\conn z\}\big\}}_\text{(e)}\nn\\
&\quad\cup\,\underbrace{\bigcup_w\Big\{\{o\conn u\ne o\}\circ\{o\conn w\}\circ
 \{w\conn u\ne w\}\circ\{w\conn z\}\Big\}.}_\text{(f)}
\end{align}
Then, again, by the BK inequality and the argument around \eqref{eq:bound-pi0} 
to derive the factor $1/2$, and applying the trivial inequality 
\eqref{eq:bound-at_least_one_line-perc} to the nonzero connections, we obtain
\begin{align}\label{eq:first-diagrams}
\mP_p(o\db u,~o\conn z)\le\underbrace{\picFirstDiagram[vanish]{o}{}{}{z}
 \delta_{u,o}+\frac12\picFirstDiagram[bubble]{o}{}{u}{z}
 +\picFirstDiagram[triangle]{o}{}{u}{z}}_{\dpst\equiv\bL{z}u}.
\end{align}

Combining this with \refeq{last-diagrams}, we obtain 
the diagrammatic bound \refeq{perc-diagrammatic-estimate-1}, as required.
\QED

\subsubsection{Bounds on $\LaceCoefficientFT{2}(0)$ and $\abs{\Laplacian{k}\LaceCoefficientFT{2}(0)}$}\label{sss:pi2}
We organize this section in a different way from the previous section for 
the case of $n=1$.  We first explain the diagrammatic bound 
\refeq{perc-diagrammatic-estimate-2} on $\LaceCoefficient{2}(x)$ for a fixed 
$x$.  Then, by using this, we prove the bounds \refeq{piperc-sumbd} for $n=2$ 
and \refeq{piperc-diffbd-even} for $m=1$.

Now we start investigating $\LaceCoefficient{2}(x)$ for a fixed $x$, which is 
defined as
\begin{equation}
  \LaceCoefficient{2}(x) =
  \sum_{b_1,b_2} pD(b_1)pD(b_2)
  \Expectation[0]\bigg[
    \ind{o\DoublyConnected\underline{b_1}}
 \Expectation[1]\Big[
 \indn{E(\overline{b_1}, \underline{b_2}; \ClusterWithout[0]{b_1}{o})}\,
 \mP_p^{\sss2}\Big(E\big(\overline{b_2}, x; \ClusterWithout[1]{b_2}
 {\overline{b_1}}\big)\Big)\Big]\bigg].
\end{equation}
First, by using \eqref{eq:bound-last-sheet}--\refeq{last-diagrams},
we obtain 
\begin{align}\lbeq{perc-lace-coefficient-2}
\LaceCoefficient{2}(x)&\le\sum_{b_1}pD(b_1)\sum_{u_2,z_2}
 \Expectation[0]\Big[\ind{o\DoublyConnected\underline{b_1}}\,
 \mP_p^{\sss1}\Big(E\big(\overline{b_1},u_2;\ClusterWithout[0]{b_1}{o}\big)
 \cap\{\overline{b_1}\conn z_2\}\Big)\Big]\times\bR{z_2}{u_2}x.
\end{align}

\begin{figure}[t]
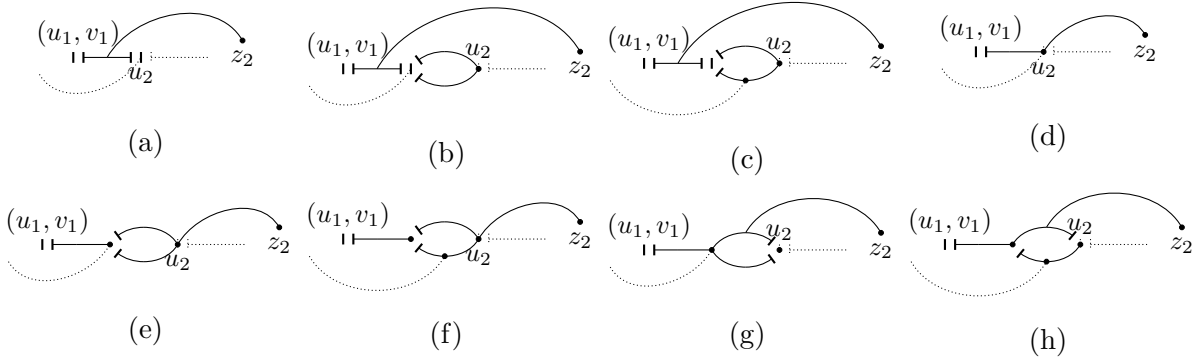

  \begin{subfigure}{0.25\linewidth}
    \centering
    \picMiddleSausagePB[mode=vanish, auxiliary=true]{\DirectedBond{u_1}{v_1}}{u_2}{z_2}
    \caption{}
    \label{fig:middle-sausage-pb-vanish}
  \end{subfigure}%
  \begin{subfigure}{0.25\linewidth}
    \centering
    \picMiddleSausagePB[mode=vertex, auxiliary=true]{\DirectedBond{u_1}{v_1}}{u_2}{z_2}
    \caption{}
    \label{fig:middle-sausage-pb-vertex}
  \end{subfigure}%
  \begin{subfigure}{0.25\linewidth}
    \centering
    \picMiddleSausagePB[mode=side, auxiliary=true]{\DirectedBond{u_1}{v_1}}{u_2}{z_2}
    \caption{}
    \label{fig:middle-sausage-pb-side}
  \end{subfigure}%
  \begin{subfigure}{0.25\linewidth}
    \centering
    \picMiddleSausageDC[modebottom=vanish, modetop=vanish, auxiliary=true]{\DirectedBond{u_1}{v_1}}{u_2}{z_2}
    \caption{}
    \label{fig:middle-sausage-dc-vanish}
  \end{subfigure}
  \begin{subfigure}{0.25\linewidth}
    \centering
    \picMiddleSausageDC[modebottom=vertex, modetop=vertex, auxiliary=true]{\DirectedBond{u_1}{v_1}}{u_2}{z_2}
    \caption{}
    \label{fig:middle-sausage-dc-vertex}
  \end{subfigure}%
  \begin{subfigure}{0.25\linewidth}
    \centering
    \picMiddleSausageDC[modebottom=side, modetop=vertex, auxiliary=true]{\DirectedBond{u_1}{v_1}}{u_2}{z_2}
    \caption{}
    \label{fig:middle-sausage-dc-side-bottom}
  \end{subfigure}%
  \begin{subfigure}{0.25\linewidth}
    \centering
    \picMiddleSausageDC[modebottom=vertex, modetop=side, auxiliary=true]{\DirectedBond{u_1}{v_1}}{u_2}{z_2}
    \caption{}
    \label{fig:middle-sausage-dc-side-top}
  \end{subfigure}%
  \begin{subfigure}{0.25\linewidth}
    \centering
    \picMiddleSausageDC[modebottom=side, modetop=side, auxiliary=true]{\DirectedBond{u_1}{v_1}}{u_2}{z_2}
    \caption{}
    \label{fig:middle-sausage-dc-side}
  \end{subfigure}
  \caption{Schematic drawings for the events in 
\eqref{eq:distinguish-middle-sausage} (cf., Figure~\ref{fig:last-first-sausage}).  
A connection to $z_2$ comes out of the backbone from $v_1$ to $u_2$ before 
the last sausage ((a)--(d)) or from the last sausage ((e)--(h)).}
\label{fig:middle-sausage}
\end{figure}

Next, we have to deal with the event $E(v_1,u_2;A)\cap\{v_1\conn z_2\}$ for a 
given set $A$.  By \refeq{distinguish-last-event}, we first obtain the relation
\begin{align}\label{eq:distinguish-middle-sausage}
E(v_1,u_2;A)\cap\{v_1\conn z_2\}\subset\bigcup_y\Big\{\{v_1\conn y\}\circ
 \big\{y\DoublyConnectedThrough{A} u_2\big\}\Big\}\cap\{v_1\conn z_2\}.
\end{align}
Next, we split the event into two depending on which vertex on the backbone from 
$v_1$ to $u_2$ a connection to $z_2$ comes out of: either before or from 
the last sausage.  Then, we split each of the two into four events depending 
on where the double connection from $y$ to $u_2$ traverses $A$.  
The resulting eight events are depicted in Figure~\ref{fig:middle-sausage}.
Finally, by the BK inequality and the argument around \eqref{eq:bound-pi0} 
to derive the factor $1/2$, and applying the trivial inequality 
\eqref{eq:bound-at_least_one_line-perc} to the nonzero connections, we obtain
\begin{gather}
\sum_{b_1}pD(b_1)\Expectation[0]\Big[\ind{o\DoublyConnected\underline{b_1}}\,
 \mP_p^{\sss1}\Big(E\big(\overline{b_1},u_2;\ClusterWithout[0]{b_1}{o}\big)
 \cap\{\overline{b_1}\conn z_2\}\Big)\Big]\nn\\
 \le\sum_{u_1,z_1}\mP_p(o\db u_1,~o\conn z_1)\times\bM{z_1}{u_1}{u_2}{z_2},
 \lbeq{bMdef}
\end{gather}
where (in some of the following diagrams, we use the identity 
$(pD*G_p)(x)=(G_p*pD)(x)$)
\begin{align}\lbeq{middle-diagrams}
\bM{z_1}{u_1}{u_2}{z_2}
&=\picMiddleDiagramPB[vanish]{u_1}{}{u_2}{}{}{z_2}\delta_{z_1,u_2}
 +\frac12\picMiddleDiagramPB[bubble]{u_1}{}{z_1}{}{u_2}{z_2}
 +\picMiddleDiagramPB[triangle]{u_1}{}{}{z_1}{u_2}{z_2}\nn\\
&\quad+\picMiddleDiagramDC[vanish=true]{u_1}{}{}{u_2}{}{z_2}\delta_{z_1,u_2}
 +\frac12\picMiddleDiagramDC[first=vanished, last=vanished]{u_1}{z_1}{}{u_2}{}
 {z_2}+\picMiddleDiagramDC[first=unvanished, last=vanished]{u_1}{}{z_1}{u_2}{}
 {z_2}\nn\\
&\quad+\picMiddleDiagramDC[first=vanished, last=unvanished]{u_1}{z_1}{}{u_2}{}
 {z_2}+\picMiddleDiagramDC[first=unvanished, last=unvanished]{u_1}{}{z_1}{u_2}{}
 {z_2}.
\end{align}


Finally, by using \refeq{first-diagrams}, we arrive at 
\begin{gather}\lbeq{perc-diagrammatic-estimate-2}
\LaceCoefficient{2}(x)\le\sum_{\substack{u_1,u_2,z_1,z_2}}\bL{z_1}{u_1}\times
 \bM{z_1}{u_1}{u_2}{z_2}\times\bR{z_2}{u_2}x,
\end{gather}
which consists of 72 ($=3\times 8\times 3$) terms.
\QED

\Proof{Sketch proof of \eqref{eq:piperc-sumbd} for $n=2$.}
The bound on $\LaceCoefficientFT{2}(0)$ is obtained by summing both sides of
  \eqref{eq:perc-diagrammatic-estimate-2} over $x$ and repeatedly using
  translation-invariance.
  For example, 
the combination of the third diagrams in \refeq{first-diagrams}, 
\refeq{middle-diagrams} and \refeq{last-diagrams} is bounded as
\begin{gather}
\sum_x\sum_{u_1,u_2,z_1,z_2}
      \picFirstDiagram[triangle]{o}{}{u_1}{z_1}
      \picMiddleDiagramPB[triangle]{u_1}{}{}{z_1}{u_2}{z_2}
      \picLastDiagram[mode=triangle, position=bottom]{u_2}{}{z_2}{x} \nn\\
\le\underbrace{\picTriangle[left]{}{o}{}}_{\leq T_p}~
      \sup_{x,\, z}\sum_{y} \picWeightedMiddlePB{o}{x}{}{}{y}{y+z}
      \underbrace{\picTriangle[right]{o}{}{}}_{\leq T_p}.
    \label{eq:example-second_lace_333}
\end{gather}
Also, the middle diagram is bounded as
\begin{align}\lbeq{middlesection1}
\sup_{x,z}\sum_y\picWeightedMiddlePB{o}{x}{}{}{y}{y+z}
&\le\sup_{x} \picInsMiddlePB{o}{x}{}~~\sup_{v,w,z}\sum_y
 \picTwoLines[bottomlinegap=true]{w}{v}{y}{y+z}\nn\\
&\le\underbrace{\sup_{x,u}\sum_{y}\picInsOpenPartOfMiddlePB{o}{x}{y}{y
 +u}}_{\le T_p}~
 \underbrace{\picInsPartOfMiddlePB{o}~}_{\le T_p}~
 \underbrace{\sup_{v,w,z}\picOpenBubble[bottomlinegap=true]{w}{v-z}{}}_{\le
 r~(\because\refeq{boundbyr})}.
\end{align}
Therefore, \eqref{eq:example-second_lace_333} is bounded above by 
$T_p T_p^2 r T_p$.

Another example is the combination of the third diagram in 
\refeq{first-diagrams}, the last diagram in \refeq{middle-diagrams} and 
the third diagram in \refeq{last-diagrams}, which is bounded as
  \begin{gather}
    \sum_{x}\sum_{u_1, u_2, z_1, z_2}
      \picFirstDiagram[triangle]{o}{}{u_1}{z_1}
      \picMiddleDiagramDC[first=unvanished, last=unvanished]{u_1}{}{z_1}{u_2}{}{z_2}
      \picLastDiagram[mode=triangle, position=bottom]{u_2}{}{z_2}{x} \notag\\
    \leq \underbrace{\picTriangle[left]{}{o}{}}_{\leq T_p}~
      \sup_{x,\, z}\sum_{y} \picWeightedMiddleDC{o}{x}{}{}{y}{y+z}
      \underbrace{\picTriangle[right]{o}{}{}}_{\leq T_p}.
    \label{eq:example-second_lace_383}
  \end{gather}
Since
\begin{align}\lbeq{middlesection2}
\sup_{x,z}\sum_{u,v,y} \picWeightedMiddleDC{o}{x}{u}{v}{y}{y+z}
&\le\underbrace{\sup_x\picOpenTriangle{o}{x}{}{}}_{= T_p}~~
 \sup_{u,u'} \picOpenTriangle[gapposition=none]{u}{u'}{}{}~~
 \underbrace{\sup_{v,v',z}\sum_y\picTwoLines[bottomlinegap=true]
 {v}{v'}{y}{y+z}}_{\le r~(\because\refeq{boundbyr})},
\end{align}
and
\begin{align}
\picOpenTriangle[gapposition=none]{u}{u'}{}{}
&=\sum_y(pD*G_p^{*2})(y-u)\,G_p(u'-y)\nn\\
&\stackrel{\mathclap{\eqref{eq:bound-at_least_one_line-perc}}}\le
 \underbrace{(pD*G_p^{*2})(u'-u)}_{\le r~(\because\refeq{boundbyr})}
 +\underbrace{\big((pD)^{*2}*G_p^{*3}\big)(u'-u)}_{\le T_p},
\end{align}
we can bound \eqref{eq:example-second_lace_383} above by 
$T_pT_p(r+T_p)rT_p$.

Applying the same analysis to the other diagrams, we obtain
\begin{align}
\LaceCoefficientFT{2}(0)&\le\bigg(1+\frac{B_p}2+T_p\bigg)\bigg(
 \underbrace{T_p+\frac12T_pB_p+T_p^2+r+\frac12rB_p+T_pr+rT_p+T_p(r+T_p)}_{=\rho}
 \bigg)\nn\\
&\quad\times r\bigg(1+\frac{B_p}2+T_p\bigg),
\end{align}
as required.
\QED

\Proof{Sketch proof of \eqref{eq:piperc-diffbd-even} for $m=1$.}
The bound on $\abs{\Laplacian{k} \LaceCoefficientFT{2}(0)}$ is obtained by 
multiplying $1 - \cos k\cdot x$ to both sides of 
\eqref{eq:perc-diagrammatic-estimate-2} and summing the resulting expression 
over $x$.  During the course, we split $1 - \cos k\cdot x$ by using the 
telescopic inequality \eqref{eq:telescope}.  For example, the combination of 
the third diagrams in \refeq{first-diagrams}, \refeq{middle-diagrams} and 
\refeq{last-diagrams} is bounded by
\begin{align}\lbeq{Deltapi2-333}
4\sum_{\substack{y_1,y_2,y_3,y_4,\\ z_1,z_2}}
 \picFirstDiagram[triangle]{o}{}{y_1}{z_1}\hskip-1pc
 \picMiddleDiagramPB[triangle]{y_1}{y_2}{}{z_1}{z_2}{y_3}\hskip-1pc
 \picLastDiagram[mode=triangle, position=bottom]{z_2}{}{y_3}{y_4}~
 \sum_{j=1}^4\big(1-\cos k\cdot(y_j-y_{j-1})\big),
\end{align}
where $y_0=o$, while the combination of the third diagram in 
\refeq{first-diagrams}, the last diagram in \refeq{middle-diagrams} and 
the third diagram in \refeq{last-diagrams} is bounded by
\begin{align}\lbeq{Deltapi2-383}
5\sum_{\substack{y_1,y_2,y_3,y_4,y_5,\\ z_1,z_2}}
 \picFirstDiagram[triangle]{o}{}{y_1}{z_1}\hskip-1pc
 \picMiddleDiagramDC[first=unvanished, last=unvanished]{y_1}{y_2}{z_1}{z_2}
 {y_3}{y_4}\hskip-1pc
 \picLastDiagram[mode=triangle, position=bottom]{z_2}{}{y_4}{y_5}~
 \sum_{j=1}^5\big(1-\cos k\cdot (y_j-y_{j-1})\big).
\end{align}
Also, for the other 70 combinations of diagrams in \refeq{first-diagrams}, 
\refeq{middle-diagrams} and \refeq{last-diagrams}, the number of 
intervals $y_j-y_{j-1}$ is at most 5.  We use this fact to uniformly bound the 
multiplicative constant in the telescopic inequality \refeq{telescope} by 5.

Now it remains to bound each combination in terms of basic diagrams.  
For example, the contribution from $1-\cos k\cdot (y_3-y_2)$ in 
\refeq{Deltapi2-333} is bounded as
\begin{align}\lbeq{Deltapi2-333-32}
&\sum_{\substack{y_1,y_2,y_3,y_4,\\ z_1,z_2}}
 \picFirstDiagram[triangle]{o}{}{y_1}{z_1}\hskip-1pc
 \picMiddleDiagramPB[triangle]{y_1}{y_2}{}{z_1}{z_2}{y_3}\hskip-1pc
 \picLastDiagram[mode=triangle, position=bottom]{z_2}{}{y_3}{y_4}~
 \big(1-\cos k\cdot(y_3-y_2)\big)\nn\\
&\le\underbrace{\picTriangle[left]{}{o}{}}_{\le T_p}~
 \underbrace{\sup_{y_1}\picInsMiddlePB{o}{y_1}{}}_{\leq T_p^2}~
 \underbrace{\sup_{y_2,z_2,v}\sum_{y_3}\picTwoLines[bottomlinegap=true]{z_2}{y_2}
 {y_3+v}{y_3}\big(1-\cos k\cdot (y_3-y_2)\big)}_{\le \hat V_p^2(k)}~
 \underbrace{\picTriangle[right]{o}{}{}}_{\le T_p}.
\end{align}
while the contribution from $1-\cos k\cdot y_1$ in \refeq{Deltapi2-383} 
is bounded as
\begin{align}\lbeq{Deltapi2-383-10}
&\sum_{\substack{y_1,y_2,y_3,y_4,y_5,\\ z_1,z_2}}
 \picFirstDiagram[triangle]{o}{}{y_1}{z_1}\hskip-1pc
 \picMiddleDiagramDC[first=unvanished, last=unvanished]{y_1}{y_2}{z_1}{z_2}
 {y_3}{y_4}\hskip-1pc
 \picLastDiagram[mode=triangle, position=bottom]{z_2}{}{y_4}{y_5}~
 (1-\cos k\cdot y_1)\nn\\
&\le\underbrace{\sup_v\sum_{y_1}
 \picOpenBubble[mode=rightopen, toplinegap=true]{o}{y_1+v}{y_1}(1-\cos k\cdot
 y_1)}_{=\hat V_p^1(k)}~
 \underbrace{\sup_{y_2,z_1}\picOpenTriangle[mode=rightopen]{}{}{z_1}{y_2}}_{=T_p}~
 \underbrace{\sup_{y_4,w} \picInsRightMiddleDC{w}{y_4}}_{\le(r+T_p)T_p}~
 \underbrace{\picTriangle[right]{o}{}{}}_{\le T_p}.
\end{align}

The other combinations can be bounded similarly.  We note that each bound 
uses one of the diagrams $\hat V_p^0(k)$, $\hat V_p^1(k)$ and $\hat V_p^2(k)$ 
($\hat V_p^3(k)$ is used only in the bounds on 
$|\hat\Delta_k\hat\pi_p^{\sss(n)}|$ for $n\ge3$).  Which one is used depends 
on which pair of two-point functions is multiplied by $1-\cos k\cdot(y_j-y_{j-1})$: $(pD*G_p)(y_j-y_{j-1})^2$, 
$(pD*G_p)(y_j-y_{j-1})G_p(y_j-y_{j-1}-v)$ for some $v$, or 
$(pD*G_p)(y_j-y_{j-1})(pD*G_p)(y_j-y_{j-1}-v)$ for some $v$, respectively.  
We complete the sketch proof of \eqref{eq:piperc-diffbd-even} for $m=1$.
\QED

\subsubsection{Bounds on $\LaceCoefficientFT{n\geq 3}(0)$ and $\abs{\Laplacian{k}\LaceCoefficientFT{n\geq 3}(0)}$}
\label{sec:estimation-perc-lace-diff}
Recall the definition of 
$\LaceCoefficient{n}(x)$ for $n\ge3$:
\begin{multline}
\LaceCoefficient{n}(x)= \sum_{b_1,\, \dots,\, b_n}\prod_{i=1}^{n} pD(b_i)\,
 \Expectation[0]\Bigg[\ind{o\DoublyConnected \underline{b_1}}\,
 \Expectation[1]\bigg[\indn{E(\overline{b_1}, \underline{b_2};
 \ClusterWithout[0]{b_1}{o})}\cdots\\
\times\Expectation[n-1]\Big[\indn{E(\overline{b_{n-1}}, \underline{b_n};
 \ClusterWithout[n-2]{b_{n-1}}{\overline{b_{n-2}}})}\,
 \mP_p^{\sss n}\Big(E\big(\overline{b_n}, x; \ClusterWithout[n-1]
 {b_n}{\overline{b_{n-1}}}\big)\Big)\Big]\cdots\bigg]\Bigg].
\end{multline}
Using \eqref{eq:bound-last-sheet} first (cf., \refeq{last-diagrams} 
and \refeq{perc-lace-coefficient-2}), then using \refeq{bMdef} for 
$n-1$ times, and finally using \refeq{first-diagrams}, we obtain the fixed-$x$ 
bound (cf., \refeq{perc-diagrammatic-estimate-2})
\begin{gather}\lbeq{perc-diagrammatic-estimate-n}
\LaceCoefficient{n}(x)\le\sum_{\substack{u_1,\dots,u_n,\\ z_1,\dots,z_n}}
 \bL{z_1}{u_1}\times\bM{z_1}{u_1}{u_2}{z_2}\times\cdots\times\bM{z_{n-1}}
 {u_{n-1}}{u_n}{z_n}\times\bR{z_n}{u_n}x.
\end{gather}

\Proof{Sketch proof of \refeq{piperc-sumbd} for $n\ge3$.}
We can follow the same line of the proof of \refeq{piperc-sumbd} for $n=2$.  
The only difference is the size of middle section (cf., \refeq{middlesection1} 
and \refeq{middlesection2}), 
and it gives rise to the factor $\rho^{n-1}$.
\QED

\Proof{Sketch proof of \eqref{eq:piperc-diffbd-even} for $m\ge2$ and 
\eqref{eq:piperc-diffbd-odd}.}
As is done in the proof of \eqref{eq:piperc-diffbd-even} for $m=1$ in 
Section~\ref{sss:pi2}, we uniformly bound the multiplicative constant in the 
telescopic inequality \refeq{telescope} by the maximum number of intervals, 
which is $2n+1$ for $|\hat\Delta_k\hat\pi_p^{\sss(n)}(0)|$.  The remaining 
task is almost the same as the previous case, except for the following two:
\begin{enumerate}[(i)]
\item
Use the telescopic inequality \refeq{telescope} along the upper sequence of 
line segments for even $n$ (see \refeq{Deltapi2-333}--\refeq{Deltapi2-383}) 
or along the lower sequence for odd $n$ (see below).
\item
Use the basic diagram $\hat V_p^3(k)$ to bound certain diagrams to which 
$1-\cos k\cdot(\cdots)$ is assigned in a peculiar way.
\end{enumerate}
For example, 
the contribution to $|\hat\Delta_k\hat\pi_p^{\sss(3)}(0)|$ from 
the combination of the third diagrams in \eqref{eq:first-diagrams}, 
\eqref{eq:middle-diagrams} and \eqref{eq:last-diagrams} is bounded as
\begin{align}
7\sum_{\substack{y_1,y_2,y_3,\\ z_1,\dots,z_6}}
&\picFirstDiagram[triangle]{o}{z_1}{y_1}{z_2}\hskip-1pc
 \picMiddleDiagramPB[triangle]{y_1}{}{}{z_2}{z_3}{y_2}\hskip-1pc
 \picReversedMiddleDiagramPB[triangle]{z_3}{z_4}{}{y_2}{y_3}{z_5}\hskip-1pc
 \picLastDiagram[mode=triangle]{y_3}{z_5}{}{z_6}\nn\\
&\times\sum_{j=1}^6\big(1-\cos k\cdot(z_j-z_{j-1})\big),
\end{align}
where $z_0=o$.  
Then, the contribution to the sum from $1-\cos k\cdot(z_3-z_2)$ is bounded by
\begin{align}
\underbrace{\picTriangle[left]{}{o}{}}_{\leq T_p}~
 \underbrace{\sup_{y_1,z_1,u}\sum_{z_2,z_3,z_4}\picWeightedMiddlePB{z_1}{y_1}
  {z_2}{z_3}{z_4}{z_4+u}\big(1-\cos k\cdot(z_3-z_2)\big)}_{=\hat V_p^3(k)}~
 \underbrace{\sup_v\picInsReversedMiddlePB{o}v}_{\leq T_p^2}~
 \underbrace{\picTriangle[right]{o}{}{}}_{\leq T_p}.
\end{align}

The other combinations can be bounded similarly, and we refrain from showing 
tedious computations.
\QED

\subsection{Diagrammatic bounds on the bootstrapping functions}
Let
\begin{align}
\hat\Pi_p\even(k)=\sum_{m=0}^\infty\hat\pi_p^{\sss(2m)}(k),&&
\hat\Pi_p\odd(k)=\sum_{m=0}^\infty\hat\pi_p^{\sss(2m+1)}(k).
\end{align}
Suppose that $\rho\equiv T_p (2r + T_p) + (r + T_p) (1 + B_p / 2 + T_p) < 1$.  
Then, by Lemma~\ref{lmm:DBPERC}, we obtain
\begin{gather}
0\le\hat\Pi_p\even(0)\le\frac{B_p}2+\frac{(1+B_p/2+T_p)^2r\rho}{1-\rho^2},
 \lbeq{piperc-sumbd-even}\\
0\le\hat\Pi_p\even(0)\le\frac{(1+B_p/2+T_p)^2r}{1-\rho^2},
 \lbeq{piperc-sumbd-odd}\\
\sup_k\frac{|\hat\Delta_k\hat\Pi_p\even(0)|}{1-\hat D(k)}\le\sum_{j=0}^3
 \varphi_j\even\bigg\|\frac{\hat{V}_p^j}{1-\hat D}\bigg\|_\infty,
\end{gather}
where
\begin{align}
\varphi_0\even&=\frac12+\bigg(1+\frac{B_p}2+T_p\bigg)\frac{10r\rho}
 {(1-\rho^2)^2}+\bigg(1+\frac{B_p}2+T_p\bigg)^2\frac{r^2(5+12\rho^2)
 +3T_pr\rho^2 (3+\rho)}{(1-\rho^2)^3},\lbeq{vphi0even}\\
\varphi_1\even&=\bigg(1+\frac{B_p}2+T_p\bigg)\frac{5T_p\rho}{(1-\rho^2)^2}
 +\bigg(1+\frac{B_p}2+T_p\bigg)^2\frac{\rho(5+3\rho^2)+T_p(r+T_p)(5+12\rho^2)}
 {(1-\rho^2)^3},\\
\varphi_2\even&=\bigg(1+\frac{B_p}2+T_p\bigg)\frac{5T_p\rho}{(1-\rho^2)^2}
 +\bigg(1+\frac{B_p}2+T_p\bigg)^2\frac{\rho(5+3\rho^2)+T_pr(5+12\rho^2)}
 {(1-\rho^2)^3},\\
\varphi_3\even&=\bigg(1+\frac{B_p}2+T_p\bigg)^2\frac{3\rho^2 (3+\rho)}
 {(1-\rho^2)^3},
\end{align}
and
\begin{align}\label{eq:perc-sum-diff-odd}
\sup_k\frac{|\hat\Delta_k\hat\Pi_p\odd(0)|}{1-\hat D(k)}\le\sum_{j=0}^3
 \varphi_j\odd\bigg\|\frac{\hat{V}_p^j}{1-\hat D}\bigg\|_\infty,
\end{align}
where
\begin{align}
\varphi_0\odd&=1+2(B_p+r)+\frac34B_p(B_p+2r)+3T_pr+\bigg(1+\frac{B_p}2
 +T_p\bigg)\frac{14r\rho^2}{(1-\rho^2)^2}\nn\\
&\quad+\bigg(1+\frac{B_p}2+T_p\bigg)^2\frac{r\rho(2r+T_p)(7+\rho^2)}
 {(1-\rho^2)^3},\\
\varphi_1\odd&=\bigg(1+\frac{B_p}2+T_p\bigg)^2\frac{\rho(\rho+2T_pr+2T_p^2)
 (7+\rho^2)}{(1-\rho^2)^3},\\
\varphi_2\odd&=T_p(8+6B_p+9T_p)+\bigg(1+\frac{B_p}2+T_p\bigg)\frac{14T_p\rho^2}
 {(1-\rho^2)^2}\nn\\
&\quad+\bigg(1+\frac{B_p}2+T_p\bigg)^2\frac{\rho^2(14+3\rho^2)+2T_pr\rho
 (7+\rho^2)}{(1-\rho^2)^3},\\
\varphi_3\odd&=\bigg(1+\frac{B_p}2+T_p\bigg)^2\frac{\rho(7+\rho^2)}
 {(1-\rho^2)^3}.\lbeq{vphi3odd}
\end{align}
Applying these bounds to \eqref{eq:g1bd-byPi}, \eqref{eq:g2bd-byPi} and 
\eqref{eq:g3bd-byPi}, we obtain the following bounds on the bootstrapping 
functions $\{g_i(p)\}_{i=1}^3$.

\begin{lmm}\label{lmm:gbds-perc}
Suppose $\rho<1$ and that $L_p,B_p,T_p,\|\hat V_p^j/(1-\hat D)\|_\infty$, 
$j=0,1,2,3$, are so small that the inequality \refeq{suffcond-perc} holds.  
Then, we have
\begin{align}
g_1(p)&\le\bigg(1-\frac{(1+B_p/2+T_p)^2r}{1-\rho^2}\bigg)^{-1},
 \label{eq:g1bd-byTVrho}\\
g_2(p)&\le\Bigg(1-\bigg(1-\frac{B_p}2-\frac{(1+B_p/2+T_p)^2r}{1-\rho}\bigg)^{-1}
 \sum_{j=0}^3\varphi_j\odd\bigg\|\frac{\hat V_p^j}{1-\hat D}\bigg\|_\infty
 \Bigg)^{-1},\label{eq:g2bd-byTVrho}\\
g_3(p)&\le\max\bigg\{g_2(p)\bigg(1-\frac{B_p}2-\frac{(1+B_p/2+T_p)^2r}{1-\rho}
 \bigg)^{-1},\,1\bigg\}^3\nn\\
&\quad\times p^2\bigg(1+B_p+\frac{2(1+B_p/2+T_p)^2r}{1-\rho}+2\sum_{j=0}^3
 (\varphi_j\even+\varphi_j\odd)\bigg\|\frac{\hat V_p^j}{1-\hat D}\bigg\|_\infty
 \bigg)^2.\label{eq:g3bd-byTVrho}
\end{align}
\end{lmm}

\Proof{Proof.}
The bound on $g_1(p)$ is easy; since 
$\hat\Pi_p(0)=\hat\Pi_p\even(0)-\hat\Pi_p\odd(0)$, we obtain
\begin{align}
g_1(p)\stackrel{\refeq{g1bd-byPi}}\le\big(1-\hat\Pi_p\odd(0)\big)^{-1}
  \stackrel{\refeq{piperc-sumbd-odd}}{\leq}\bigg(1-\frac{(1+B_p/2+T_p)^2
 r}{1-\rho^2}\bigg)^{-1}.
\end{align}
Also, since $-\hat\Delta_k\hat\Pi_p(0)=|\hat\Delta_k\hat\Pi_p\even(0)|
-|\hat\Delta_k\hat\Pi_p\odd(0)|$ 
and $\hat I(k)\equiv1+\hat\Pi_p(k)\ge1-\hat\Pi_p\even(0)-\hat\Pi_p\odd(0)$ 
($>0$ as long as the inequality \refeq{suffcond-perc} holds), we obtain
\begin{align}\lbeq{g2bd-prebyBWr}
g_2(p)&\stackrel{\refeq{g2bd-byPi}}\le\sup_k\bigg(1-\frac1{1-\hat\Pi_p\even(0)
 -\hat\Pi_p\odd(0)}\frac{|\hat\Delta_k\hat\Pi_p\odd(0)|}{1-\hat D(k)}\bigg)^{-1}
 \nn\\
&~\:\le\Bigg(1-\bigg(1-\frac{B_p}2-\frac{(1+B_p/2+T_p)^2r}{1-\rho}\bigg)^{-1}
 \sum_{j=0}^3\varphi_j\odd\bigg\|\frac{\hat V_p^j}{1-\hat D}\bigg\|_\infty
 \Bigg)^{-1}.
\end{align}
For $g_3(p)$, since 
$\hat G_p(k)=\hat I_p(k)\hat A_p(k)\equiv \hat{I}_p(k)/(1-\hat J_p(k))$ for 
percolation and $|\hat G_p(k)|\le g_2(p)\hat S_1(k)\equiv g_2(p)/(1-\hat D(k))$, 
we obtain
\begin{align}\lbeq{g3bd-prebyBWr}
&g_3(p)\nn\\
&\stackrel{\refeq{g3bd-byPi}}\le\sup_{k,l}\frac{1-\hat D(k)}{\hat U(k,l)}
 \Bigg(\frac{\hat S_1(l+k)+\hat S_1(l-k)}2\hat S_1(l)\bigg(\frac{g_2(p)}{1
 -\hat\Pi_p\even(0)-\hat\Pi_p\even(0)}\bigg)^2\frac{|\hat\Delta_k\hat J_p(l)|}
 {1-\hat D(k)}\nn\\
&\hskip4pc+4\hat S_1(l+k)\hat S_1(l-k)\bigg(\frac{g_2(p)}{1-\hat\Pi_p\even(0)
 -\hat\Pi_p\even(0)}\bigg)^3\frac{-\hat\Delta_l\widehat{|J_p|}(0)}{1-\hat D(l)}
 \frac{-\hat\Delta_k\widehat{|J_p|}(0)}{1-\hat D(k)}\Bigg)\nn\\
&\stackrel{\refeq{U-def}}\le\max\bigg\{\frac{g_2(p)}{1-\hat\Pi_p\even(0)-\hat\Pi_p
 \even(0)},\,1\bigg\}^3\max\bigg\{\sup_{k,l}\frac{|\hat\Delta_k\hat J_p(l)|}
 {1-\hat D(k)},~\bigg(\sup_k\frac{-\hat
  \Delta_k\widehat{|J_p|}(0)}{1-\hat D(k)}\bigg)^2~\bigg\}.
\end{align}
Since $J_p=pD+\Pi_p*pD$ for percolation, we have
\begin{align}
\frac{|\hat\Delta_k\hat J_p(l)|}{1-\hat D(k)}&=\frac{p}{1-\hat D(k)}\bigg|\sum_x
 (1-\cos k\cdot x)e^{il\cdot x}\big(D(x)+(\Pi_p*D)(x)\big)\bigg|\nn\\
&\le p\bigg(\sum_x\frac{1-\cos k\cdot x}{1-\hat D(k)}D(x)+\sum_{x,y}\frac{1
 -\cos k\cdot x}{1-\hat D(k)}|\Pi_p(y)|D(x-y)\bigg)\nn\\
&\le p\bigg(1+\sum_{x,y}\frac{1-\cos k\cdot x}{1-\hat D(k)}\big(\Pi_p\even(y)
 +\Pi_p\odd(y)\big)D(x-y)\bigg),
\end{align}
which is larger than 1, since $p\ge1$.
By the telescopic inequality \eqref{eq:telescope}, the sum over $x,y$ is 
bounded as 
\begin{align}
&\sum_{x,y}\frac{1-\cos k\cdot x}{1-\hat D(k)}\big(\Pi_p\even(y)+\Pi_p\odd(y)
 \big)D(x-y)\nn\\
&\le2\sum_y\frac{(1-\cos k\cdot y)(\Pi_p\even(y)+\Pi_p\odd(y))}{1-\hat D(k)}
 \sum_xD(x-y)\nn\\
&\quad+2\sum_y\big(\Pi_p\even(y)+\Pi_p\odd(y)\big)\sum_x\frac{(1-\cos k\cdot
 (x-y))D(x-y)}{1-\hat D(k)}\nn\\
&\le2\bigg(\frac{|\hat\Delta_k\hat\Pi_p\even(0)|}{1-\hat D(k)}+\frac{|\hat
 \Delta_k\hat\Pi_p\odd(0)|}{1-\hat D(k)}\bigg)+2\big(\hat\Pi_p\even(0)+\hat
 \Pi_p\odd(0)\big).
\end{align}
As a result,
\begin{align}
\frac{|\hat\Delta_k\hat J_p(l)|}{1-\hat D(k)}\le p\Bigg(1+2\big(\hat\Pi_p\even
 (0)+\hat\Pi_p\odd(0)\big)+2\bigg(\frac{|\hat\Delta_k\hat\Pi_p\even(0)|}{1-\hat
 D(k)}+\frac{|\hat\Delta_k\hat\Pi_p\odd(0)|}{1-\hat D(k)}\bigg)\Bigg).
\end{align}
It is not difficult to check if 
$-\hat\Delta_k\widehat{|J_p|}(0)/(1-\hat D(k))$, which is nonnegative, obeys
the same bound.  Therefore, we obtain
\begin{align}
g_3(p)&\le\max\bigg\{g_2(p)\bigg(1-\frac{B_p}2-\frac{(1+B_p/2+T_p)^2r}{1-\rho}
 \bigg)^{-1},\,1\bigg\}^3\nn\\
&\quad\times p^2\bigg(1+B_p+\frac{2(1+B_p/2+T_p)^2r}{1-\rho}+2\sum_{j=0}^3
 (\varphi_j\even+\varphi_j\odd)\bigg\|\frac{\hat V_p^j}{1-\hat D}\bigg\|_\infty
 \bigg)^2,
\end{align}
as required.
\QED

\subsection{Bounds on diagrams in terms of random-walk quantities}\label{ss:TVOHbd}
In this subsection, we evaluate the diagrams 
for $p\in[1,\pc)$ and complete the proof of 
Propositions~\ref{prp:f-initial}--\ref{prp:f-bootstrapping} for percolation.  

First, we evaluate the diagrams for $p\in(1,\pc)$ under the bootstrapping 
assumptions.

\begin{lmm}\label{lmm:TpVpOpHpbd-perc}
Let $d\ge7$ and $p\in(1,\pc)$ and suppose that $g_i(p)\le K_i$, $i=1,2,3$,
for some constants $\{K_i\}_{i=1}^3$.  Then, we have
\begin{align}
T_p&\le K_1^2K_2^3\vep_3,\lbeq{Tpbd-perc}\\
\bigg\|\frac{\hat V_p^0}{1-\hat D}\bigg\|_\infty&\le K_1^2K_2\vep_1
 +5K_1^2K_2K_3\vep_3,\lbeq{Opbd-perc}\\
\bigg\|\frac{\hat V_p^1}{1-\hat D}\bigg\|_\infty&\le K_1+\frac{5K_1^2}{2^d}
 +5K_1^3K_2\vep_1+6K_1^3K_2^2\vep_2+20K_1^2K_2K_3\vep_3,\lbeq{tildeVpbd}\\
\bigg\|\frac{\hat V_p^2}{1-\hat D}\bigg\|_\infty&\le\frac{K_1^2}{2^d}
 +K_1^3K_2\vep_1+2K_1^3K_2^2\vep_2+10K_1^2K_2K_3\vep_3,\lbeq{Vpbd-perc}\\
\bigg\|\frac{\hat V_p^3}{1-\hat D}\bigg\|_\infty&\le5K_1^5K_2^7K_3(1+3\vep_1
 +2\vep_2+\vep_3)\vep_3^2.\lbeq{Hpbd-perc}
\end{align}
\end{lmm}

\Proof{Proof.}
The inequality \refeq{Tpbd-perc} has already been explained in 
\refeq{Lpbd-gen}--\refeq{BpTpbd-gen}.

To prove \refeq{Opbd-perc}, we first use the trivial inequality 
$G_p\le\delta+pD*G_p$ to obtain 
\begin{align}\lbeq{V0bd-1st}
\frac{\hat V_p^0(k)}{1-\hat D(k)}
&\le\sum_{x\sim o}\frac{1-\cos k\cdot x}{1-\hat D(k)}pD(x)\,(pD*G_p)(x)\nn\\
&\quad+\sum_x\frac{1-\cos k\cdot x}{1-\hat D(k)}\big((pD)^{*2}*G_p\big)(x)
 \,(pD*G_p)(x).
\end{align}
By symmetry (cf., \refeq{symm-improve}), the first term is bounded by
\begin{align}
p\,\sup_{x\sim o}(pD*G_p)(x)=\underbrace{\frac{p}{2^d}\sum_{x\sim o}(pD*G_p)
 (x)}_{=((pD)^{*2}*G_p)(o)}=L_p\stackrel{\refeq{LpBpB'pbd-saw}}\le K_1^2K_2
 \vep_1.
\end{align}
For the second term, we use the Fourier representation to obtain
\begin{align}\lbeq{V0bd-2nd}
&\sum_x\big((pD)^{*2}*G_p\big)(x)\frac{(pD*G_p)(x)(1-\cos k\cdot x)}{1-\hat
 D(k)}\nn\\
&=\int_{\Td}\frac{\mathrm{d}^dl}{(2\pi)^d}\big(p\hat D(l)\big)^2\hat G_p(l)
 \frac{-\hat\Delta_k(p\hat D(l)\hat G_p(l))}{1-\hat D(k)}\nn\\
&\le K_1^2K_2K_3\int_{\Td}\frac{\mathrm{d}^dl}{(2\pi)^d}\hat D(l)^2\hat S_1(l)
 \frac{\hat U(k,l)}{1-\hat D(k)}.
\end{align}
Recall the definition \refeq{U-def} of $\hat U(k,l)$, in which we have three 
terms: $\frac12\hat S_1(l+k)\hat S_1(l)$, $\frac12\hat S_1(l-k)\hat S_1(l)$ 
and $4\hat S_1(l+k)\hat S_1(l-k)$.  By inversion, we have, e.g.,
\begin{align}\lbeq{e.g.U}
&\int_{\Td}\frac{\mathrm{d}^dl}{(2\pi)^d}\hat D(l)^2\hat S_1(l)\hat S_1(l+k)
 \hat S_1(l-k)\nn\\
&=\sum_{x,y}S_1(x)S_1(y)e^{ik\cdot (x-y)}\underbrace{\int_{\Td}
 \frac{\mathrm{d}^dl}{(2\pi)^d}\hat D(l)^2\hat S_1(l)e^{il(x+y)}}_{=
 (D^{*2}*S_1)(x+y)}\nn\\
&\le(D^{*2}*S_1^{*3})(o)=\vep_3.
\end{align}
It is not difficult to check if the other combinations obey the same bound.  
This completes the proof of \refeq{Opbd-perc}.

Next, we prove \refeq{Vpbd-perc} before showing \refeq{tildeVpbd}.  
First, by the trivial inequality $G_p\le\delta+pD*G_p$, we obtain
\begin{align}
\frac{\hat V_p^2(k)}{1-\hat D(k)}
&\le\sup_x\sum_y\frac{pD(y)(1-\cos k\cdot y)}{1-\hat D(k)}(pD*G_p)(x-y)\nn\\
&\quad+\sup_x\sum_y\frac{((pD)^{*2}*G_p)(y)(1-\cos k\cdot y)}{1-\hat D(k)}
 (pD*G_p)(x-y).
\end{align}
Since $\|pD*G_p\|_\infty\le p/2^d+L_p$, the first term is bounded by 
$K_1(K_1/2^d+K_1^2K_2\vep_1)$.  For the second term, we use the telescopic 
inequality \refeq{telescope} to obtain
\begin{align}
&2\sum_{y,z}\frac{pD(z)(pD*G_p)(y-z)}{1-\hat D(k)}\Big((1-\cos k\cdot z)+\big(
 1-\cos k\cdot(y-z)\big)\Big)(pD*G_p)(x-y)\nn\\
&\le2\sum_z\frac{pD(z)(1-\cos k\cdot z)}{1-\hat D(k)}\big((pD)^{*2}*G_p^{*2}
 \big)(x-z)\nn\\
&\quad+2\sum_{y'}\frac{(pD*G_p)(y')(1-\cos k\cdot y')}{1-\hat D(k)}
 \big((pD)^{*2}*G_p\big)(x-y'),
\end{align}
where we have used the replacement $y'=y-z$.  Since 
$((pD)^{*2}*G_p^{*2})(x-z)\le B_p$, which is due to the Schwarz inequality, 
the first term is bounded by $2pB_p\le 2K_1^3K_2^2\vep_2$ (cf., 
\refeq{LpBpB'pbd-saw}).  On the other hand, by the Fourier representation, 
the second term is bounded by (cf., \refeq{V0bd-2nd}--\refeq{e.g.U})
\begin{align}\lbeq{FS-5factor1}
&2\int_{\Td}\frac{\mathrm{d}^dl}{(2\pi)^d}\frac{|\hat\Delta_k(p\hat D(l)\hat
 G_p(l))|}{1-\hat D(k)}(p\hat D(l))^2\hat G_p(l)\nn\\
&\le2 K_1^2K_2K_3\int_{\Td}\frac{\mathrm{d}^dl}{(2\pi)^d}\frac{\hat U(k,l)}
 {1-\hat D(k)}\hat D(l)^2\hat S_1(l)\le10 K_1^2K_2K_3\vep_3.
\end{align}
This completes the proof of \refeq{Vpbd-perc}.

Similarly, by using the trivial inequality $G_p\le\delta+pD*G_p$ and the 
telescopic inequality \refeq{telescope}, we have
\begin{align}
\frac{\hat V_p^1(k)}{1-\hat D(k)}
&\le\sup_x\sum_y\frac{pD(y)(1-\cos k\cdot y)}{1-\hat D(k)}\underbrace{G_p
 (x-y)}_{\le\|G_p\|_\infty}\nn\\
&\quad+2\sup_x\sum_y\frac{pD(y)(1-\cos k\cdot y)}{1-\hat D(k)}\underbrace{(pD
 *G_p^{*2})(x-y)}_{\le r}\nn\\
&\quad+2\underbrace{\sup_x\sum_y\frac{(pD*G_p)(y)(1-\cos k\cdot y)}{1-\hat
 D(k)}(pD*G_p)(x-y)}_{=\|\hat V_p^2/(1-\hat D)\|_\infty}.
\end{align}
By $\|G_p\|_\infty\le1+p/2^d+L_p\le1+K_1/2^d+K_1^2K_2\vep_1$ and using 
\refeq{rpbd} and \refeq{Vpbd-perc}, we obtain \refeq{tildeVpbd}.

It remains to show the bound \refeq{Hpbd-perc} of order $\vep_3^2$.  
This is an improvement from a naive bound of order $\vep_3$, and is a 
result of repeated applications of the H\"older and Schwarz inequalities, 
as explained now.  First, we recall the definition of $\hat V_p^3(k)$:
\begin{align}
\hat V_p^3(k)=\sup_{x,y}\sum_{\{v_j\}_{j=1}^5}&G_p(v_1)\,(pD*G_p)(v_2-v_1)
 \big(1-\cos k\cdot(v_2-v_1)\big)\nn\\
&\times(pD*G_p)(v_3-v_2)\,(pD*G_p)(v_4-v_1)\,(pD*G_p)(v_4-v_2)\nn\\
&\times(pD*G_p)(v_5-v_4)\,(pD*G_p)(x-v_5)\,G_p(y+v_3-v_5).
\end{align}
By using the Fourier representation and then the assumptions $g_j(p)\le K_j$ 
for $j=1,2,3$, the above sum over $\{v_j\}_{j=1}^5$ is bounded above as
\begin{align}\lbeq{barHprebd}
&\int\prod_{j=1}^3\frac{\mathrm{d}^dl_j}{(2\pi)^d}\hat G_p(l_1)\Big(-\hat
 \Delta_k\big(p\widehat{D*G_p}(l_2)\big)\Big)p\widehat{D*G_p}(l_3)\,p
 \widehat{D*G_p}(l_1-l_2)\nn\\
&\hskip5pc\times p\widehat{D*G_p}(l_3-l_2)\,p\widehat{D*G_p}(l_1-l_3)\,
 p\widehat{D*G_p}(l_1)\,\hat G_p(l_3)\,e^{-il_1\cdot x+il_3\cdot y}\nn\\
&\le p^5\int\prod_{j=1}^3\frac{\mathrm{d}^dl_j}{(2\pi)^d}|\hat G_p(l_1)|\,
 \big|\hat\Delta_k\big(p\hat D(l_2)\hat G_p(l_2)\big)\big|\,|\hat D(l_3)
 \hat G_p(l_3)|\,|\hat D(l_1-l_2)\hat G_p(l_1-l_2)|\nn\\
&\hskip5pc\times|\hat D(l_3-l_2)\hat G_p(l_3-l_2)|\,|\hat D(l_1-l_3)
 \hat G_p(l_1-l_3)|\,|\hat D(l_1)\hat G_p(l_1)|\,|\hat G_p(l_3)|\nn\\
&\le K_1^5K_2^7K_3\int\prod_{j=1}^3\frac{\mathrm{d}^dl_j}{(2\pi)^d}\hat
 S_1(l_1)\,\hat U(k,l_2)\,|\hat D(l_3)|\hat S_1(l_3)\,|\hat D(l_1-l_2)|\hat
 S_1(l_1-l_2)\nn\\
&\hskip5pc\times|\hat D(l_3-l_2)|\hat S_1(l_3-l_2)\,|\hat D(l_1-l_3)|\hat
 S_1(l_1-l_3)\,|\hat D(l_1)|\hat S_1(l_1)\,\hat S_1(l_3),
\end{align}
where we have used the abbreviation $\int=\iiint_{(\Td)^3}$ 
and the fact that $\hat S_1\ge0$.

To investigate the above integral, we introduce the notation, such as 
$\hat S_{1-2}=\hat S_1(l_1-l_2)$ and $\hat S_{2+}=\hat S_1(l_2+k)$ (n.b. 
the new subscripts are not the values of $p$).  By repeated applications of 
the H\"older and Schwarz inequalities and periodicity, 
the contribution from, e.g., $\hat S_{2+}\hat S_{2-}$ in $\hat U(k,l_2)$ is 
bounded as
\begin{align}
&\int\prod_{j=1}^3\frac{\mathrm{d}^dl_j}{(2\pi)^d}\hat S_1\,\hat S_{2+}\hat
 S_{2-}\,|\hat D_3|\hat S_3\,|\hat D_{1-2}|\hat S_{1-2}\,|\hat D_{3-2}|\hat
 S_{3-2}\,|\hat D_{1-3}|\hat S_{1-3}\,|\hat D_1|\hat S_1\,\hat S_3\nn\\
&\le\bigg(\int\prod_{j=1}^3\frac{\mathrm{d}^dl_j}{(2\pi)^d}\Big(\hat S_{2+}\,
 |\hat D_{3-2}|\hat S_{3-2}\,|\hat D_{1-3}|\hat S_{1-3}\Big)^3\bigg)^{1/3}\nn\\
&\quad\times\bigg(\int\prod_{j=1}^3\frac{\mathrm{d}^dl_j}{(2\pi)^d}\Big(|\hat
 D_1|\hat S_1^2\,|\hat D_3|\hat S_3^2\,\hat S_{2-}\,|\hat D_{1-2}|\hat
 S_{1-2}\Big)^{3/2}\bigg)^{2/3}\nn\\
&=\bigg(\underbrace{\int_{\Td}\frac{\mathrm{d}^dl_2}{(2\pi)^d}\hat S_{2+}^3}_{
 \equiv\bigtriangledown}\underbrace{\int_{\Td}\frac{\mathrm{d}^dl_3}{(2\pi)^d}
 |\hat D_{3-2}|^3\hat S_{3-2}^3}_{\le\vep_3}\underbrace{\int_{\Td}\frac{
 \mathrm{d}^dl_1}{(2\pi)^d}|\hat D_{1-3}|^3\hat S_{1-3}^3}_{\le\vep_3}
 \bigg)^{1/3}\nn\\
&\quad\times\bigg(\int_{\Td}\frac{\mathrm{d}^dl_1}{(2\pi)^d}|\hat D_1|^{3/2}
 \hat S_1^3\int_{\Td}\frac{\mathrm{d}^dl_3}{(2\pi)^d}|\hat D_3|^{3/2}\hat S_3^3
 \underbrace{\int_{\Td}\frac{\mathrm{d}^dl_2}{(2\pi)^d}\hat S_{2-}^{3/2}\,
 |\hat D_{1-2}|^{3/2}\hat S_{1-2}^{3/2}}_{\le\bigtriangledown^{1/2}
 \vep_3^{1/2}~(\because\text{Schwarz})}\bigg)^{2/3}\nn\\
&\le\bigtriangledown^{2/3}\vep_3\bigg(\int_{\Td}\frac{\mathrm{d}^dl_1}{(2
 \pi)^d}|\hat D_1|^{3/2}\hat S_1^3\bigg)^{4/3}\nn\\
&\le\bigtriangledown^{2/3}\vep_3\Bigg(\bigg(\int_{\Td}\frac{\mathrm{d}^dl_1}
 {(2\pi)^d}\hat D_1^2\hat S_1^3\bigg)^{3/4}\bigg(\int_{\Td}\frac{\mathrm{d}^d
 l_1}{(2\pi)^d}\hat S_1^3\bigg)^{1/4}\Bigg)^{4/3}\nn\\
&\le\bigtriangledown\,\vep_3^2.
\end{align}
By the identity $S_1=\delta+D*S_1$, we further obtain
\begin{align}\lbeq{S3bd}
\bigtriangledown=\int_{\Td}\frac{\mathrm{d}^dl}{(2\pi)^d}\hat S_{1}(l)^3
 =1+3\vep_1+2\vep_2+\vep_3.
\end{align}
The contributions from the other terms in $\hat U(k,l_2)$ obey the same 
bound.  This completes the proof of \refeq{Hpbd-perc}, hence the proof of Lemma~\ref{lmm:TpVpOpHpbd-perc}.
\QED

Next, we evaluate the diagrams at $p=1$ by using the trivial inequality 
$G_1(x)\le S_1(x)$.  Here, we do not need the bootstrapping assumptions.

\begin{lmm}\label{lmm:T1V1O1H1bd-perc}
Let $d\ge7$ and $p=1$.  Then, we have
\begin{align}
T_1&\le\vep_3,\lbeq{T1bd-perc}\\
\bigg\|\frac{\hat V_1^0}{1-\hat D}\bigg\|_\infty&\le\vep_1+5\vep_3,
 \lbeq{O1bd-perc}\\
\bigg\|\frac{\hat V_1^1}{1-\hat D}\bigg\|_\infty&\le 1+\frac5{2^d}
 +5\vep_1+6\vep_2+20\vep_3,\lbeq{tildeV1bd}\\
\bigg\|\frac{\hat V_1^2}{1-\hat D}\bigg\|_\infty&\le\frac1{2^d}+\vep_1+2\vep_2
 +10\vep_3,\lbeq{V1bd-perc}\\
\bigg\|\frac{\hat V_1^3}{1-\hat D}\bigg\|_\infty&\le5(1+3\vep_1+2\vep_2+\vep_3)
 \vep_3^2.\lbeq{H1bd-perc}
\end{align}
\end{lmm}

\Proof{Proof.}
The inequality \refeq{T1bd-perc} has already been explained in 
\refeq{L1bd-gen}--\refeq{B1T1bd-gen}.  Similarly, we can show the other 
inequalities by using the trivial inequality $G_1(x)\le S_1(x)$.  For example 
(cf., \refeq{V0bd-1st}--\refeq{V0bd-2nd}),
\begin{align}\lbeq{V01bd-1st}
\frac{\hat V_1^0(k)}{1-\hat D(k)}
&\le\sum_{x\sim o}\frac{D(x)(1-\cos k\cdot x)}{1-\hat D(k)}(D*S_1)(x)\nn\\
&\quad+\sum_x(D^{*2}*S_1)(x)\frac{(D*S_1)(x)(1-\cos k\cdot x)}{1-\hat D(k)}\nn\\
&\le\underbrace{(D^{*2}*S_1)(o)}_{=\vep_1}+\int_{\Td}\frac{\mathrm{d}^dl}{(2
 \pi)^d}\hat D(l)^2\hat S_1(l)\frac{|\hat\Delta_k(\hat D(l)\hat S_1(l))|}{1-\hat
 D(k)}.
\end{align}
Since
\begin{align}
\big|\hat\Delta_k\big(\hat D(l)\hat S_1(l)\big)\big|=\big|\hat\Delta_k\big(\hat
 S_1(l)-1\big)\big|=|\hat\Delta_k\hat S_1(l)|\le\hat U(k,l),
\end{align}
the integral in \refeq{V01bd-1st} is bounded by $5\vep_3$ (cf., \refeq{e.g.U}). 
To avoid redundancy, we refrain from showing the other inequalities.  
This completes the proof of Lemma~\ref{lmm:T1V1O1H1bd-perc}.
\QED

\Proof{Proof of Proposition \ref{prp:f-initial}}
Since $\vep_1$, $\vep_2$ and $\vep_3$ are finite for $d\ge7$ (see
Table~\ref{table:vep123} in Section~\ref{ss:bcc}) and
decreasing in $d$ (because $D^{*2n}(o)\equiv\big(\binom{2n}n2^{-2n}\big)^d$ 
on $\Ld$ is decreasing in $d$), we have
\begin{align}\lbeq{r1bd-perc}
  r = \SupNorm{D} + L_1 + B_1 \stackrel{\refeq{L1B1B'1bd-saw}}
  \leq 2^{-d} + \vep_1 + \vep_2 \le
  \begin{cases}
    0.0326 & [d=7],\\
    0.0146 & [d=8],\\
    0.0068 & [d\ge 9],
  \end{cases}
\end{align}
and, by \refeq{L1B1B'1bd-saw}, \refeq{T1bd-perc} and also \refeq{r1bd-perc},
\begin{align}
  \rho &= T_1 (2r + T_1) + (r + T_1) \left( 1 + \frac{B_1}{2} + T_1\right)\nn\\
  &\le \vep_3(2r+\vep_3)+(r+\vep_3)\left(1 + \frac{\vep_2}2 +\vep_3 \right)
  \le
  \begin{cases}
    0.0967 & [d=7],\\
    0.0279 & [d=8],\\
    0.0111 & [d\ge 9].
  \end{cases}
\end{align}
In addition, by \refeq{piperc-sumbd-even}--\refeq{vphi3odd} and 
Lemma~\ref{lmm:T1V1O1H1bd-perc}, we have
\begin{align}
\sum_{n=0}^\infty\hat\pi_1^{\sss(n)}(0)+\sup_k\sum_{n=0}^\infty\frac{-\hat
 \Delta_k\hat\pi_1^{\sss(n)}(0)}{1-\hat D(k)}  \le
  \begin{cases}
    2.7700 & [d=7],\\
    0.3623 & [d=8],\\
    0.1124 & [d\ge9],
  \end{cases}
\end{align}
which implies that the inequality \refeq{suffcond-perc} holds for all $d\ge8$ 
(but not for $d=7$).  Then, by Lemma~\ref{lmm:gbds-perc}, we obtain
\begin{align}
g_1(1)&\le\bigg(1-\frac{(1+\vep_2/2+\vep_3)^2r}{1-\rho^2}\bigg)^{-1}\le
  \begin{cases}
    1.0154 & [d=8],\\
    1.0070 & [d\ge9],
  \end{cases}\label{eq:g1(1)bd-perc}\\
g_2(1)&\le\Bigg(1-\bigg(1-\frac{\vep_2}2-\frac{(1+\vep_2/2+\vep_3)^2r}{1-\rho}
 \bigg)^{-1}\sum_{j=0}^3\varphi_j\odd\bigg\|\frac{\hat V_1^j}{1-\hat D}
 \bigg\|_\infty\Bigg)^{-1}\nn\\
 &\le
 \begin{cases}
 1.1049 & [d=8],\\
 1.0272 & [d\ge9],
 \end{cases}\label{eq:g2(1)bd-perc}\\
g_3(1)&\le \Big[\text{the bounds in \refeq{g2(1)bd-perc}}\Big]^3\times\bigg(
 1-\frac{\vep_2}2-\frac{(1+\vep_2/2+\vep_3)^2r}{1-\rho}\bigg)^{-3}\nn\\
 &\qquad\times \bigg(1+\vep_2+\frac{2(1+\vep_2/2+\vep_3)^2r}{1-\rho}+2\sum_{j
 =0}^3(\varphi_j\even+\varphi_j\odd)\bigg\|\frac{\hat V_1^j}{1-\hat D}
 \bigg\|_\infty\bigg)^2\nn\\
 &\le
 \begin{cases}
 4.2433 & [d=8],\\
 1.6673 & [d\ge 9].
 \end{cases}\lbeq{g3(1)bd-perc}
\end{align}
Proposition~\ref{prp:f-initial} for percolation holds as long as
$K_1>1.0154$, $K_2>1.1049$, $K_3>4.2433$ for $d=8$, and
$K_1>1.0070$, $K_2>1.0272$, $K_3>1.6673$ for $d\ge9$.
\QED

\Proof{Proof of Proposition \ref{prp:f-bootstrapping}}
Let
\begin{align}
  K_1 = 1.01,&&
  K_2 = 1.09,&&
  K_3 = 2.70,
\end{align}
so that Proposition~\ref{prp:f-initial} holds for $d\ge 9$.
Then, by Table~\ref{table:vep123} in Section~\ref{ss:bcc}, we obtain
\begin{equation}\lbeq{rpbd-perc}
  r\stackrel{\refeq{LpBpB'pbd-saw}}\le K_12^{-d}+K_1^2K_2\vep_1+K_1^2K_2^2\vep_2
  \le 0.0077,
\end{equation}
and, by \refeq{LpBpB'pbd-saw}, \refeq{Tpbd-perc} and \refeq{rpbd-perc},
\begin{align}
\rho&\le K_1^2 K_2^3 \varepsilon_3 \left( 2r + K_1^2 K_2^3 \varepsilon_3\right) + \left( r + K_1^2 K_2^3 \varepsilon_3\right) \left( 1 + \frac{K_1^2 K_2^2 \varepsilon_2}{2} + K_1^2 K_2^3 \varepsilon_3\right)\nn\\
&\le 0.0134.
\end{align}
In addition, by \refeq{piperc-sumbd-even}--\refeq{vphi3odd} and 
Lemma~\ref{lmm:TpVpOpHpbd-perc}, we have
\begin{align}
\sum_{n=0}^\infty\hat\pi_p^{\sss(n)}(0)+\sup_k\sum_{n=0}^\infty\frac{-\hat
 \Delta_k\hat\pi_p^{\sss(n)}(0)}{1-\hat D(k)}\le0.2151,
\end{align}
which implies that the inequality \refeq{suffcond-perc} holds.  
Then, similarly to \refeq{g1(1)bd-perc}--\refeq{g3(1)bd-perc}, we obtain
\begin{align}
  g_1(p) &\le \left( 1 - \frac{(1 + K_1^2 K_2^2\vep_2 / 2 + K_1^2 K_2^3 \vep_3)^2 r}{1 - \rho^2}\right)^{-1} \le 1.0080 < K_1,\\
  g_2(p) &\le \Bigg(1-\bigg(1-\frac{K_1^2 K_2^2 \vep_2}2-\frac{(1+K_1^2 K_2^2 \vep_2/2+K_1^2 K_2^3 \vep_3)^2r}{1-\rho}\bigg)^{-1}
 \sum_{j=0}^3\varphi_j\odd\bigg\|\frac{\hat V_p^j}{1-\hat D}\bigg\|_\infty
 \Bigg)^{-1}\nn\\
 &\le 1.0810 < K_2,\\
  g_3(p) &\le (1.081)^3 \bigg(1-\frac{K_1^2 K_2^2 \vep_2}2-\frac{(1+K_1^2 K_2^2 \vep_2/2+K_1^2 K_2^3 \vep_3)^2r}{1-\rho}\bigg)^{-3}K_1^2 \nn\\
    &\times \bigg(1+K_1^2 K_2^2 \vep_2+\frac{2(1+K_1^2 K_2^2 \vep_2/2+K_1^2 K_2^3 \vep_3)^2r}{1-\rho}+2\sum_{j=0}^3
 (\varphi_j\even+\varphi_j\odd)\bigg\|\frac{\hat V_p^j}{1-\hat D}\bigg\|_\infty
 \bigg)^2 \nn\\
    &\le 2.6606 < K_3.
\end{align}
This completes the proof of Proposition~\ref{prp:f-bootstrapping} for percolation.
\QED

\subsection{Further discussion}\label{ss:discussion-perc}
We have been able to prove convergence of the lace expansion for percolation 
on $\mL^{\!d\ge9}$ in full detail. Compared with the analysis for SAW, 
the analysis for percolation is more 
involved.  However, as compared to the NoBLE analysis on $\mZ^{\,d\ge11}$ 
\cite{fh15a,fh15b}, the current analysis is much simpler, shorter and 
more transparent.
This is due to the simple structure of the BCC lattice $\Ld$ and the choice of 
the bootstrapping functions $\{g_i(p)\}_{i=1}^3$.

To go down to the desired 7 dimensions, we must improve our analysis in 
various aspects.  Some of the key elements we can think of are almost identical 
to those for SAW already mentioned in Section~\ref{ss:discussion-saw}, with 
slight modifications as follows.
\begin{enumerate}
\item
The largest contribution comes from $|\hat\Delta_k\hat\pi_p^{\sss(0)}(0)|$ and 
$|\hat\Delta_k\hat\pi_p^{\sss(1)}(0)|$, and their common leading term is 
proportional to the diagram function $\hat V_p^0(k)$.  To improve its bound, 
we may introduce extra diagrams, such as 
$B'_p\equiv\|(pD)^{*4}*G_p^{*2}\|_\infty$ 
and $T'_p\equiv\|(pD)^{*4}*G_p^{*3}\|_\infty$, as is done for SAW, or even 
longer diagrams, such as $T^{\sss(n)}_p\equiv\|(pD)^{*2n}*G_p^{*3}\|_\infty$.  
Although its RW counterpart $(D^{*2n}*S_1^{*3})(o)$ is decreasing in $n$, 
the bound on $T^{\sss(n)}_p$ may attain the minimum at some $n_*\in\mN$, 
due to the exponentially growing factor $p^{2n}$.  So far, we have not 
investigated a result of using $T^{\sss(n_*)}_p$, since introducing such new 
diagrams increases the number of terms to deal with, which may cause extra 
complication.
\item
Similarly to the case of $\|\hat W_p/(1-\hat D)\|_\infty$ for SAW, we used the 
Schwarz inequality to bound $\|\hat V_p^j/(1-\hat D)\|_\infty$, $j=0,1,2,3$, in 
Lemmas~\ref{lmm:TpVpOpHpbd-perc}--\ref{lmm:T1V1O1H1bd-perc}.  As a 
result, the relatively large factor 5 appeared (see, e.g., \refeq{FS-5factor1}), 
as explained in the third footnote.  It would be of great help if we could do 
away with the Schwarz inequality to achieve a better bound on 
\refeq{trig-identity}.
\item
As is the case for SAW, we ignored the contributions from $\hat\Pi_p\even(0)$ 
and $|\hat\Delta_k\hat\Pi_p\even(0)|$ in 
\refeq{g1bd-byTVrho}--\refeq{g2bd-byTVrho}.  
If we include their effect into computation, then $g_1(p)$ could be much closer 
to 1 (cf. \refeq{g1bd-byPi}) and $g_2(p)$ could be even smaller than 1 (cf. 
\refeq{g2bd-byPi}), and as a result, we could achieve convergence of the lace 
expansion on $\mL^{\!d\ge7}$.  However, to make use of those even terms, 
we must also control lower bounds on $g_1(p)$ and $g_2(p)$, and to do so, 
we need nontrivial lower bounds on the lace-expansion coefficients.  
Achieving this goal without causing too much complication would be a 
challenging task.
\item
Instead of estimating $\hat G_p(k)$ by $\hat S_1(k)\equiv(1-\hat D(k))^{-1}$ 
uniformly in $k\in\Td$, which may not be efficient in the ultraviolet regime, 
we may split the estimate of $G_p(x)$ between for large $x$ and for small $x$.  
For small $x$, estimating $G_p(x)$ by $S_1(x)$ is expected to be a bit too 
crude.  Therefore, estimating the lace-expansion coefficients in the 
ultraviolet regime by naively using the BK inequality would be too primitive.  
Here, we may need incorporate the ultraviolet regularization of \cite{an84} or 
large-field analysis in the rigorous renormalization group for spin systems.
\end{enumerate}

\section*{Acknowledgements}
The work of SH is supported by the Ministry of Education, Culture, Sports,
Science and Technology through Program for Leading Graduate Schools
(Hokkaido University ``Ambitious Leader's Program'').
The work of AS is supported in part by the the JSPS KAKENHI Grant 
Number 15K13440.
We are grateful to Lung-Chi Chen, Robert Fitzner, Markus Heydenreich,
Remco van der Hofstad and Gordon Slade for continual encouragement.

\end{document}